\documentclass[aos]{imsart}

\RequirePackage[OT1]{fontenc}
\usepackage{amsmath,epsfig,amssymb,amsfonts,amsthm,verbatim,epstopdf, enumitem,hyperref,lipsum,float,multirow,bbm}
\usepackage{color,graphicx,lscape,longtable,mathabx,natbib}
\newcommand{\indep}{\;\, \rule[0em]{.03em}{.67em} \hspace{-.25em}
\rule[0em]{.65em}{.03em} \hspace{-.25em}\rule[0em]{.03em}{.67em}\;\,}

\newtheorem{Th}{\underline{\bf Theorem}}

\newtheorem{Pro}{Proposition}
\newtheorem{Lem}{\underline{\bf Lemma}}

\def\bse{\begin{eqnarray*}}
\def\ese{\end{eqnarray*}}
\def\be{\begin{eqnarray}}
\def\ee{\end{eqnarray}}
\def\bsq{\begin{equation*}}
\def\esq{\end{equation*}}
\def\bq{\begin{equation}}
\def\eq{\end{equation}}

\def\var{\hbox{var}}

\def\wh{\widehat}
\def\wt{\widetilde}
\def\eff{_{\rm eff}}

\def\n{\nonumber}

\def\sumi{\sum_{i=1}^n}
\def\sumj{\sum_{j=1}^n}
\def\trans{^{\rm T}}

\def\bb{\boldsymbol\beta}

\def\bg{\boldsymbol\gamma}

\def\0{{\bf 0}}
\def\A{{\bf A}}
\def\a{{\bf a}}

\def\g{{\bf g}}

\def\h{{\bf h}}

\def\I{{\bf I}}

\def\K{{\bf K}}

\def\T{{\bf T}}
\def\bP{{\bf P}}

\def\bS{{\bf S}}
\def\u{{\bf u}}

\def\v{{\bf v}}

\def\T{{\bf T}}

\def\X{{\bf X}}
\def\x{{\bf x}}
\def\I{{\bf I}}

\def\z{{\bf z}}

\def\blam{{\boldsymbol \lambda}}

\def\bq{\begin{equation}}
\def\eq{\end{equation}}
\def\pr{\hbox{pr}}
\def\wh{\widehat}
\def\wt{\widetilde}
\def\trans{^{\rm T}}

\def\log{\hbox{log}}

\def\squarebox#1{\hbox to #1{\hfill\vbox to #1{\vfill}}}

\def\var{\hbox{var}}

\def\bse{\begin{eqnarray*}}
\def\ese{\end{eqnarray*}}
\def\be{\begin{eqnarray}}
\def\ee{\end{eqnarray}}
\def\bsq{\begin{equation*}}
\def\esq{\end{equation*}}
\def\bq{\begin{equation}}
\def\eq{\end{equation}}
\def\pr{\hbox{pr}}
\def\wh{\widehat}
\def\wt{\widetilde}

\def\log{\hbox{log}}

\def\trans{^{\rm T}}

\def\boxit#1{\vbox{\hrule\hbox{\vrule\kern6pt\vbox{\kern6pt#1\kern6pt}\kern6pt\vrule}\hrule}}

\allowdisplaybreaks 
\begin{document}
	
	\begin{frontmatter}
\title{Efficient Estimation for Dimension Reduction with Censored Data}
\author{Ge Zhao, Yanyuan Ma and Wenbin Lu{\thanks{Ge Zhao is graduate student (E-mail: gzz13@psu.edu) and Yanyuan Ma is professor (E-mail: yzm63@psu.edu),
			Department of Statistics, Penn State University, University Park, PA,
			16802.  
			Wenbin Lu is professor (E-mail: lu@stat.ncsu.edu), 
			Department of Statistics, North Carolina State University, Raleigh, NC 27695.
			The work is supported by NSF grant DMS-1608540.
}}}     
\begin{abstract}
We propose a general index model for survival data, which generalizes
many commonly used semiparametric survival models and belongs to the
framework of dimension reduction. Using a combination 
of geometric approach in semiparametrics and martingale treatment in
survival data analysis, we devise estimation procedures that are
feasible and do not require covariate-independent censoring as assumed
in many dimension reduction methods for censored survival data. We
establish the root-$n$ consistency and asymptotic normality of the proposed
estimators and derive the most efficient estimator in this class for
the general index model. Numerical experiments are carried 
out to demonstrate the empirical performance of the proposed
estimators and an application to an AIDS data further illustrates the
usefulness of the work.
\end{abstract}


\end{frontmatter}

\section{Introduction}

Cox proportional hazards model \citep{cox1972} is probably the most
widely used semiparametric model for analyzing survival data. In the Cox
model, covariate effect is described by a single linear
combination of covariates in an exponential function and is
multiplicative in modeling the hazard function. Although this special
way of modeling the hazard function permits a convenient estimation
procedure, such as the maximum partial likelihood estimation
\citep{Cox1975}, it has its limitations. As widely studied in the
literature, there are many situations where the Cox model
may not be proper. Due to the limitations of the Cox model,
many other semiparametric survival models have been proposed in the literature, such as the
accelerated failure time model \citep{bj1979}, proportional odds model
\citep{mccullagh1980} and linear transformation model \citep{dd1988}, etc.
Despite of all these efforts, the link between the summarized covariate
effect, typically in the form of a linear combination of
covariates, and the possibly transformed event time remains to have a
predetermined form and hence can be restrictive sometimes. 

The single index feature of the above mentioned semiparametric
survival models is appealing since the covariates effect has a nice
interpretation. It also naturally achieves dimension reduction when
there is a large number of covariates. However, the specific model
form to link the covariate index to the event time may be restrictive,
and it is often difficult to check the goodness-of-fit of the specific
link function form.
To achieve a model that is flexible yet is feasible in practice, we
borrow and extend the idea of linear summary of the covariate effects,
while free up the specific functional relation between the event time and the linear
summaries. Thus, we propose the following general index model 
\be\label{eq:model}
\pr(T\le t\mid\X)=\pr(T\le t\mid \bb_0\trans\X), \quad t > 0
\ee
where $T$ is the survival time of interest, $\X$ is the $p$-dimensional baseline covariates, and $\bb_0\in{\cal R}^{p\times d}$ is the regression coefficient matrix, with  
$p>d$. Several properties of model
(\ref{eq:model}) is worth mentioning. 1) First of all, instead of a
single linear summary, we allow $d$ linear summaries described by the
$d$ columns of $\bb_0$. This increases the flexibility of how the
covariate effects are combined. We can view 
this as a generalization from single index to multi index 
covariate summary.
Imagine an extreme case when $d=p$,
this model degenerates to the restriction free case where the
dependence of $T$ on $\X$ is arbitrary. Of course, in practice, when
$d$ is large, the estimation will encounter difficulties and it is not
feasible to carry out the analysis. However conceptually this provides
a way of appreciating the flexibility of the model. In addition, we
will see that in practice, when $d$ is often smaller than $p$, this
model framework allows us to find and incorporate the suitable
number of indices $d$. 2) Second, we do not specify any functional form of
the conditional probability. Thus, the conditional probability in (\ref{eq:model})
is simply a function of both $t$ and $\bb_0\trans\X$. This relaxes
both the exponential form of the covariate relation and the
multiplicative form of the hazard function in the Cox model and is
also much more flexible than other popular semiparametric survival
models, such as the accelerated failure time and linear 
transformation models. Despite of the flexibility of the model in
(\ref{eq:model}), we show that through properly incorporating
semiparametric treatment and martingale techniques, estimation and
inference is still possible. 3) In addition, the analysis can be carried
out under the usual conditional independent censoring assumption, where the
censoring time is allowed to depend on the covariates. 

The proposed general index model and associated semiparametric
estimation method naturally provide a dimension reduction tool for
survival data. It has a few advantages over existing
dimension reduction methods for survival data. 1) First, many existing
dimension reduction methods for survival data require a stronger
assumption on the censoring time, such as the covariate-independent
censoring assumption \citep[e.g.][]{Lietal1999,LuLi2011}, or requires
nonparametric estimation of the conditional survival function of
censored survival times \citep{xia2010dimension} or censoring times \citep{Lietal1999}
given all the covariates, 
which may suffer from the curse of dimensionality.
All these drawbacks are avoided here.
2) Second, most of existing methods 
\citep{xia2010dimension, Lietal1999}
are constructed based on general inverse probability weighted estimation
techniques in one way or another, and are thus not efficient. In contrast, our proposed
method is built on the semiparametric theory \citep{tsiatis2006} and
achieves the optimal semiparametric efficient estimator.

The rest of the paper is organized as the following. In Section
\ref{sec:method}, we develop the estimation procedures for both the
index parameters in $\bb$ and functional relation between event time
and the multiple indices. In Section \ref{sec:asym}, we establish the
large sample properties to enable inference. We perform extensive
numerical experiments in Section \ref{sec:numeric}, where both
simulation and analysis of an AIDS data are included. We conclude the
paper with a discussion in Section \ref{sec:discuss}, while relegate
all the technical details in an Appendix.

\section{Methodology Development}\label{sec:method}

Define $Z=\min(T,C)$ and $\Delta=I(T\le C)$, where $C$ is the
censoring time. Assume $C\indep T\mid \X$ and the 
relation between $T$ and $\X$ follows the model in (\ref{eq:model}).
The observed 
data consist of $(\X_i, Z_i,\Delta_i)$, $i=1, \dots, n$, which are
independent copies of $(\X, Z, \Delta)$.  Note that even without 
censoring, $\bb_0$ in (\ref{eq:model}) is not identifiable because 
for any $d\times d$ full rank matrix $\A$, $\bb_0$ and $\bb_0\A$
suite the model (\ref{eq:model}) equally well. Thus, 
we fix a parameterization of $\bb_0$ by assuming the upper $d\times d$ block of $\bb_0$ 
to be the identity matrix $\I_d$. This ensures the unique
identification of $\bb_0$. Here we consider a fixed $d$, and our focus will be in estimating the
lower block of $\bb_0$, which has dimension $(p-d)\times d$. We then proceed to estimate the
conditional distribution function in (\ref{eq:model}).
For convenience, 
write $\X=(\X_u\trans, \X_l\trans)\trans$, where $\X_u\in{\cal R}^d$ and 
$\X_l\in{\cal R}^{p-d}$. 
Note that under the assumption of $C\indep
T\mid \X$ and (\ref{eq:model}), 
we have 
\bse
E\{f_1(C)f_2(T)\mid\bb_0\trans\X\}
&=&E[E\{f_1(C)f_2(T)\mid\X\}\mid\bb_0\trans\X]\\
&=&E[E\{f_1(C)\mid\X\}E\{f_2(T)\mid\bb_0\trans\X\}\mid\bb_0\trans\X]\\
&=&E\{f_1(C)\mid\bb_0\trans\X\}E\{f_2(T)\mid\bb_0\trans\X\}
\ese
for any functions $f_1, f_2$, hence
$C\indep T\mid\bb_0\trans\X$.

To facilitate further development, in describing the censoring
process, write
$S_c(z,\x)=\pr (C\ge z\mid \X=\x)$,
$\Lambda_c(z,\x)=-\log S_c(z,\x)$, 
$\lambda_c(z,\x)=\partial\Lambda_c(z,\x)/\partial z$ and 
$f_c(z,\x)=-\partial S_c(z,\x)/\partial z$. Similarly, to describe the
event process, 
for any parameter vector $\bb$, define
$S(z,\bb\trans\x)=\pr (T\ge z\mid \bb\trans\X=\bb\trans\x)$,
$f(z,\bb\trans\x)=-\partial S(z,\bb\trans\x)/\partial z$,
$\Lambda(z,\bb\trans\x)=-\log S(z,\bb\trans\x)$ and
$\lambda(z,\bb\trans\x)=\partial \Lambda(z,\bb\trans\x)/\partial z$.
Note that the functions $S, \Lambda, \lambda$ and $f$ are in fact
different when the parameter $\bb$ 
changes, so the more precise notations are
$S(z,\bb\trans\x,\bb), \Lambda(z,\bb\trans\x,\bb), \lambda(z,\bb\trans\x,\bb)$
and $f(z,\bb\trans\x,\bb)$. Here, without causing confusion, we omit the last
parameter for notational simplicity. 
Using these notation, the pdf of the model described in (\ref{eq:model}) is
\be\label{eq:pdf}
f_{\X,Z,\Delta}(\x,z,\delta,\bb,\lambda,\lambda_c,f_X)
=f_\X(\x)
\lambda(z,\bb\trans\x)^\delta
e^{-\int_0^z\lambda(s,\bb\trans\x)ds}
\lambda_c(z,\x)^{1-\delta}
e^{-\int_0^z\lambda_c(s,\x)ds},
\ee
where $f_\X(\x)$ is the pdf of $\X$. We assume the true data
generation process is based on
$f_{\X,Z,\Delta}(\x,z,\delta,\bb_0,\lambda_0,\lambda_{c0},f_{X0})$. 
Note that following our notation, 
$\lambda_0(z,\bb_0\trans\x)=\lambda(z,\bb_0\trans\x)$.

We view (\ref{eq:pdf}) as a semiparametric model, where $\bb$ is a
finite dimensional parameter of interest and all the remaining unknown
components of the model are treated as infinite dimensional nuisance
parameters. In survival analysis, the most popular approaches to
estimation are martingale based estimators \citep{fh1991} and
nonparametric maximum 
likelihood estimators (NPMLE) \citep{zl2007}. Here we find that NPMLE
does not suit well without adaption due to the inseparable relation
between the hazard function and the covariates. Martingale approach
may enable us to obtain one specific estimator for $\bb$, while we aim
at obtaining a more comprehensive understanding of the estimation of
$\bb$. Thus we use a less conventional approach by adopting the
geometrical treatment in semiparametrics. Similar practice has been
performed in \cite{tsiatis2006} to rediscover the partial likelihood
estimator for Cox proportional hazard model. To this end, we first
characterize the nuisance tangent space as described in Proposition
\ref{pro:nuisance}. The proof utilizes properties of martingale
integration and the details are given in Appendix
\ref{sec:pronuisance}. Define $M(t,\bb_0\trans\X)\equiv N(t)- \int_0^t
Y(s)\lambda_0(s,\bb_0\trans\X)ds$ and $M_c(t,\X)\equiv N_c(t)-\int_0^t
Y(s)\lambda_c(s,\X)ds$, where $N(t)=\Delta I(Z \le t)$, $N_c(t) =
(1-\Delta)I(Z \le t) $ and $Y(t)=I(Z\ge t)$. Then
$M(t, \bb_0\trans\X)$ and $M_c(t,\X)$ are mean-zero martingale
processes.

\begin{Pro}\label{pro:nuisance}
	
	The nuisance tangent space 
	$\Lambda=\Lambda_1\oplus\Lambda_2\oplus\Lambda_3$, where
	\bse
	\Lambda_1&=&\left[\a(\X): E\{\a(\X)\}=\0, \a(\X)\in{\cal R}^{(p-d)d}\right],\\
	\Lambda_2&=&\left\{\int_0^\infty\h(s,\bb_0\trans\X)dM(s,\bb_0\trans\X):\forall\h(Z,\bb_0\trans\X)\in{\cal R}^{(p-d)d}\right\},\\
	\Lambda_3&=&\left\{\int_0^\infty\h(s,\X)dM_c(s,\X):\forall\h(Z,\X)\in{\cal R}^{(p-d)d}\right\}.
	\ese
\end{Pro}

Having found the nuisance tangent space, we can now proceed to
identify the efficient score function through projecting the 
score function onto $\Lambda$ and calculating the residual. 
The score function is
defined as $\bS_{\bb}(\Delta,Z,\X)\equiv\partial\log
f_{\X,Z,\Delta}(\x,z,\delta,\bb,\lambda,\lambda_c,f_X)/\partial\bb$.
Let
$\blam_1(s,\bb\trans\X)=\partial\lambda(s,\bb\trans\X)/\partial(\bb\trans\X)$
and
$\blam_{10}(s,\bb_0\trans\X)=\partial\lambda_0(s,\bb_0\trans\X)/\partial(\bb_0\trans\X)$. 
Straightforward calculation yields
\bse
\bS_{\bb}(\Delta,Z,\X) 
=\int_0^\infty
\frac{\blam_{1}(s,\bb\trans\X)}{\lambda(s,\bb\trans\X)}\otimes\X_l dM(s,\bb\trans\X).
\ese
We can verify that $\bS_{\bb_0}(\Delta,Z,\X)\perp\Lambda_1$ and
$\bS_{\bb_0}(\Delta,Z,\X)\perp\Lambda_3$ due to the martingale
properties. 
Thus to look for the efficient score, we only need to project
$\bS_{\bb}(\Delta,Z,\X)$ onto $\Lambda_2$ and calculate its
residual. To this end,
we search for $\h^*(s,\bb_0\trans\X)$ so that
\bse
\bS\eff(\Delta,Z,\X)&=&\bS_{\bb_0}(\Delta,Z,\X)-\int_0^\infty\h^*(s,\bb_0\trans\X)dM(s,\bb_0\trans\X)\\
&=&\int_0^\infty\left\{
\frac{\blam_{10}(s,\bb_0\trans\X)}{\lambda_0(s,\bb_0\trans\X)}\otimes\X_l-
\h^*(s,\bb_0\trans\X)\right\}
dM(s,\bb_0\trans\X)
\ese
is orthogonal to $\Lambda_2$. 
This entails that for any $\h(s,\bb_0\trans\X)$, 
\bse
0&=&E\left[
\int_0^\infty\h\trans(s,\bb_0\trans\X)dM(s,\bb_0\trans\X)
\int_0^\infty\left\{
\frac{\blam_{10}(s,\bb_0\trans\X)}{\lambda_0(s,\bb_0\trans\X)}\otimes\X_l-
\h^*(s,\bb_0\trans\X)\right\}
dM(s,\bb_0\trans\X)\right]\\
&=&E\left[
\int_0^\infty \h\trans(s,\bb_0\trans\X)\left\{
\frac{\blam_{10}(s,\bb_0\trans\X)}{\lambda_0(s,\bb_0\trans\X)}\otimes\X_l-
\h^*(s,\bb_0\trans\X)\right\}Y(s)\lambda_0(s,\bb_0\trans\X)ds\right].
\ese
By letting $\h(s,\bb_0\trans\X)=I(s=t)\a(\bb_0\trans\X)$ for any
$\a(\bb_0\trans\X)$, we obtain that
\bse
0&=&E\left[\left\{
\frac{\blam_{10}(t,\bb_0\trans\X)}{\lambda_0(t,\bb_0\trans\X)}\otimes\X_l-
\h^*(t,\bb_0\trans\X)\right\}Y(t)\lambda_0(t,\bb_0\trans\X)\mid\bb_0\trans\X\right]\\
&=&E\left[\left\{
\frac{\blam_{10}(t,\bb_0\trans\X)}{\lambda_0(t,\bb_0\trans\X)}\otimes\X_l-
\h^*(t,\bb_0\trans\X)\right\}Y(t)\mid\bb_0\trans\X\right].
\ese
Note that
\be\label{eq:useful}
\frac{E\left\{ \X_lY(t)\mid\bb_0\trans\X\right\}}
{E\left\{Y(t)\mid\bb_0\trans\X\right\}}
=\frac{E\left\{\X_l 
	S_c(t,\X)\mid\bb_0\trans\X\right\}}
{E\left\{S_c(t,\X)\mid\bb_0\trans\X\right\}}.
\ee
This leads to
\bse
\h^*(t,\bb_0\trans\X) 
\frac{\blam_{10}(t,\bb_0\trans\X)}{\lambda_0(t,\bb_0\trans\X)}\otimes \frac{E\left\{\X_l 
	S_c(t,\X)\mid\bb_0\trans\X\right\}}
{E\left\{S_c(t,\X)\mid\bb_0\trans\X\right\}}.
\ese
Thus the efficient score is
\be\label{eq:effscore}
\bS\eff(\Delta,Z,\X)=
\int_0^\infty 
\frac{\blam_{10}(s,\bb_0\trans\X)}{\lambda_0(s,\bb_0\trans\X)}
\otimes\left[\X_l-
\frac{E\left\{\X_l 
	S_c(s,\X)\mid\bb_0\trans\X\right\}}
{E\left\{S_c(s,\X)\mid\bb_0\trans\X\right\}}\right]dM(s,\bb_0\trans\X).
\ee

Several properties of the efficient score is worth pointing out. 
First of all, 
$E\{\bS\eff(\Delta,Z,\X)\mid\X\}=\0$ is ensured by
$E\{dM(t,\bb_0\trans\X)\mid\X\}=0$, hence to preserve the mean zero property, we
can replace $\blam_{10}(s,\bb_0\trans\X)/\lambda_0(s,\bb_0\trans\X)$ by any
function of $s$ and $\bb_0\trans\X$, say $\g(s,\bb_0\trans\X)$, and still
obtain
\bse
E\int_0^\infty\g(s,\bb_0\trans\X)
\otimes\left[\X_l-
\frac{E\left\{\X_l 
	S_c(s,\X)\mid\bb_0\trans\X\right\}}
{E\left\{S_c(s,\X)\mid\bb_0\trans\X\right\}}\right]dM(s,\bb_0\trans\X)=\0.
\ese
This implies that if we are only aiming at a consistent estimator, we
can use an arbitrary function $\g(s,\bb_0\trans\X)$ to replace 
${\blam_{10}(s,\bb_0\trans\X)}/{\lambda_0(s,\bb_0\trans\X)}$ in the
efficient score to get a more general martingale integration.

Second, using (\ref{eq:useful}), by first taking expectation
conditional on $\bb\trans\X$, we can verify that
\bse
&&E\int_0^\infty\g(s,\bb_0\trans\X) 
\otimes\left[\X_l-
\frac{E\left\{\X_l 
	S_c(s,\X)\mid\bb_0\trans\X\right\}}
{E\left\{S_c(s,\X)\mid\bb_0\trans\X\right\}}\right]
Y(s)\lambda_0(s,\bb_0\trans\X)ds\\
&=&E\int_0^\infty\g(s,\bb_0\trans\X) 
\otimes\left[E\{\X_lY(s)\mid\bb_0\trans\X\}-
\frac{E\left\{\X_l 
	S_c(s,\X)\mid\bb_0\trans\X\right\}}
{E\left\{S_c(s,\X)\mid\bb_0\trans\X\right\}}
E\{Y(s)\mid\bb_0\trans\X\}\right]\\
&&\times\lambda_0(s,\bb_0\trans\X)ds\\
&=&\0.
\ese
As a consequence,
\bse
E\int_0^\infty\g(s,\bb_0\trans\X) 
\otimes\left[\X_l-
\frac{E\left\{\X_l 
	S_c(s,\X)\mid\bb_0\trans\X\right\}}
{E\left\{S_c(s,\X)\mid\bb_0\trans\X\right\}}\right]dN(s)=\0.
\ese
This implies that we can construct estimating equations of the form
\be\label{eq:general}
\sumi \Delta_i\g(Z_i,\bb\trans\X_i) 
\otimes\left[\X_{li}-
\frac{\wh E\left\{\X_{li} 
	Y_i(Z_i)\mid\bb\trans\X_i\right\}}
{\wh E\left\{Y_i(Z_i)\mid\bb\trans\X_i\right\}}\right]=\0
\ee
for any $\g$,
where
\be 
\wh E\left\{Y_i(Z_i)\mid\bb\trans\X_i\right\}
&=& \frac{\sumj K_h(\bb\trans\X_j-\bb\trans\X_i)I(Z_j\ge Z_i)}
{\sumj K_h(\bb\trans\X_j-\bb\trans\X_i)},\label{eq:expectY}\\
\wh E\left\{\X_{li}Y_i(Z_i)\mid\bb\trans\X_i\right\}
&=& \frac{\sumj K_h(\bb\trans\X_j-\bb\trans\X_i)\X_{lj}I(Z_j\ge Z_i)}
{\sumj K_h(\bb\trans\X_j-\bb\trans\X_i)}.\label{eq:expectXY}
\ee
Here $E\{Y_i(Z_i)\mid \bb\trans\X_i\}\equiv E\{Y_i(t)\mid\bb\trans\X_i\}|_{t=Z_i}$
and similarly for other terms, $K(\cdot)$ is a kernel function and $K_h(\cdot)=K(\cdot/h)/h$.

Third, we can further relax the estimating equation form to
\be\label{eq:silly}
\sumi \Delta_i\g(Z_i,\bb\trans\X_i) 
\otimes\left[\a(\X_{li})-
\frac{\wh E\left\{\a(\X_{li}) 
	Y_i(Z_i)\mid\bb\trans\X_i\right\}}
{\wh E\left\{Y_i(Z_i)\mid\bb\trans\X_i\right\}}\right]=\0
\ee
by taking advantage of the fact that
\bse
E\Delta\g(Z,\bb_0\trans\X) 
\otimes\left[\a(\X_{l})-
\frac{E\left\{\a(\X_{l}) 
	Y(Z)\mid\bb_0\trans\X\right\}}
{E\left\{Y(Z)\mid\bb_0\trans\X\right\}}\right]=\0
\ese
for any $\a(\X_l)$.

Fourth,
when we choose to estimate 
$\blam_{10}(s,\bb_0\trans\X)/\lambda_0(s,\bb_0\trans\X)$, 
using, for example, the smoothed Kaplan-Meier estimator,
we can then obtain the
efficient estimator from solving 
\be\label{eq:eff}
\sumi \Delta_i\frac{\wh\blam_1(Z_i,\bb\trans\X_i)}{\wh\lambda(Z_i,\bb\trans\X_i)}
\otimes\left[\X_{li}-
\frac{\wh E\left\{\X_{li} 
	Y_i(Z_i)\mid\bb\trans\X_i\right\}}
{\wh E\left\{Y_i(Z_i)\mid\bb\trans\X_i\right\}}\right]=\0.
\ee
Here we can use a 
local Nelson-Aalen estimator to estimator the cumulative hazard function
\bse
\widehat{\Lambda}(Z,\bb\trans\X) = \sum_{Z_i \le Z} \frac{\Delta_iK_h(\bb\trans\X_i-\bb\trans\X)}{\sumj I(Z_j\ge
	Z_i)K_h(\bb\trans\X_j-\bb\trans\X)}.
\ese
The local Nelson-Aalen estimator of 
hazard function can be obtained from
\be
\widehat{\lambda}(Z,\bb\trans\X) &=& \int_0^{\infty} K_b(t-Z) d\widehat{\Lambda}(t|\bb\trans\X)\nonumber\\
&=& \sum_{i=1}^nK_b(Z_i-Z) \frac{\Delta_iK_h(\bb\trans\X_i-\bb\trans\X)}{\sumj I(Z_j\ge
	Z_i)K_h(\bb\trans\X_j-\bb\trans\X)},\label{eq:lambda}
\ee
and we estimate the derivative
\be
\wh\blam_1(Z,\bb\trans\X)&=& 
\partial\wh{\lambda}(Z,\bb\trans\X)/\partial(\bb\trans\X)\nonumber\\
&=&-\sum_{i=1}^nK_b(Z_i-Z) \frac{\Delta_i\K_h'(\bb\trans\X_i-\bb\trans\X)}{\sumj I(Z_j\ge
	Z_i)K_h(\bb\trans\X_j-\bb\trans\X)}\nonumber\\
& & + \sum_{i=1}^nK_b(Z_i-Z)\Delta_iK_h(\bb\trans\X_i-\bb\trans\X) \frac{\sumj I(Z_j\ge
	Z_i)\K_h'(\bb\trans\X_j-\bb\trans\X)}{\{\sumj I(Z_j\ge Z_i)K_h(\bb\trans\X_j-\bb\trans\X)\}^2}.\notag\\\label{eq:lambda1}
\ee
Here 
$\K_h'(\v)=\partial K_h(\v)/\partial\v$ is the first derivative of $K_h$ with respect to its 
variables, which is a vector, 
and $b$ is a bandwidth. For any vector or matrix $\a$, let $\a^{\otimes2}=\a\a\trans$.

Among the different constructions of  consistent estimators, the
estimator obtained from (\ref{eq:eff}) will be shown to achieve the
smallest possible variability, hence this estimator is efficient and
is what we recommend.
The efficient estimator will be the focus of our study. We provide the detailed
algorithm of the efficient estimation procedure below. 

\begin{enumerate}
	\item Obtain an initial estimator of $\bb$ through,
	for example, hmave \citep{xia2010dimension}. Denote the result $\wt{\bb}$.
	
	\item Replacing
	$E\{Y(Z)\mid\bb\trans\X\}$,
	$E\{\X_{l} 
	Y(Z)\mid\bb\trans\X\}$, 
	$\lambda(Z,\bb\trans\X)$ and $\blam_1(Z,\bb\trans\X)$
	with their nonparametric estimated versions given in (\ref{eq:expectY}), (\ref{eq:expectXY}),
	(\ref{eq:lambda}) and (\ref{eq:lambda1})  
	respectively. Write the resulting
	estimators as 
	$\wh
	E\{\X_{l} 
	Y(Z)\mid\bb\trans\X\}$, $\wh
	E\{Y(Z)\mid\bb\trans\X\}$, $\wh\lambda(Z,\bb\trans\X)$
	and $\wh\blam_1(Z,\bb\trans\X)$.
	
	\item Plug 
	$\wh
	E\{\X_{l} 
	Y(Z)\mid\bb\trans\X\}$, $\wh E\{Y(Z)\mid\bb\trans\X\}$,
	$\wh\lambda(Z,\bb\trans\X)$ and $\wh\blam_1(Z,\bb\trans\X)$
	into (\ref{eq:eff}) and solve the estimating equation to obtain the
	efficient estimator $\wh\bb$, using $\wt\bb$ as starting value. 
	
\end{enumerate}
In performing the nonparametric estimations 
, bandwidths need to
be selected. Because the final estimator is  insensitive to the bandwidths, as
indicated by condition \ref{assum:bandwidth}, Lemma \ref{lem:pre},
Theorems \ref{th:consistency} and \ref{th:eff}, where a range of
different 
bandwidths all lead to the same asymptotic property, we
suggest that one should select the corresponding bandwidths by taking the sample
size $n$ to its suitable power to satisfy \ref{assum:bandwidth}, and
then multiply a constant to scale it.  For example, when $d=1$, we let
$h$ be $n^{-1/4-1/32}$ multiplying the standard deviation of $\wt\bb\trans\X_i$, let $b$ be $n^{-1/4+1/8}$ multiplying the standard
deviation of $Z_i$.

\section{Asymptotics}\label{sec:asym}

We will show that the efficient estimator described in Section
\ref{sec:method} is root-$n$ consistent, asymptotically normally
distributed  and achieves the optimal efficiency.
Let the parameter space of $\bb$ be
$\cal B$. We first list
some regularity conditions. 

\begin{enumerate}[label=C\arabic*]
	
	\item\label{assum:kernel} ({\itshape The kernel function.})
	The univariate kernel function $K(x)$ is symmetric, monotonically
	decreasing when $x>0$ and differentiable, with bounded derivative. In addition,
	$\int x^j K(x)dx=0$, for $1\le j<\nu$, $0<\int x^\nu K(x)dx<\infty$,
	$\int K^2(x)dx<\infty$, 
	$\int x^2 K^2(x)dx<\infty$,
	$\int K'^2(x)dx<\infty$, 
	$\int x^2K'^2(x)dx<\infty$,
	$\int K''^2(x)dx<\infty$, 
	$\int x^2K''^2(x)dx<\infty$. 
	The $d$-dimension kernel function is a product of $d$ univariate
	kernel functions, that is $K(\u)=\prod_{j=1}^{d}K(u_j)$ for
	$\u=(u_1,...,u_d)\trans$. For simplicity, we use the same $K$ for both
	univariate and multivariate kernel functions. 
	\item\label{assum:bandwidth}({\itshape The bandwidths.})
	The bandwidths satisfy $h\to0$, $b\to0$,
	$nh^{d+2}b\to\infty$ and 
	$nh^{2\nu}\to0$, where $2\nu> d+1$.
	
	\item\label{assum:bounded}({\itshape The boundedness.})
	The parameter space ${\cal B}$ is bounded.

	\item\label{assum:fbeta}({\itshape The density functions of covariates.})
	For all $\bb\in{\cal B}$,
	the probability density function of $\bb\trans\X$,
	$f_{\bb\trans\X}(\bb\trans\x)$, has a compact support and has
	four derivatives. The function $f_{\bb\trans\X}(\bb\trans\x)$ is 
	bounded away from zero and infinity on the support, and its first four derivatives
	are bounded uniformly on the support.
	
	\item\label{assum:exi}({\itshape The smoothness.})
	For all $\bb\in{\cal B}$, 
	the absolute value of $E\{\X_jI(Z_j\ge Z)\mid\bb\trans\x\}$, $E\{I(Z_j\ge Z)\mid\bb\trans\x\}$,
	and their first four
	derivatives are bounded uniformly component wise.
	The absolute value of $E\{\X_j\X_j\trans I(Z_j\ge 
	Z)\mid\bb\trans\x\}$
	and its first two
	derivatives are bounded uniformly component wise.
	
	\item\label{assum:survivalfunction}({\itshape The survival function.})
	For all $\bb\in{\cal B}$, the survival function of the event process
	$S(t,\bb\trans\x)$, the conditional expectation of the survival function of the censoring
	processes 
	$E\{S_c(t,\X)\mid \bb\trans\x\}$ and the probability  density
	function of the survival process 
	$f(t,\bb\trans\x)$ satisfy $\partial^{i+j}S(t,\bb\trans\x)/\partial
	t^i\partial (\bb\trans\x)^j$, $\partial^{i+j}E\{S_c(t,\X)\mid
	\bb\trans\x\}/\partial t^i\partial (\bb\trans\x)^j$,
	$\partial^{i+j}f(t,\bb\trans\x)/\partial t^i\partial
	(\bb\trans\x)^j$ exist and  are bounded and bounded away from zero, for all $i\ge0, j\ge0, i+j\le 4$.

	\item\label{assum:unique}({\itshape The uniqueness.})
	The equation
	\bse
	E\left( \Delta\frac{\blam_1(Z,\bb\trans\X)}{\lambda(Z,\bb\trans\X)}
	\otimes\left[\X_{l}-
	\frac{ E\left\{\X_{l} 
		Y(Z)\mid\bb\trans\X\right\}}
	{E\left\{Y(Z)\mid\bb\trans\X\right\}}\right]\right)
	=\0
	\ese
	has a unique solution on $\cal B$.
	Because the true parameter $\bb_0$ satisfies the equation, 
	hence the unique solution is $\bb_0$.
	
\end{enumerate}

These conditions are quite commonly imposed in nonparametrics,
survival analysis and estimating equations and are generally mild.
Conditions \ref{assum:kernel} and \ref{assum:bandwidth} contain some
basic requirements on the kernel function and the bandwidths, which
are common in kernel related works. The boundedness of the parameter
space  $\cal B$ in \ref{assum:bounded} is satisfied in most
practical problems.
Condition \ref{assum:fbeta}-\ref{assum:survivalfunction} impose certain
boundedness condition of event time, censoring time, covariates, their
expectations and corresponding derivatives.
The unique solution requirement in \ref{assum:unique} is needed to
ensure the convergence of the estimator, which could be further
relaxed to local uniqueness if needed.

Before presenting the main results, we summarize several preliminary
results first. These results highlight the theoretical properties of the
kernel based estimators of several conditional expectations, as well
as the estimation properties of the hazard function and its
derivative, hence are of their own interest. These properties also play an
important role in the proof of the Theorems \ref{th:consistency} and
\ref{th:eff}.
The proofs of Lemma \ref{lem:pre} and both theorems  are respectively
in the Appendice \ref{sec:prooflempre}, \ref{sec:proofthconsistency} and 
\ref{sec:prooftheff}.

\begin{Lem}\label{lem:pre}
	Under the regularity conditions \ref{assum:kernel}-\ref{assum:unique}  listed above, 
	\be
	\wh E\left\{Y(Z)\mid\bb\trans\X\right\}
	&=&E\{Y(Z)\mid\bb\trans\X\}+O_p\{(nh)^{-1/2}+h^2\},\label{eq:lemeq1}\\
	\wh E\left\{\X Y(Z)\mid\bb\trans\X\right\}
	&=&E\{\X Y(Z)\mid\bb\trans\X\}+O_p\{(nh)^{-1/2}+h^2\}, \label{eq:lemeq2}\\
	\frac{\partial}{\partial\bb\trans\X}\wh E\left\{Y(Z)\mid\bb\trans\X\right\}
	&=&\frac{\partial}{\partial\bb\trans\X} E\{Y(Z)\mid\bb\trans\X\}+O_p\{(nh^3)^{-1/2}+h^2\}\label{eq:lemeq3}\\
	\frac{\partial}{\partial\bb\trans\X}\wh E\left\{\X Y(Z)\mid\bb\trans\X\right\}
	&=&\frac{\partial}{\partial\bb\trans\X}E\{\X Y(Z)\mid\bb\trans\X\}+O_p\{(nh^3)^{-1/2}+h^2\} \label{eq:lemeq4}\\
	\widehat{\lambda}(z,\bb\trans\x)&=&\lambda(z,\bb\trans\x)+O_p\{(nhb)^{-1/2}+h^2+b^2\}\label{eq:lemeq5}\\
	\wh\blam_1(z,\bb\trans\x)&=&\blam_1(z,\bb\trans\x)+O_p\{(nbh^3)^{-1/2}+h^2+b^2\}\label{eq:lemeq6}
	\ee
	uniformly for all $Z,\bb\trans\X$.
\end{Lem}

\begin{Th}\label{th:consistency}
	The estimator obtained from solving (\ref{eq:eff}) is consistent,
	i.e. 
	$\wh\bb-\bb_0\to\0$ in probability when $n\to\infty$.
\end{Th}

\begin{Th}\label{th:eff}
	The estimator obtained from solving (\ref{eq:eff}) satisfies
	\bse
	\sqrt{n}(\wh\bb-\bb_0)\to N(\0, [E\{\bS\eff ^{\otimes2} (\Delta,
	Z,\X)\}]^{-1})
	\ese
	in distribution when $n\to\infty$. Here $\bS\eff(\Delta, Z, \X)$ is given
	in (\ref{eq:effscore}). Thus, the estimator is efficient. 
\end{Th}

Note that because $\bS\eff$ is a martingale, we have
\bse
&&E\{\bS\eff ^{\otimes2} (\Delta,
Z,\X)\}\\
&=&E\int_0^\infty 
\left(\frac{\blam_{10}(s,\bb_0\trans\X)}{\lambda_0(s,\bb_0\trans\X)}
\otimes\left[\X_l-
\frac{E\left\{\X_l 
	S_c(s,\X)\mid\bb_0\trans\X\right\}}
{E\left\{S_c(s,\X)\mid\bb_0\trans\X\right\}}\right]\right)^{\otimes2}
\lambda(s,\bb_0\trans\X)Y(s)ds\\
&=&E\int_0^\infty 
\left(\frac{\blam_{10}(s,\bb_0\trans\X)}{\lambda_0(s,\bb_0\trans\X)}
\otimes\left[\X_l-
\frac{E\left\{\X_l 
	S_c(s,\X)\mid\bb_0\trans\X\right\}}
{E\left\{S_c(s,\X)\mid\bb_0\trans\X\right\}}\right]\right)^{\otimes2}
dN(s).
\ese
Therefore,  a natural estimator of the estimation variance is the inverse of
\bse
\frac{1}{n}\sumi\delta_i
\left(\frac{\wh\blam_{1}(z_i,\wh\bb\trans\x_i)}{\wh\lambda(z_i,\wh\bb\trans\x_i)}
\otimes\left[\x_{il}-
\frac{\wh E\left\{\X_{l} 
	S_c(z_i,\X)\mid\wh\bb\trans\x_i\right\}}
{\wh E\left\{S_c(z_i,\X)\mid\wh\bb\trans\x_i\right\}}\right]\right)^{\otimes2}.
\ese

\section{Numerical Experiments}\label{sec:numeric}
\subsection{Simulation}

To evaluate the finite sample performance of our method, we perform
four simulation studies. 
In the first study, we generate event times from
\bse
T = \Phi\left[5\epsilon\left\{ \exp\left(3\bb\trans\X\right)+1\right\}-2\right]
\ese
where $\Phi$ is the cumulative distribution function (cdf) of the
standard normal distribution, $\epsilon$ has an exponential
distribution with parameter 1,  and $\X$ follows a standard normal distribution
independent with $\epsilon$. We consider $d=1, p=7$ and the true
parameter values are taken 
to be $\bb=(1, 0, -1, 0, 1, 0, -1)\trans$. We further generate the
covariate dependent 
censoring times using
$C = \Phi(2X_2+2X_3)+U$
where $U$ denotes a random variable uniformly distributed on
$(0,c_1)$, where $c_1$ is a constant controlling the proportion of censoring. 

In the second study, we generate the event times from
\bse
T = \exp(\bb\trans\X+\epsilon)
\ese
where $\epsilon$ follows an exponential distribution with parameter 1
and each component in $\X$ follows independent uniform  
distribution on $(0,\,\,\,1)$. We consider $d=1, p=7$ and set the true parameter
value to be $\bb=(1,1.3, -1.3,1,-0.5,0.5,-0.5)\trans$.
We generate the censoring time from a uniform distribution on
$(0,c_2)$, where different values of $c_2$ are used to achieve various censoring rate. 

In the third study, we generate the event times from
\bse
T=\exp\left\{5-10(1-\bb\trans\X)^2+\epsilon\right\}
\ese
where $\epsilon\sim \text{Normal}(0,1)$, and each component of $\X$ is
independently distributed with uniform distribution on $(0,\,1)$. We consider $d=1,p=10$ and
set the true parameter value to be
$\bb=(1, -0.6, 0, -0.3, -0.1, 0, 0.1, 0.3, 0,0.6)\trans$.
The censoring time is generated
from $C=U\bb_c\trans\X $ where $\bb_c=(0,0,0,1,1,0,0,0,0,0)\trans$
and $U$ is uniformly distributed on $(0,c_3)$, and $c_3$ is a constant
controlling the censoring proportion.

In the last simulation study, we increase $d$ to 2 to further evaluate the
performance of the proposed method.
We set the event times
\bse
T=\exp\left\{5-10\sum_{j=1}^2(1-\bb_j\trans\X)^2+\epsilon\right\}
\ese
where  $\epsilon\sim \text{Normal}(0,1)$ and each component of $\X$ is
independently distributed with uniform distribution on $(0,\,1)$,
$\bb_j,j=1,2,$ denotes the $j$th column of $\bb$ with $p=6$. The
censoring time is generated from a uniform distribution on 
$(0,c_4)$, where $c_4$ controls the censoring rate.

These studies are designed to resemble and extend the
simulation studies considered in \cite{xia2010dimension}, which
proposed hmave, the best method  so far in the literature for dimension
reduction under censored data. In fact,  we compared our results with
those from hmave. In all the studies, 
we generated 1000 data sets. In the first study, sample size 
$n=100$ was considered. We set the sample sizes to
$n=500$ for the second study and to $n=200$ for all the remaining studies.

\begin{figure}[H]
	\includegraphics[trim={0 0 0 0},width=19cm]{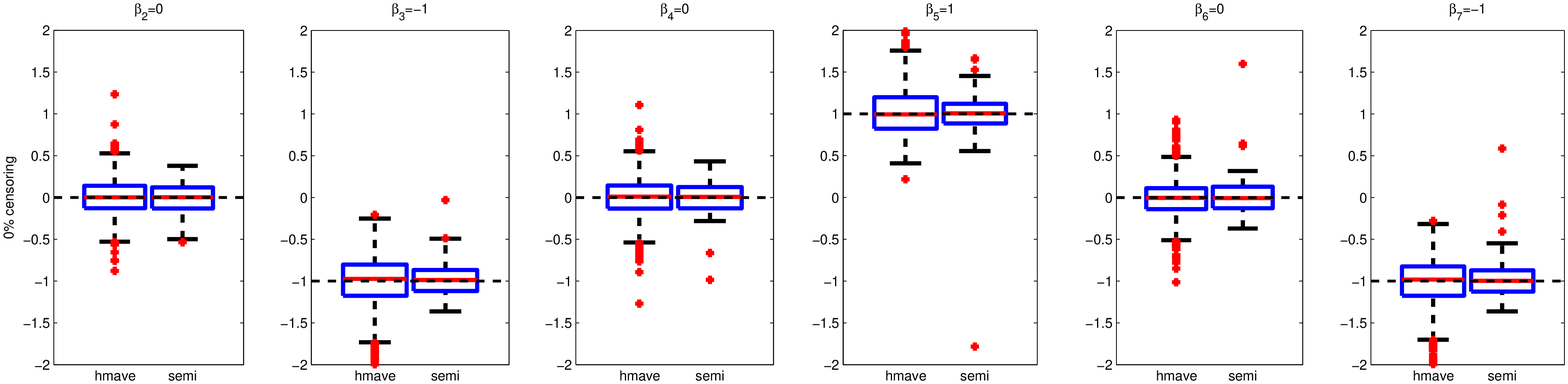}
	\includegraphics[trim={0 0 0 0},width=19cm]{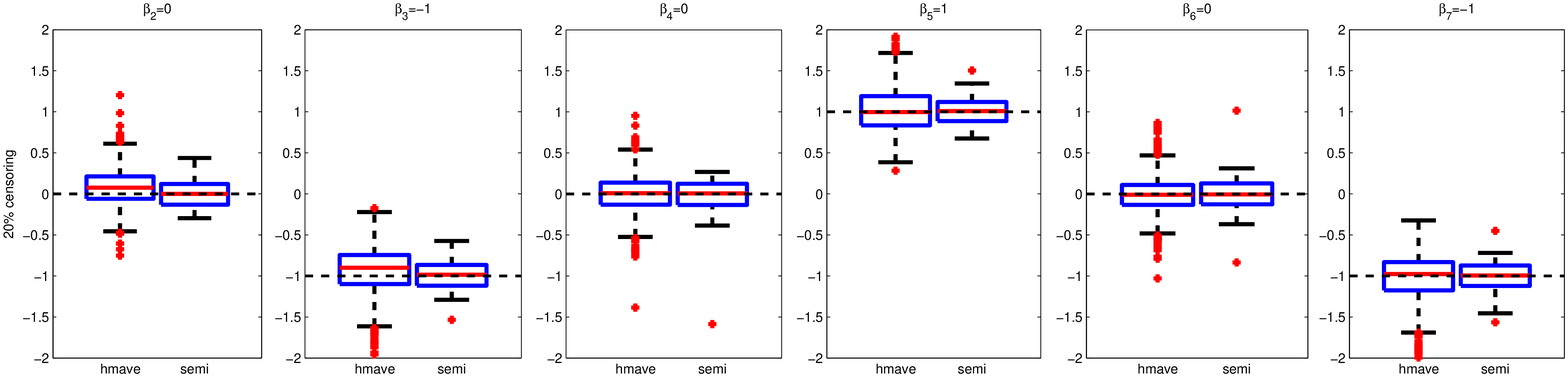}
	\includegraphics[trim={0 0 0 0},width=19cm]{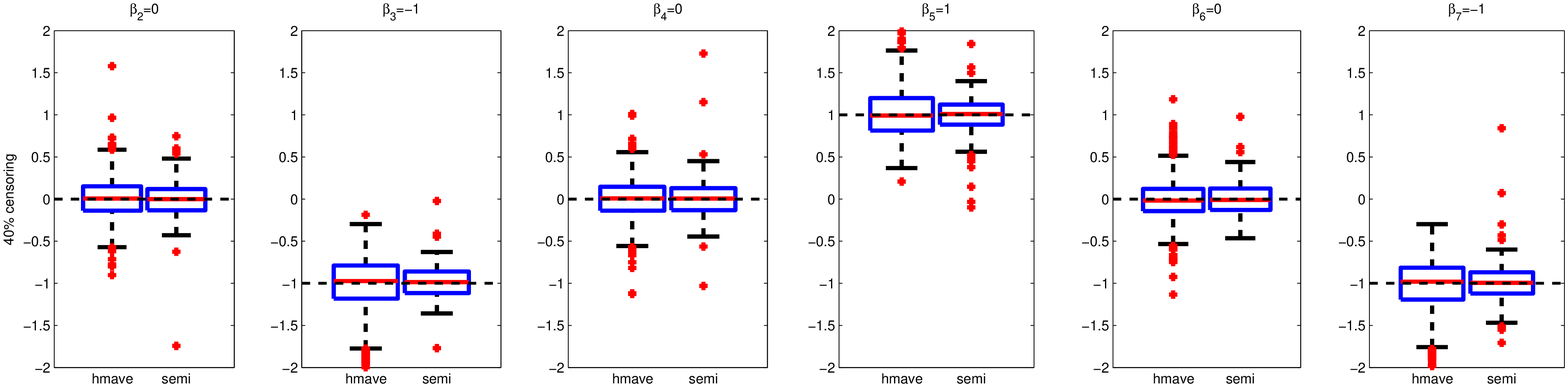}
	\caption{Boxplot of hmave and the semiparametric methods of study 1. 
		First row: no censoring; Second row: 20\% censoring rate; Third row:
		40\% censoring rate. Dashed line: True $\bb$.}
	\label{fig:simu1}
\end{figure}

The results of the first simulation study are given in Table \ref{tab:simu1}
and Figure \ref{fig:simu1}, where we considered three different
censoring rates, 0\%, 20\% and 40\% separately. 
From these results, we
can see that the semiparametric method we proposed generally performs
better, and often is much better than  hmave, in that it has much
smaller absolute biases as well as smaller sample standard errors. The
semiparametric method also yields smaller difference between the
estimated projection matrix $\wh\bP\equiv
\wh\bb(\wh\bb\trans\wh\bb)^{-1}\wh\bb\trans$ and the true projection
matrix $\bP\equiv \bb(\bb\trans\bb)^{-1}\bb\trans$, in that both the mean
and variance of the largest singular value of $\wh\bP-\bP$ are much
smaller based on the 
semiparametric method than based on hmave. The same results are also
presented in Figure \ref{fig:simu1} to provide a quick visual
inspection. 

\begin{table}[H]
	\centering
	\caption{Results of  study 1, based on 1000
		simulations with sample size 100. ``bias'' is the absolute
		value of $\text{mean}(\wh\bb)-\bb$ of each component in $\bb$, ``sd'' is the
		sample standard errors of the corresponding estimation.
		The last column is the mean and standard errors of the largest singular value of
		$\wh\bP-\bP$. }
	\label{tab:simu1}
	\vspace{.2cm}
	\begin{tabular}{cc|cccccc|c}\hline
		&& $\beta_2$ &$\beta_3$ &$\beta_4$ &$\beta_5$
		&$\beta_6$      &$\beta_7$ & $\lambda_{\max}$\\
		&true& $0$ &$-1$ &$0$ &$1$ &$0$      &$-1$ &  \\\hline
		&&\multicolumn{6}{c}{No censoring}&\\
		hmave&bias& 0.0010 & 0.0256 & 0.0099 & 0.0058 & 0.0052 & 0.0182 & 0.2208 \\
		& sd  & 0.1700 & 0.2344 & 0.1710 & 0.2320 & 0.1643 & 0.2264&  0.1578 \\
		semi	& bias & 0.0025& 0.0129 & 0.0071 & 0.0059 & 0.0071 & 0.0033  & 0.0903\\
		& sd  & 0.1298 & 0.1333 & 0.1314& 0.1277 & 0.1337 & 0.1335& 0.0622\\
		\hline
		&&\multicolumn{6}{c}{20\% censoring}&\\		
		hmave& bias  & 0.0747 & 0.0994 & 0.0095 & 0.0042 & 0.0099 & 0.0228&  0.2256 \\
		& sd   & 0.1688 & 0.2236 & 0.1663 & 0.2281 & 0.1612 & 0.2217  &     0.1560    \\
		semi	& bias & 0.0003 & 0.0143 & 0.0064 & 0.0079 & 0.0055 & 0.0054 &  0.0928  \\
		& sd & 0.1301 & 0.1339 & 0.1300 & 0.1268 & 0.1320 & 0.1294 &   0.0574\\
		\hline
		&&\multicolumn{6}{c}{40\% censoring}&\\		
		hmave& bias  & 0.0056 & 0.0261 & 0.0078 & 0.0079 & 0.0169 & 0.0189 &  0.2314 \\
		& sd  & 0.1812 & 0.2502 & 0.1784 & 0.2462 & 0.1707 & 0.2416&  0.1604  \\
		semi	& bias  & 0.0012 & 0.0130 & 0.0064 & 0.0090 & 0.0103 & 0.0056&  0.0948  \\
		& sd  & 0.1345 & 0.1352 & 0.1353 & 0.1305 & 0.1351 & 0.1354 &  0.0694\\
		\hline 
	\end{tabular}
\end{table}

\begin{figure}[H]

	\includegraphics[trim={0 0 0 0},width=19cm]{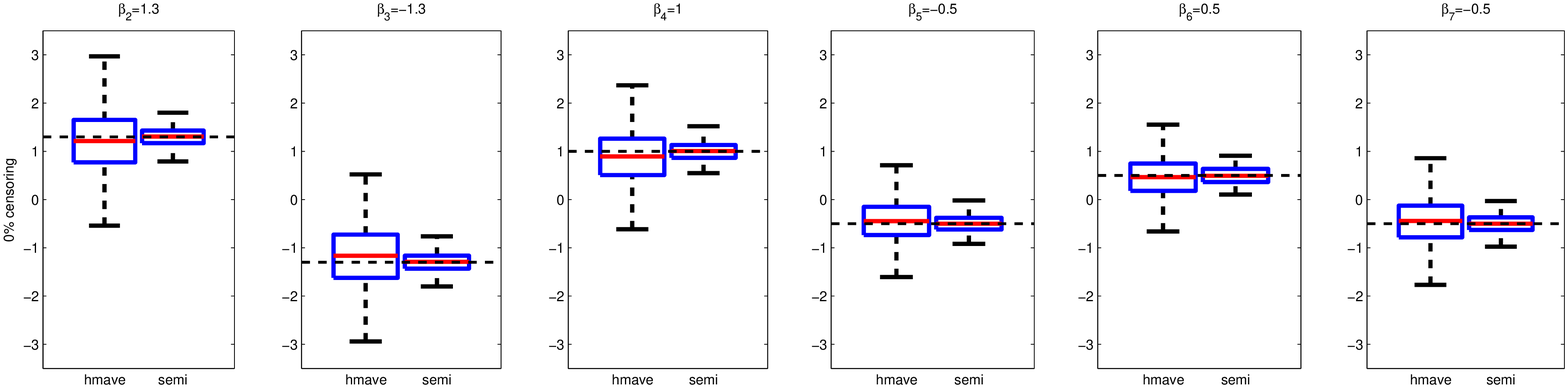}
	\includegraphics[trim={0 0 0 0},width=19cm]{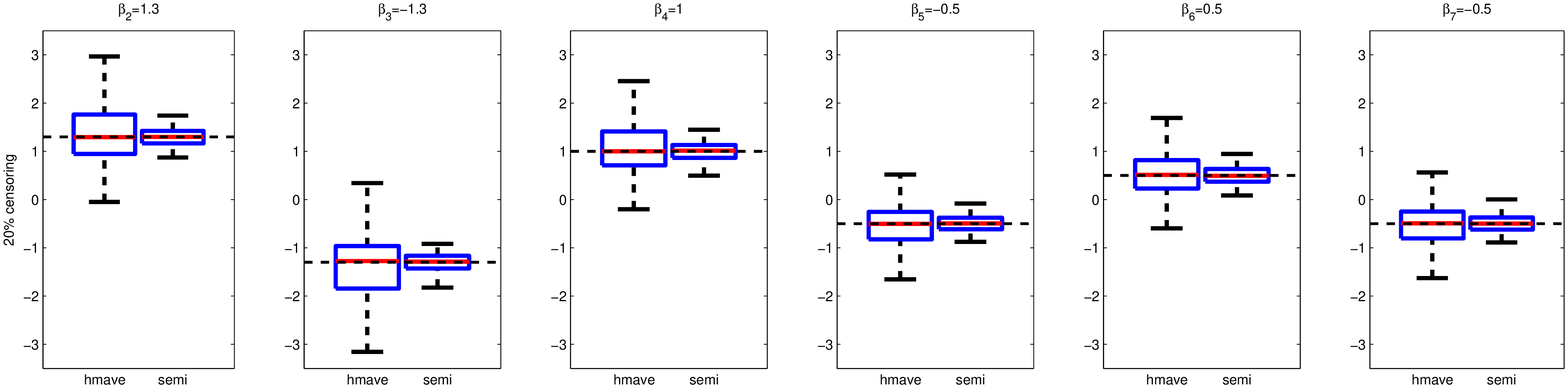}
	\includegraphics[trim={0 0 0 0},width=19cm]{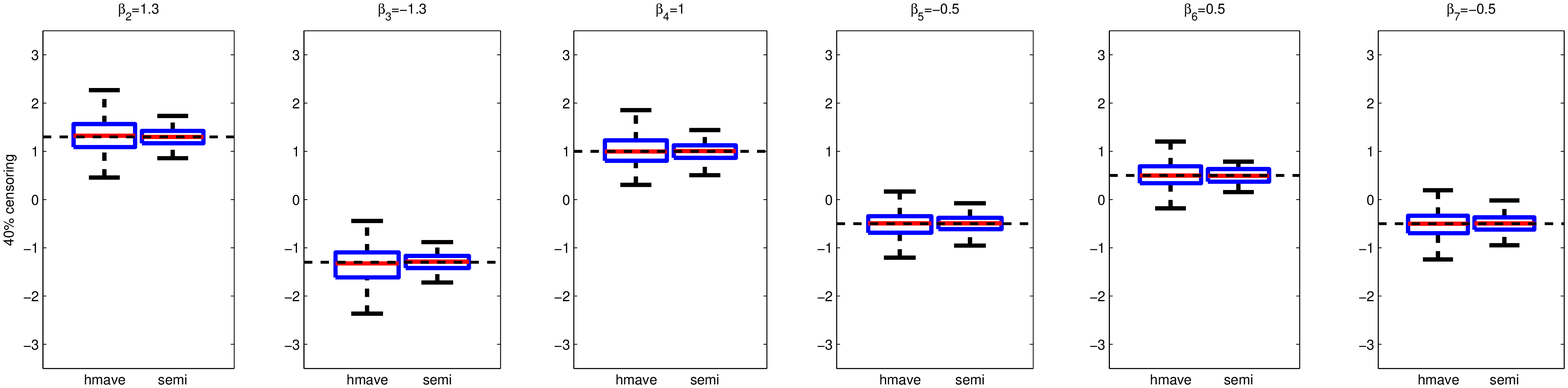}
	\caption{Boxplot of hmave and the semiparametric methods of study 2. 
		First row: no censoring; Second row: 20\% censoring rate; Third row:
		40\% censoring rate. Dashed line: True $\bb$.}
	\label{fig:simu2}
\end{figure}

\begin{figure}[H]
\includegraphics[trim={0 0 0 0},width=19cm]{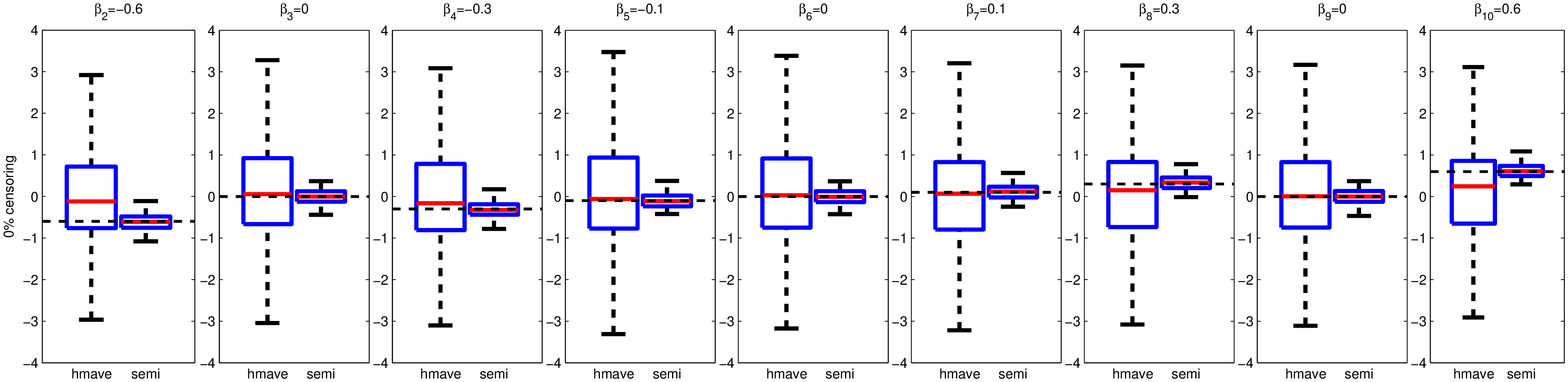}
\includegraphics[trim={0 0 0 0},width=19cm]{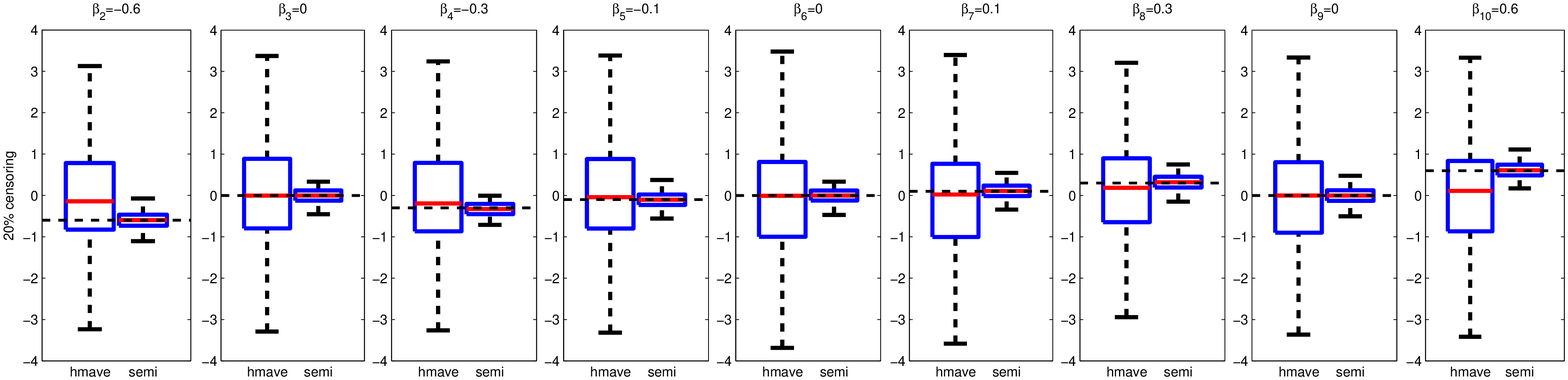}
\includegraphics[trim={0 0 0 0},width=19cm]{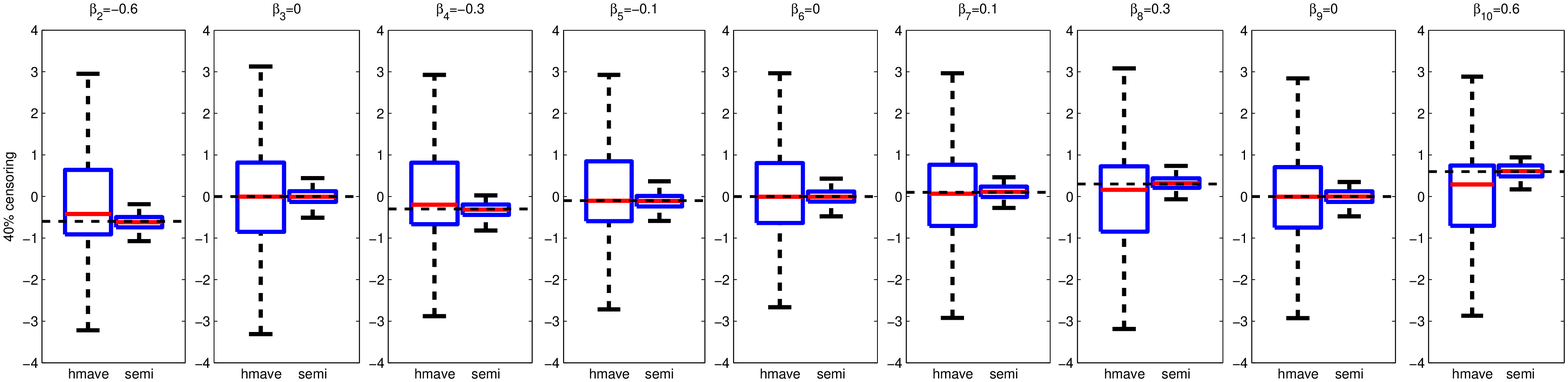}
	\caption{Boxplot of hmave and the semiparametric methods of study 3. 
		First row: no censoring; Second row: 20\% censoring rate; Third row:
		40\% censoring rate. Dashed line: True $\bb$.}
	\label{fig:simu3}
\end{figure}

\begin{figure}[H]

	\includegraphics[trim={0 0 0 0},width=19cm]{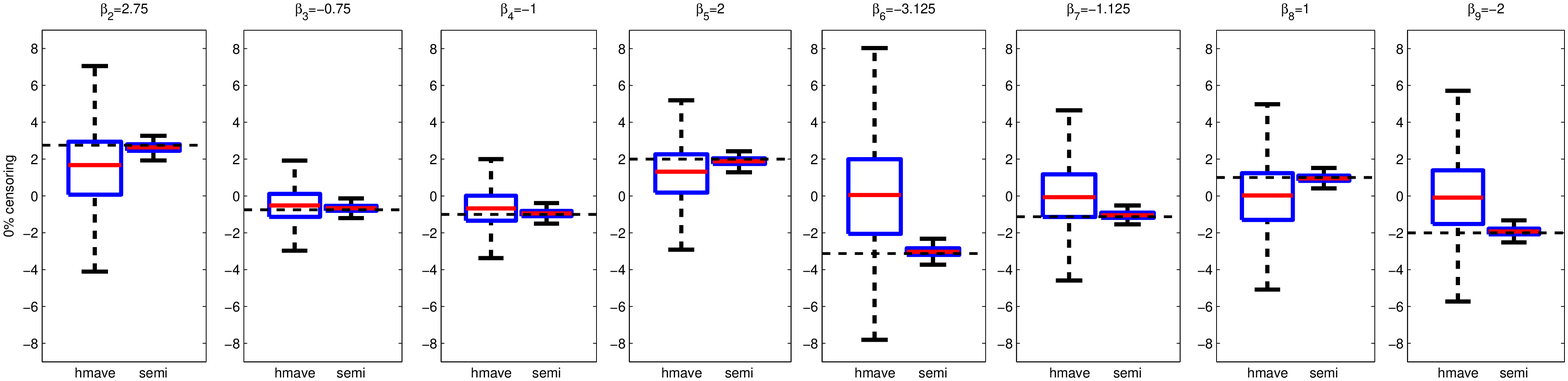}
	\includegraphics[trim={0 0 0 0},width=19cm]{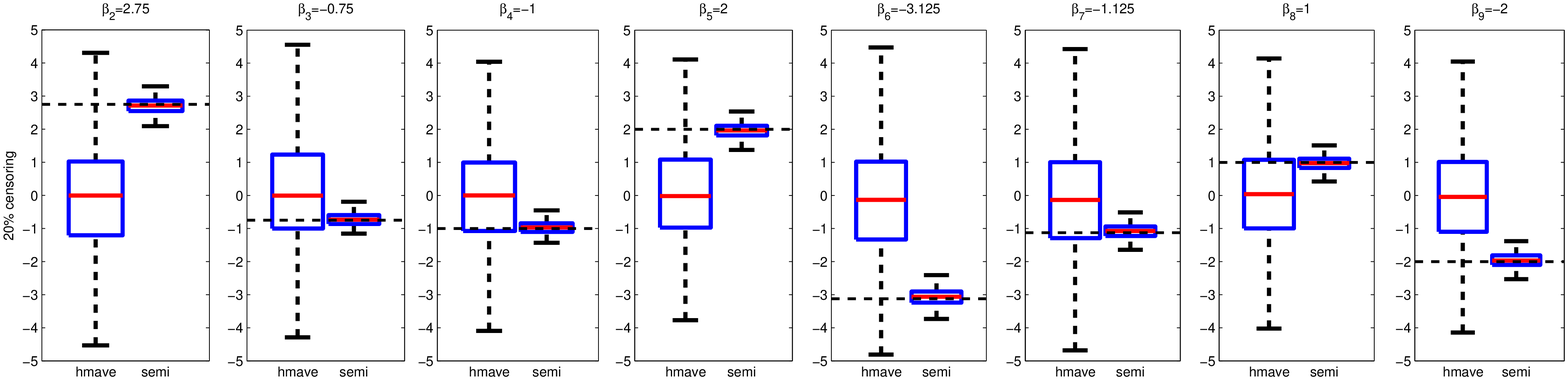}
	\includegraphics[trim={0 0 0 0},width=19cm]{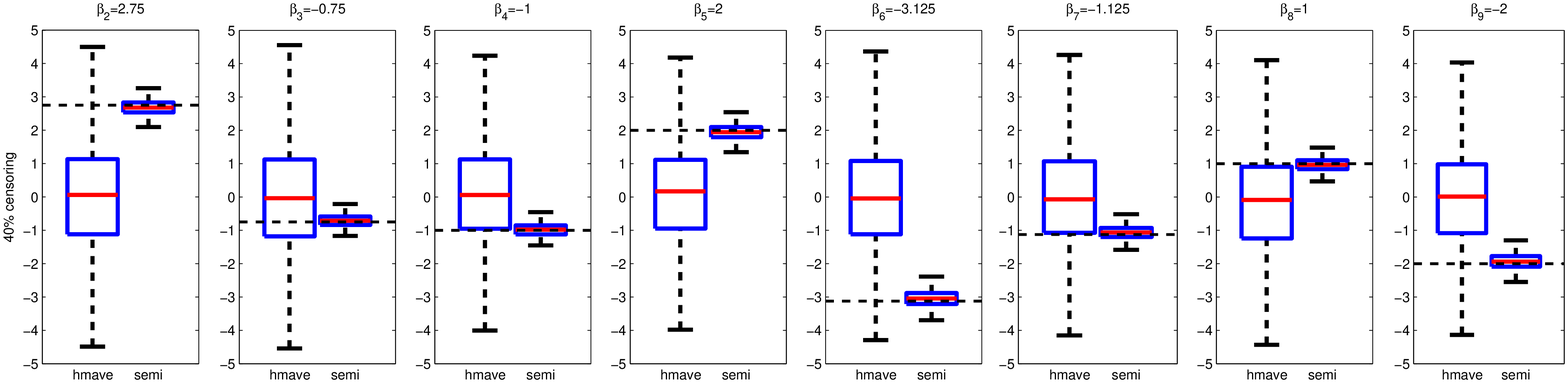}
	\caption{Boxplot of hmave and the semiparametric methods of study 4. 
		First row: no censoring; Second row: 20\% censoring rate; Third row:
		40\% censoring rate. Dashed line: True $\bb$.}
	\label{fig:simu4}
\end{figure}

The results of the second study study are presented in Tables  \ref{tab:simu2},
and Figures \ref{fig:simu2}, where
the same conclusion can be drawn as in the first study.
The superiority of the semiparametric method to hmave is even more prominent in the
third study, as reflected in Table \ref{tab:simu3} and
Figure \ref{fig:simu3}. Here, 
the semiparametric method is substantially more accurate in estimating
each component in $\bb$, yielding smaller biases and variances. 
The largest singular value of the difference between the estimated
and true projection matrices is also much smaller for the
semiparametric method in comparison with hmave.
When we increased $d$ to 2 in the last simulation, the semiparametric
method continues to generate satisfactory results, see
Table \ref{tab:simu4} and Figure \ref{fig:simu4}. In this case, the
performance of hmave is rather concerning, possibly caused by the
difficulties associated with multiple indices.

\begin{table}[H]
	\centering
	\caption{Results of  study 2, based on 1000
		simulations with sample size 200. ``bias'' is the absolute
		value of $\text{mean}(\wh\bb)-\bb$ of each component in $\bb$, ``sd'' is the
		sample standard errors of the corresponding estimation.
		The last column is the mean and standard errors of the largest singular value of
		$\wh\bP-\bP$. } \label{tab:simu2}
	\vspace{0.1cm}
	\begin{tabular}{cc|cccccc|c}\hline
		&& $\beta_2$ &$\beta_3$ &$\beta_4$ &$\beta_5$
		&$\beta_6$      &$\beta_7$ & $\lambda_{\max}$\\
		&true& $1.3$ &$-1.3$ &$1$ &$-0.5$ &$0.5$      &$-0.5$ &  \\\hline
		&&\multicolumn{6}{c}{No censoring}&\\
		hmave& bias  & 0.1958 & 0.3428 & 0.1644 & 0.0137 & 0.3699 & 0.1672 & 0.4247 \\
		& sd   & 7.0919 & 7.6671 & 9.8179 & 6.2261 & 10.976 & 5.6033 &  0.1442 \\
		semi	& bias  & 0.2483 & 0.0617 & 0.0929 & 0.2354 & 0.1273 & 0.0646 &  0.2915 \\
		& sd  & 5.3549 & 4.3952 & 2.5283 & 5.9618 & 3.2538 & 1.9437&  0.1133\\
		\hline
		&&\multicolumn{6}{c}{20\% censoring}&\\		
		hmave& bias & 0.5650 & 0.5097 & 0.3841 & 0.1947 & 0.2155 & 0.2656 & 0.3212\\
		& sd   & 4.2947 & 3.4980 & 2.4377 & 1.5136 & 2.2722 & 2.4448&  0.1356 \\
		semi	& bias    & 0.0864 & 0.0289 & 0.0167 & 0.0109 & 0.0537 & 0.1300 &   0.1448   \\
		& sd & 1.9872 & 0.5161 & 1.2467 & 0.8233 & 0.9729 & 3.8902&  0.0909\\
		\hline
		&&\multicolumn{6}{c}{40\% censoring}&\\		
		hmave& bias   & 0.0998 & 0.1024 & 0.0566 & 0.0381 & 0.0386 & 0.0371& 0.1991  \\
		& sd  & 0.4603 & 0.4678 & 0.3877 & 0.3085 & 0.2939 & 0.3085 &  0.1303   \\
		semi	& bias   & 0.0158 & 0.0107 & 0.0053 & 0.0128 & 0.0063 & 0.0062  &   0.0781   \\
		& sd   & 0.3976 & 0.4678 & 0.3492 & 0.5517 & 0.2718 & 0.2873& 0.0720\\
		\hline 
	\end{tabular}
\end{table}

\begin{table}[H]
	\centering
	\caption{Results of  study 3, based on 1000
		simulations with sample size 200. ``bias'' is the absolute
		value of $\text{mean}(\wh\bb)-\bb$ of each component in $\bb$, ``sd'' is the
		sample standard errors of the corresponding estimation.
		The last column is the mean and standard errors of the largest singular value of
		$\wh\bP-\bP$. } \label{tab:simu3}
	\vspace{0.1cm}
	\begin{tabular}{cc|ccccccccc|c}\hline
		&& $\beta_2$ &$\beta_3$ &$\beta_4$ &$\beta_5$
		&$\beta_6$      &$\beta_7$&$\beta_8$&$\beta_9$&$\beta_{10}$ & $\lambda_{\max}$\\
		&true& $-0.6$ &  0 & -0.3 &  -0.1 & 0 & 0.1 & 0.3 &  0 & 0.6 &  \\\hline
		&&\multicolumn{9}{c}{No censoring}&\\
		hmave& bias & 0.3711 & 0.3643 & 1.3056 & 0.5394 & 0.0913 & 0.1954 & 0.4499 & 0.2092 & 0.3660&   0.8706   \\
		& sd    & 17.267 & 10.472 & 35.110 & 16.246 & 9.449 & 16.931 & 26.510 & 12.342 & 17.304 &   0.2953  \\
		semi	& bias & 0.0124 & 0.0044 & 0.0175 & 0.0035 & 0.0026 & 0.0111 & 0.0323 & 0.0013 & 0.0180&   0.2337     \\
		& sd  & 0.1639 & 0.1538 & 0.1523 & 0.1563 & 0.1585 & 0.1535 & 0.1590 & 0.1543 & 0.1631& 0.0637\\
		\hline
		&&\multicolumn{9}{c}{20\% censoring}&\\		
		hmave& bias& 1.4974 & 2.3355 & 1.6021 & 0.4699 & 2.3620 & 1.6553 & 0.5311 & 1.6596 & 2.9426&   0.8822  \\
		& sd    & 41.451 & 44.735 & 66.302 & 49.280 & 40.673 & 47.488 & 58.485 & 57.025 & 72.228 &   0.2952  \\
		semi	& bias  & 0.0035 & 0.0003 & 0.0239 & 0.0063 & 0.0023 & 0.0072 & 0.0184 & 0.0035 & 0.0177 &    0.2148 \\
		& sd & 0.1600 & 0.1584 & 0.1722 & 0.1716 & 0.1555 & 0.1615 & 0.1531 & 0.1633 & 0.1595  &   0.0691
		\\
		\hline
		&&\multicolumn{9}{c}{40\% censoring}&\\		
		hmave& bias& 0.5909 & 5.5304 & 3.2835 & 0.8370 & 1.775 & 5.6482 & 4.8272 & 1.2442 & 0.6951&  0.8382 \\
		& sd & 20.946 & 146.58 & 68.000 & 25.442 & 38.877 & 145.82 & 90.530 & 58.178 & 23.032&   0.2900     \\
		semi	& bias & 0.0209 & 0.0004 & 0.0198 & 0.0062 & 0.0041 & 0.0138 & 0.0244 & 0.0021 & 0.0166 & 0.2656   \\
		& sd    & 0.1500 & 0.1514 & 0.1518 & 0.1506 & 0.1492 & 0.1560 & 0.1453 & 0.1503 & 0.1579&  0.0556	\\
		\hline 
	\end{tabular}
\end{table}

\begin{table}[H]
	\centering
	\caption{Results of  study 4, based on 1000
		simulations with sample size 200. ``bias'' is the absolute
		value of $\text{mean}(\wh\bb)-\bb$ of each component in $\bb$, ``sd'' is the
		sample standard errors of the corresponding estimation.
		The last column is the mean and standard errors of the largest singular value of
		$\wh\bP-\bP$. } \label{tab:simu4}
	\vspace{0.1cm}
	\begin{tabular}{cc|cccccccc|c}\hline
		&& $\beta_{3,1}$ &$\beta_{4,1}$ &$\beta_{5,1}$ &$\beta_{6,1}$
		&$\beta_{3,2}$      &$\beta_{4,2}$&$\beta_{5,2}$&$\beta_{6,2}$& $\lambda_{\max}$\\
		&true& 2.75 &-0.75&-1 &2.0 &-3.125&-1.125  & 1.0 &-2.0& \\\hline
		&&\multicolumn{8}{c}{No censoring}&\\
		hmave& bias  & 5.8883 & 3.4252 & 0.4704 & 2.9062 & 43.257 & 38.079 & 8.0696 & 33.551 &  0.8351   \\
		& sd   & 130.6 & 109.3 & 24.214 & 101.81 & 1091.9 & 939.06 & 174.9 & 855.8&  0.1116     \\
		semi	& bias     & 0.1609 & 0.0914 & 0.0609 & 0.1388 & 0.1257 & 0.1002 & 0.0372 & 0.0887 &   0.1791    \\
		& sd & 0.3393 & 0.2143 & 0.2545 & 0.2837 & 0.3674 & 0.2619 & 0.2559 & 0.3121 &    0.0788\\
		\hline
		&&\multicolumn{8}{c}{20\% censoring}&\\		
		hmave& bias   & 4.2011 & 2.8472 & 1.7685 & 3.4267 & 3.5963 & 0.4313 & 0.6443 & 0.9892 &   0.9273 \\
		& sd     & 64.499 & 39.786 & 56.389 & 41.570 & 20.496 & 24.564 & 23.323 & 22.339  &   0.1038  \\
		semi	& bias & 0.0846 & 0.0358 & 0.0311 & 0.0596 & 0.1133 & 0.0726 & 0.0305 & 0.0724 & 0.0971  \\
		& sd  & 0.3433 & 0.2068 & 0.2018 & 0.2646 & 0.4051 & 0.2690 & 0.2374 & 0.2952 &   0.0777\\
		\hline
		&&\multicolumn{8}{c}{40\% censoring}&\\		
		hmave& bias& 3.2328 & 0.5529 & 2.0826 & 1.0118 & 4.2734 & 3.6817 & 2.0900 & 2.9417 & 0.9363  \\
		& sd    & 14.712 & 19.789 & 17.661 & 19.938 & 26.892 & 73.550 & 29.838 & 31.115&  0.1071   \\
		semi	& bias  & 0.0986 & 0.0555 & 0.0246 & 0.0808 & 0.1420 & 0.0868 & 0.0451 & 0.0950&    0.0915  \\
		& sd    & 0.3604 & 0.2173 & 0.2168 & 0.2991 & 0.4864 & 0.2627 & 0.2645 & 0.3099 &  0.0898\\
		\hline 
	\end{tabular}
\end{table}

We also performed an additional experiment to 
further assess the finite sample performance of
the asymptotic results established in Section \ref{sec:asym}. To this
end, we generate covariates $\X$ from a standard normal 
distribution and event times $T$ from a distribution with hazard function
\bse
\lambda(t|\X) = \lambda _{0}(t)\left\{\sum_{j=1}^2\exp\left(\bb_j\trans\X\right)\right\},
\ese
where the baseline hazard $\lambda_{0}(t) = t$ and the dimension of $\bb$ is $d=2,p=6$.
We use the parameter values $\bb=\{(2.75,-0.75,-1,2)\trans;$
$(-3.125,-1.125,1,-2)\trans\}\trans$, and
adopt the same censoring process as in
the second study to yield 40\% censoring rate.
We carry out 1000 simulations and consider sample
sizes $n=100$, 500 and 1000.  
The estimation results, together with sample standard errors,
average of the estimated standard deviations and coverage
probabilities of the 95\% confidence intervals are given in Table
\ref{tab:simulast}. These
results indicate that the large sample properties of the estimator
requires more sample size than 1000. However, the general trend is
that when sample size increases, the results are approaching what we
expect based on the asymptotic results, in that the sample standard
errors and their estimated versions are becoming closer to each other,
and the 95\% coverage probabilities are getting closer to the nominal
level. The phenomenon that asymptotic result requires very large
sample size to illustrate itself is quite common in survival data
analysis and is not unique to our semiparametric method. Due to the
limited sample size in practice, we recommend to use bootstrap to
assess estimation variability. 

\begin{table}[H]
	\centering
	\caption{Results of  study 5, based on 1000
		simulations with sample size 100, 500, 1000 respectively. ``bias'' is the absolute
		value of $\text{mean}(\wh\bb)-\bb$ of each component in $\bb$, ``sd'' is the
		sample standard errors of the corresponding estimation, ``$\wh\sigma$''
		is the mean of the estimated standard errors of $\wh\bb$ component,
		``95\%'' is the sample coverage of the 95\% confidence intervals.  
	}\label{tab:simulast}
	\vspace{0.1cm}
	\begin{tabular}{c|cccccccc}\hline
		& $\beta_{31}$ & $\beta_{41}$ & $\beta_{51}$ & $\beta_{61}$  & $\beta_{32}$ & $\beta_{42}$ & $\beta_{52}$ & $\beta_{62}$ \\
		& 2.75 & -0.75 & -1 & 2.0 & -3.125 & -1.125 & 1.0 & -2.0 \\\hline
		&\multicolumn{8}{c}{$n=100$}\\
		bias  & 0.3995 & 0.5031 & 0.2066 & 0.3799 & 0.5515 & 0.5349 & 0.1757 & 0.3395 \\
		sd  & 0.5760 & 0.4236 & 0.4673 & 0.5608 & 0.6163 & 0.4376 & 0.4772 & 0.5377 \\
		$\wh {\sigma}$ & 0.3868 & 0.3188 & 0.3312 & 0.3427 & 0.3956 & 0.3131 & 0.3331 & 0.3602 \\
		95\%  & 0.7100 & 0.6577 & 0.8051 & 0.7034 & 0.6416 & 0.6292 & 0.8089 & 0.7414 \\
		\hline
		&\multicolumn{8}{c}{$n=500$}\\
		bias   & 0.1790 & 0.1258 & 0.0714 & 0.1338 & 0.2100 & 0.1489 & 0.07386 & 0.1340 \\
		sd  & 0.2741 & 0.1714 & 0.2177 & 0.2380 & 0.2979 & 0.1897 & 0.2202 & 0.2244 \\
		$\wh {\sigma}$ & 0.1585 & 0.1371 & 0.1644 & 0.1659 & 0.2683 & 0.2179 & 0.2538 & 0.2558 \\
		95\% & 0.6663 & 0.8022 & 0.8298 & 0.7566 & 0.8127 & 0.8773 & 0.9125 & 0.8764 \\
		\hline
		&\multicolumn{8}{c}{$n=1000$}\\
		bias  & 0.0611 & 0.0492 & 0.0188 & 0.0423 & 0.0695 & 0.0467 & 0.0209 & 0.0448 \\
		sd  & 0.1951 & 0.1451 & 0.1555 & 0.1538 & 0.1867 & 0.1433 & 0.1650 & 0.1711 \\
		$\wh {\sigma}$ & 0.1062 & 0.1113 & 0.1134 & 0.1190 & 0.1823 & 0.1712 & 0.1749 & 0.1740 \\
		95\% & 0.8060 & 0.8830 & 0.8783 & 0.8621 & 0.9268 & 0.9705 & 0.9515 & 0.9287 \\
		\hline
		
	\end{tabular}
	
\end{table}

\subsection{AIDS Application}

We apply the proposed method to analyze the HIV data from AIDS
Clinical Trials Group Protocol 175 (ACTG175)
(\cite{hammer1996trial}). In this study, 2137 HIV-infected subjects
were randomized to receive one of  four treatments: zidovudine
(ZDV) monotherapy, ZDV plus didanosine, ZDV plus zalcitabine and ddI
monotherapy. As in \cite{Gengetal2015} and \cite{Jiangetal2016}, the
survival time of interest was chosen as the time to having a larger
than $50\%$ decline in the CD4 count, or progressing to AIDS or death,
whichever comes first. Besides the treatments, there were 12 covariates included in our
study, specifically, patient age 
in years at baseline ($X_1$), patient weight in kilogram at baseline
($X_2$), hemophilia indicator ($X_3$), homosexual activity ($X_4$),
history of IV drug use ($X_5$), Karnofsky score on a scale of 0-100
($X_6$), race ($X_7$), gender ($X_8$), antiretroviral history ($X_9$),
symptomatic indicator 
($X_{10}$),  number of CD4 at baseline ($X_{11}$), number of CD8 at
baseline ($X_{12}$), treatment indicator ($X_{13}$), 
where we coded $X_{13}=0$ for treatment ZDV+ddl and
$X_{13}=1$ for treatment ZDV+Zal.
As in \cite{Jiangetal2016}, we only analyzed data from the two
composite treatments: ZDV plus didanosine and ZDV plus zalcitabine,
which had been shown to have significantly better survival than the
other two treatments \citep{Gengetal2015}. This subset of data
contains 1046 subjects with the censoring rate around 75\%. In
addition, each covariate is standardized respectively with no obvious
outliers and no missing values.

To determine the proper reduced space dimension $d$, we employ the
Validated Information Criterion (VIC) \citep{ma2015validated}, where
the $d$ corresponding to the smallest VIC value is selected. In the
example,  the VIC value at $d=1$ is 90.38. Further, when $d\ge 2$, the
VIC values are all greater than 180.7 which result from the penalty
term alone. Hence we choose $d=1$ as the final model. 
Table \ref{tab:apply} contains the estimated coefficient $\wh
\bb$'s under the selected model, with the corresponding estimation
standard errors and $p$-values. Here, we implemented the
semiparametric estimator to obtain these results due to its superior
theoretical and numerical performance. 

\begin{figure}[H]
	\centering
	\includegraphics[trim={1cm 0 0 0},height=8cm,width=8cm]{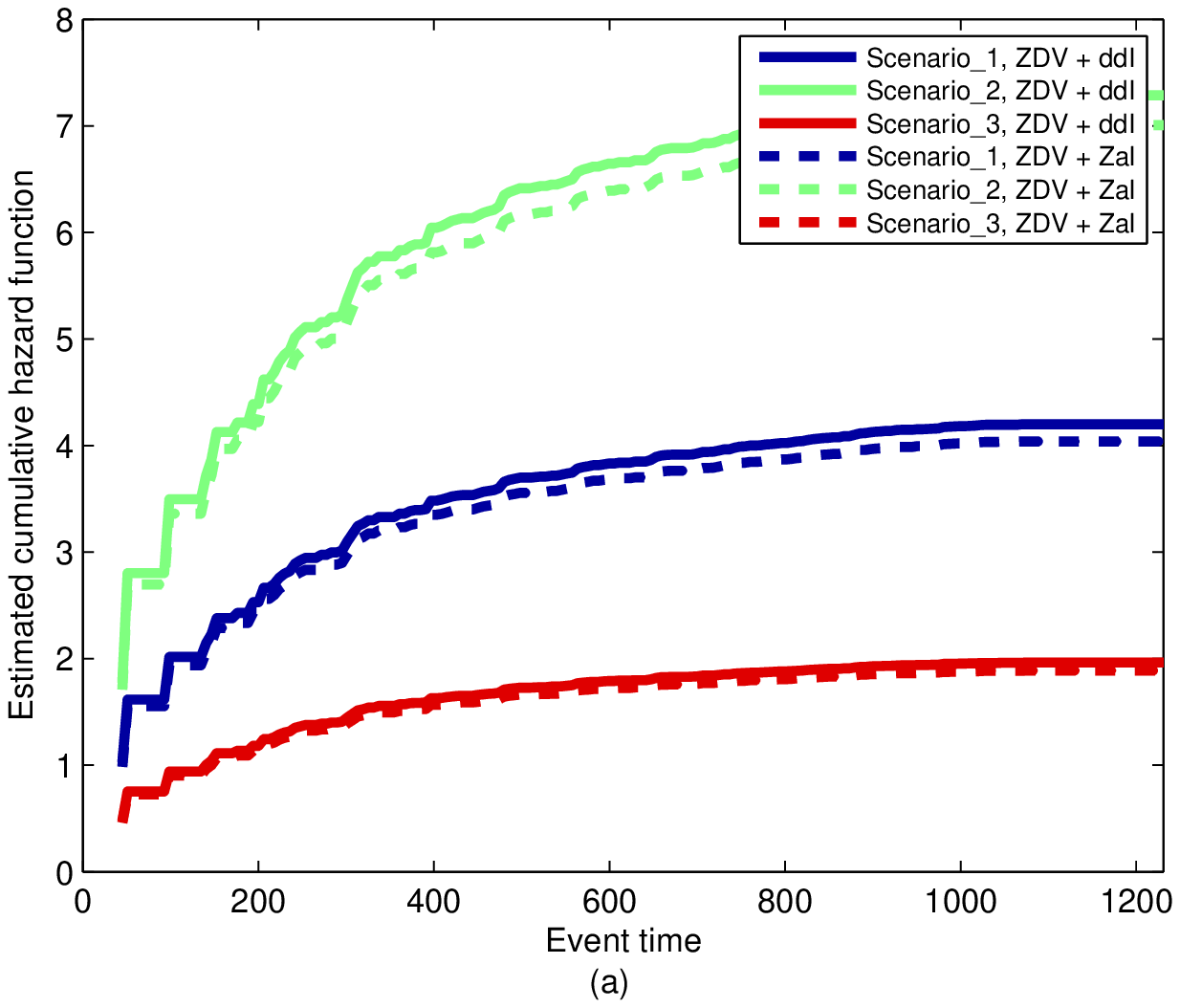}
	\includegraphics[trim={1cm 0 0 0},height=8cm,width=8cm]{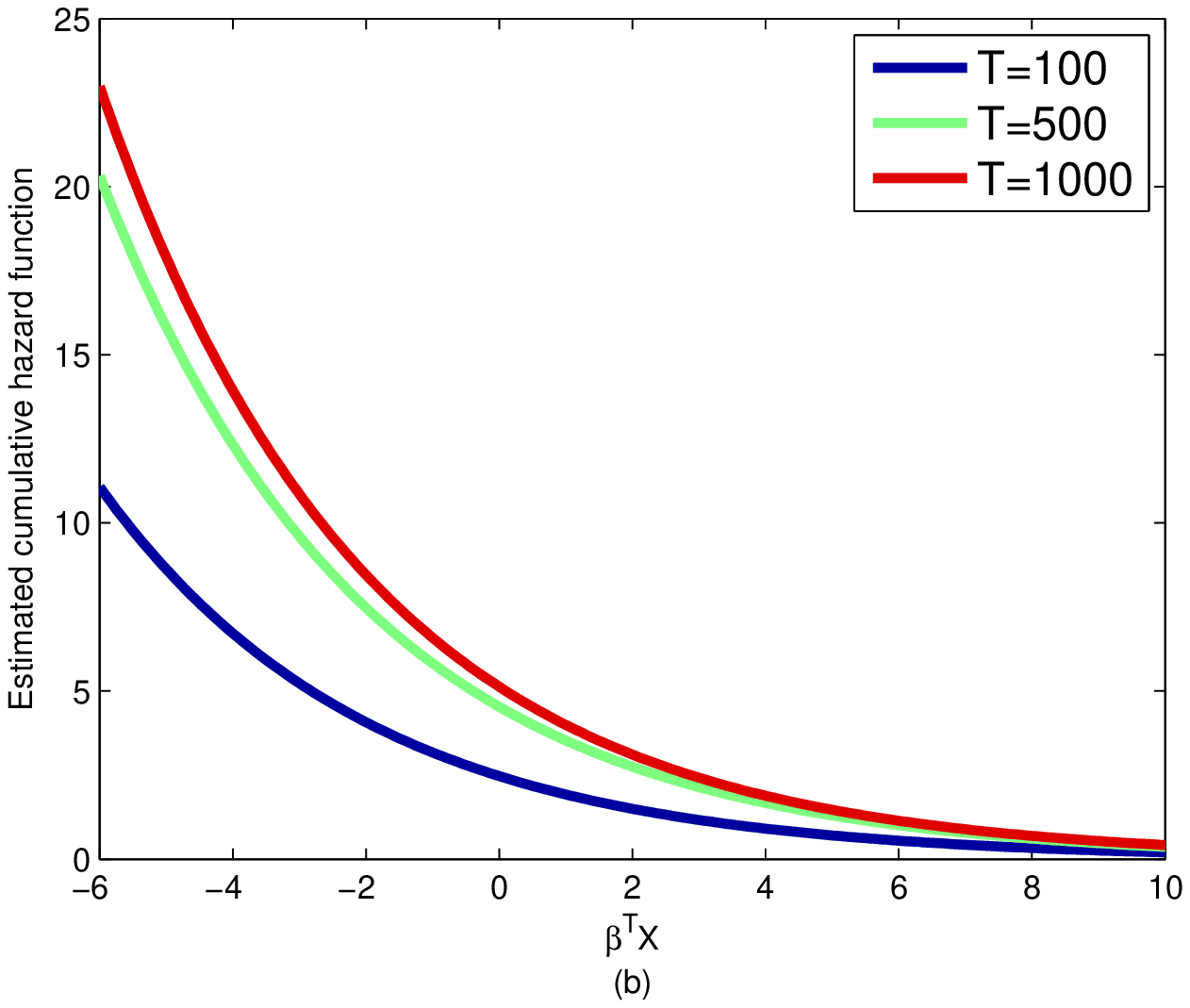}\\
	\includegraphics[trim={1cm 0 0 0},height=8cm,width=10cm]{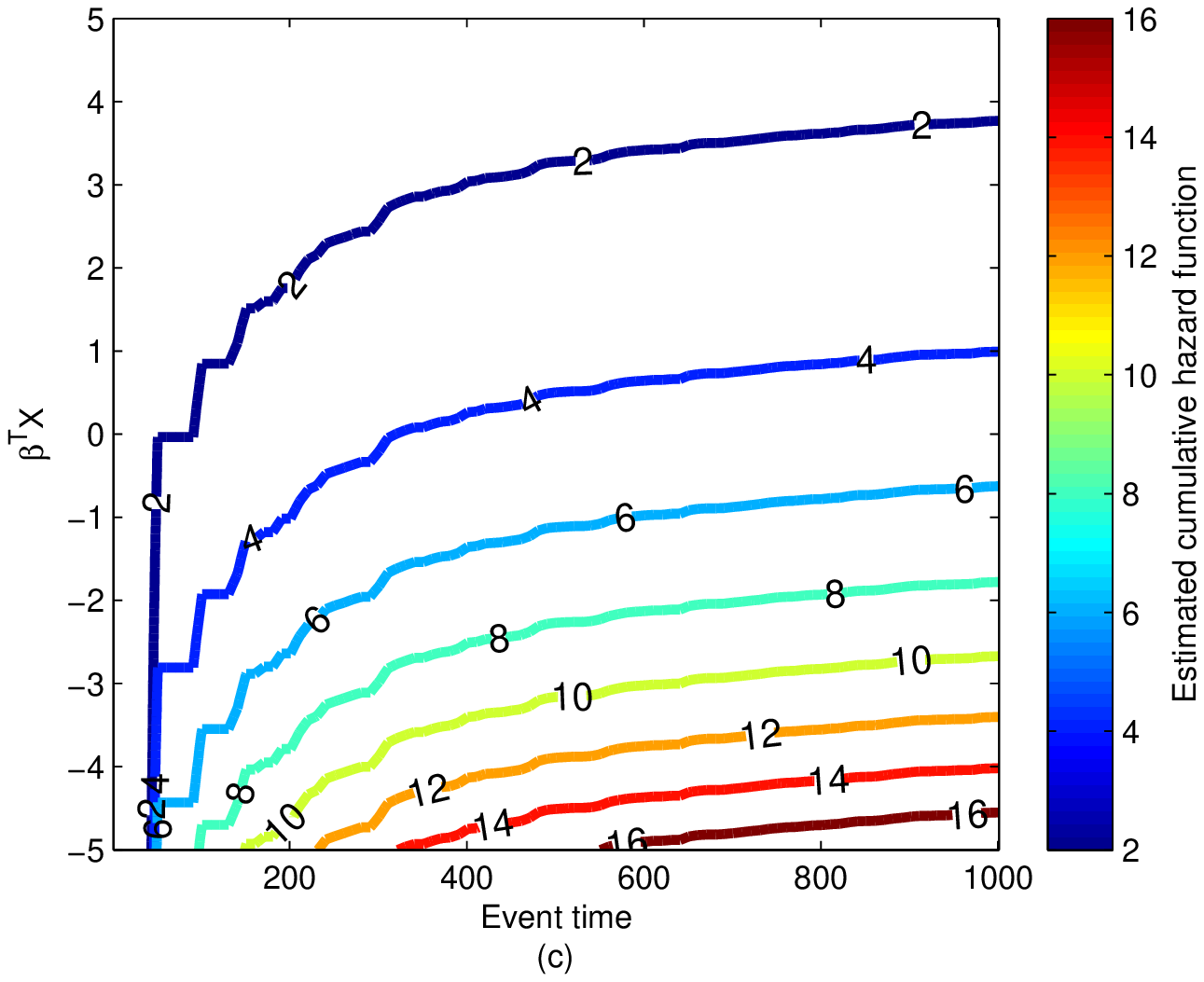}
	\caption{{Estimated cumulative hazard function $\wh \Lambda$ in
			AIDS data. (a). Comparisons of $\wh \Lambda$ as a function
			of $t$ between treatments ZDV+ddl and ZDV+Zal when other
			covariates are fixed at three indices. (b). $\wh \Lambda$ as
			a function of $\wh\bb\trans\X$ at $T=100,500,1000$. (c). Contour
			plot of $\wh\Lambda$ as a function of $T$ and $\wh\bb\trans\X$.}}
	\label{fig:application}
\end{figure}

The results in Table \ref{tab:apply}
indicate that in forming  the index described by $\wh\bb_{\cdot,1}$,
all covariates are significant except 
hemophilia indicator ($X_3$), gender ($X_8$) and number of CD4 at
baseline ($X_{11}$).  
The estimated cumulative hazard functions are also reported in Figure
\ref{fig:application}, where it is plotted as a function of time
(upper left panel), a
function of the covariate index $\bb\trans\x$ (upper right panel)  and as a function of
both (bottom panel).  Specifically, in plotting the cumulative hazard
as a function of time $t$, we fixed the covariate index at three 
different sets of covariate values, respectively 
$X_{1:12}=(40,60,1,0,0,80,0,0,0,1,200,800)\trans$,
$X_{1:12}=(20,70,1,0,0,80,0,1,0,1,200,800)\trans$,
$X_{1:12}=(60,70,1,0,0,20,0,0,0,1,200,200)\trans$, in combination with
the treatment indicator of both $X_{13}=0$ and $X_{13}=1$.
Based on the plots, the 
estimated cumulative hazard of treatment ZDV+ddl is slightly larger
than that of treatment ZDV+Zal in all scenarios. 
In plotting the
estimated cumulative hazard $\wh\Lambda$ as a function of the
index $\wh\bb\trans\x$, we fixed the time at
$t=100,500$ and 1000.
Finally, we also plotted the cumulative hazard as a function of two
variables $t$ and $\bb\trans\x$ using the contour plot, where the
hazard values are explicitly written out on each contour. 

\begin{table}[H]
	\centering
	
	\caption{The estimated coefficients, standard errors and
		p-value in AIDS data.
	}\label{tab:apply}
	\vspace{0.1cm}
	\begin{tabular}{c|cccccccccccc}\hline
		&&&&&&&&&&&&\\[-0.8em]
		& $\wh\bb_{2,1}$ & $\wh\bb_{3,1}$ & $\wh\bb_{4,1}$ & $\wh\bb_{5,1}$ & $\wh\bb_{6,1}$ & $\wh\bb_{7,1}$ & $\wh\bb_{8,1}$ & $\wh\bb_{9,1}$ & $\wh\bb_{10,1}$ & $\wh\bb_{11,1}$ & $\wh\bb_{12,1}$& $\wh\bb_{13,1}$\\
		&&&&&&&&&&&&\\[-0.8em]\hline
		&&&&&&&&&&&&\\[-0.5em]
		est&     0.115&    -0.002&  0.093&  0.088& -0.090& 0.231&      -0.003& -0.178&       0.058& -0.031& 0.201&     0.156 \\
		std&   0.039  &  0.039 &   0.039 &   0.037 &   0.043  &  0.046   & 0.036 &   0.046  &  0.035  &  0.042  &  0.033  &  0.038 \\
		$p$-value &        0.003&  0.965&        0.017&   0.017& 0.036&   0.001&  0.928&       0.001& 0.100&       0.457&      0.001&  0.001
		\\\hline
	\end{tabular}
	
\end{table}

\section{Discussion}\label{sec:discuss}

We have considered a very general model for analyzing time to event
data subject to censoring. The model allows the event times to link to
the covariate indices in an unspecified fashion. Because both the
number of indices and the functional form of the linkage to the
indices are data determined, conceptually the model is maximally
flexible. In practice, relatively low number of indices are expected
to avoid curse of dimensionality. The work is conducted without requiring covariate
independent censoring. Instead, it only requires event
independent censoring conditional on covariates, which is the minimum
requirement for identification. We derived a class of estimators
which are consistent and asymptotically normal. We also proposed a procedure
to construct the semiparametric efficient estimator that achieves the
optimal estimation variability among all possible consistent
estimators. 

\section*{Appendix}
\renewcommand{\theequation}{A.\arabic{equation}}
\renewcommand{\thesubsection}{A.\arabic{subsection}}
\setcounter{equation}{0}

\subsection{Proof of Proposition \ref{pro:nuisance}}\label{sec:pronuisance}
The result of $\Lambda_1$ is obvious. 
To obtain $\Lambda_2$,
let
$\h(t,\bb_0\trans\X,\bg)=\partial
\log\lambda(t,\bb_0\trans\X,\bg)/\partial\bg$,
where $\lambda(t,\bb_0\trans\X,\bg)$ is a submodel of $\lambda(t,\bb_0\trans\X)$.
Hence
\bse
\frac{\partial\log f(\X,Z,\Delta)}{\partial \bg}
&=&\Delta\frac{\partial\log\lambda(Z,\bb_0\trans\X,\bg)}{\bg}-\int_0^Z\frac{\partial
	\lambda(s,\bb_0\trans\X,\bg)}{\partial\bg} ds\\
&=&\Delta\h(Z,\bb_0\trans\X,\bg)-\int_0^Z\h(s,\bb_0\trans\X,\bg)\lambda(s,\bb_0\trans\X)ds\\
&=&\int_0^\infty\h(s,\bb_0\trans\X,\bg)dM(s,\bb_0\trans\X).
\ese
Because
$\lambda_0(t,\bb_0\trans\X)$ can be any positive function,
$\h(s,\bb_0\trans\X,\bg)$ can be any function. We denote it
$\h(s,\bb_0\trans\X)$. This leads to the form of $\Lambda_2$.

Similar derivation leads to
$\Lambda_3$.

It is easy to verify that $\Lambda_1\perp\Lambda_2$ and 
$\Lambda_1\perp\Lambda_3$. Because $C\indep T\mid\X$, the 
martingale integrations associated with $M(t,\bb_0\trans\X)$ and
$M_C(z,\X)$ are also independent conditional on $\X$, hence
$\Lambda_2\perp\Lambda_3$. 
This completes the proof. 
\qed

\subsection{Proof of Lemma \ref{lem:pre}}\label{sec:prooflempre}
For notational convenience, we prove the results for $d=1$ and assume
the first component of $\bb$ is 1. 
The first four results are obtained from the uniform convergence
property of the kernel estimation \citep{macksilverman1982} under conditions \ref{assum:kernel}-\ref{assum:bandwidth}. 
Specifically, to derive the first four results, we first establish the following
preliminary conclusion for any $Z$ and $\bb\trans\X$,
\be
\frac{1}{n}\sumj K_h(\bb\trans\X_j-\bb\trans\x)=f_{\bb\trans\X}(\bb\trans\x)+O_p(n^{-1/2}h^{-1/2}+h^2).\label{eq:fbeta}
\ee
To see this, we compute the absolute bias of the left hand side of
(\ref{eq:fbeta}) as
\bse
&&\left|E\left\{\frac{1}{n}\sumj K_h(\bb\trans\X_j-\bb\trans\x)\right\}-f_{\bb\trans\X}(\bb\trans\x)\right|\\
&=&\left|EK_h(\bb\trans\X_j-\bb\trans\x)-f_{\bb\trans\X}(\bb\trans\x)\right|\\
&=&\left|\int\frac{1}{h}K\left(\frac{\bb\trans\x_j-\bb\trans\x}{h}\right)f_{\bb\trans\X}(\bb\trans\x_j)d\bb\trans\x_j -f_{\bb\trans\X}(\bb\trans\x)\right|\\
&=&\left|\int K(u)f_{\bb\trans\X}(\bb\trans\x+hu)du -f_{\bb\trans\X}(\bb\trans\x)\right|\\
&=&\left|\int K(u)\left\{f_{\bb\trans\X}(\bb\trans\x)+f_{\bb\trans\X}'(\bb\trans\x)hu+\frac{1}{2}f''_{\bb\trans\X}(\bb\trans\x^*)h^2u^2\right\}du -f_{\bb\trans\X}(\bb\trans\x)\right|\\
&\le&\frac{h^2}{2}\sup_{\bb\trans\x} |f''_{\bb\trans\X}(\bb\trans\x)|\int u^2K(u)du,
\ese
where throughout the text, $\bb\trans\x^*$ is on the line connecting $\bb\trans\x$ and $\bb\trans\x+hu$,
and the variance to be
\bse
&&\var\left\{\frac{1}{n}\sumj K_h(\bb\trans\X_j-\bb\trans\x)\right\}\\
&=&\frac{1}{n}\var K_h(\bb\trans\X_j-\bb\trans\x)\\
&=&\frac{1}{n}\left[EK_h^2(\bb\trans\X_j-\bb\trans\x)-\left\{EK_h(\bb\trans\X_j-\bb\trans\x)\right\}^2\right]\\
&=&\frac{1}{n}\left[\int\frac{1}{h^2}K^2\{(\bb\trans\x_j-\bb\trans\x)/h\}f_{\bb\trans\x}(\bb\trans\x_j)d\bb\trans\x_j-f_{\bb\trans\x}^2(\bb\trans\x)+O(h^2)\right]\\
&=&\frac{1}{nh}\int
K^2(u)f_{\bb\trans\x}(\bb\trans\x+hu)du-\frac{1}{n}f_{\bb\trans\x}^2(\bb\trans\x)+O(h^2/n)
\\
&\le&\frac{1}{nh}f_{\bb\trans\x}(\bb\trans\x)\int K^2(u)du
+\frac{h}{2n}\sup_{\bb\trans\x}|f_{\bb\trans\x}''(\bb\trans\x)|\int u^2K^2(u)du
+\frac{1}{n}|f_{\bb\trans\x}^2(\bb\trans\x)|+O(h^2/n).
\ese
Therefore, we have
\bse
\frac{1}{n}\sumj K_h(\bb\trans\X_j-\bb\trans\x)=f_{\bb\trans\X}(\bb\trans\x)+O_p(n^{-1/2}h^{-1/2}+h^2)
\ese
uniformly under conditions \ref{assum:kernel}, \ref{assum:bandwidth} and \ref{assum:fbeta}.

Following similar derivations, we have
\bse
&&E\left\{-\frac{1}{n}\sumj \X_j I(Z_j\ge Z)  K_h'(\bb\trans\X_j-\bb\trans\x)\right\}\\
&=&E\left\{-\X_j I(Z_j\ge Z)  K_h'(\bb\trans\X_j-\bb\trans\x)\right\}\\
&=&-\int\frac{1}{h^2}\x_j I(c\ge Z)I(t\ge Z)  K'\left\{(\bb\trans\x_j-\bb\trans\x)/h\right\}\\
&&\times f_c(c,\x_j)f(t,\bb\trans\x_j)f_{\X\mid\bb\trans\X}(\x_j\mid\bb\trans\x_j)f_{\bb\trans\X}(\bb\trans\x_j)dcdtd\x_j d\bb\trans\x_j\\
&=&-\int\frac{1}{h^2}\x_j K'\left\{(\bb\trans\x_j-\bb\trans\x)/h\right\}
S_c(Z,\x_{jl},\bb\trans\x_j)S(Z,\bb\trans\x_j)f_{\X_{l}\mid\bb\trans\X}(\x_{jl}\mid\bb\trans\x_j)\\
&&\times f_{\bb\trans\X}(\bb\trans\x_j)d\x_{jl} d\bb\trans\x_j\\
&=&-\frac{1}{h}\int (\bb\trans\x+hu,\x_{jl}\trans)\trans
K'(u)S_c(Z,\x_{jl},\bb\trans\x+hu)S(Z,\bb\trans\x+hu)f_{\X_l\mid\bb\trans\X}(\x_{jl}\mid\bb\trans\x+hu)\\
&&\times f_{\bb\trans\X}(\bb\trans\x+hu)d\x_{jl} du\\
&=&-\frac{1}{h}\int(\bb\trans\x,\x_{jl}\trans)\trans
S_c(Z,\x_{jl},\bb\trans\x)S(Z,\bb\trans\x)f_{\X_l\mid\bb\trans\X}(\x_{jl}\mid\bb\trans\x)f_{\bb\trans\X}(\bb\trans\x)d\x_{jl}
\int K'(u)du\\
&&-\frac{1}{h}\int\frac{\partial}{\partial\bb\trans\x} \left\{(\bb\trans\x,\x_{jl}\trans)\trans S_c(Z,\x_{jl},\bb\trans\x)S(Z,\bb\trans\x)f_{\X_l\mid\bb\trans\X}(\x_{jl}\mid\bb\trans\x)f_{\bb\trans\X}(\bb\trans\x)\right\}d\x_{jl} \\
&&\times hu K'(u)du\\
&&-\frac{1}{2h}\int\frac{\partial^2}{\partial(\bb\trans\x)^2}
\left\{(\bb\trans\x,\x_{jl}\trans)\trans
S_c(Z,\x_{jl},\bb\trans\x)S(Z,\bb\trans\x)f_{\X_l\mid\bb\trans\X}(\x_{jl}\mid\bb\trans\x)f_{\bb\trans\X}(\bb\trans\x)\right\}d\x_{jl}\\
&&\times \int h^2u^2 K'(u)du\\
&&-\frac{1}{6h}\iint \frac{\partial^3}{\partial(\bb\trans\x)^3}
\left\{(\bb\trans\x,\x_{jl}\trans)\trans
S_c(Z,\x_{jl},\bb\trans\x^*)S(Z,\bb\trans\x^*)f_{\X_l\mid\bb\trans\X}(\x_{jl}\mid\bb\trans\x^*)f_{\bb\trans\X}(\bb\trans\x^*)\right\}d\x_{jl}
\\
&&\times h^3u^3 K'(u)du\\
&=&
\int\frac{\partial}{\partial\bb\trans\x} \left\{(\bb\trans\x,\x_{jl}\trans)\trans S_c(Z,\x_{jl},\bb\trans\x)S(Z,\bb\trans\x)f_{\X_l\mid\bb\trans\X}(\x_{jl}\mid\bb\trans\x)f_{\bb\trans\X}(\bb\trans\x)\right\}d\x_{jl}\\
&&-\frac{h^2}{6}\int \int\frac{\partial^3}{\partial(\bb\trans\x)^3}
\left\{(\bb\trans\x,\x_{jl}\trans)\trans
S_c(Z,\x_{jl},\bb\trans\x^*)S(Z,\bb\trans\x^*)f_{\X_l\mid\bb\trans\X}(\x_{jl}\mid\bb\trans\x^*)f_{\bb\trans\X}(\bb\trans\x^*)\right\}d\x_{jl}
\\
&&\times u^3 K'(u)du\\
&=&
\frac{\partial}{\partial \bb\trans\x}
S(Z,\bb\trans\x)f_{\bb\trans\X}(\bb\trans\x)E\{\X_j
S_c(Z,\X_j)\mid\bb\trans\x\}\\
&&-\frac{h^2}{6}\int
\frac{\partial^3}{\partial(\bb\trans\x)^3}
S(Z,\bb\trans\x^*)f_{\bb\trans\X}(\bb\trans\x^*) E
\left\{\X_j
S_c(Z,\X_j)\mid\bb\trans\x^*\right\}
u^3 K'(u)du\\
&=&
\frac{\partial}{\partial \bb\trans\x}
f_{\bb\trans\X}(\bb\trans\x)E\{\X_jI(Z_j\ge Z)\mid\bb\trans\x\}\\
&&-\frac{h^2}{6}\int
\frac{\partial^3}{\partial(\bb\trans\x)^3}
f_{\bb\trans\X}(\bb\trans\x^*) E
\left\{\X_jI(Z_j\ge Z)\mid\bb\trans\x^*\right\}
u^3 K'(u)du.
\ese
Hence the absolute bias is
\bse
&&\left|
E\left\{-\frac{1}{n}\sumj \X_j I(Z_j\ge Z)
K_h'(\bb\trans\X_j-\bb\trans\x)\right\}
-\frac{\partial}{\partial \bb\trans\x}
f_{\bb\trans\X}(\bb\trans\x)E\{\X_j(Z_j\ge
Z)\mid\bb\trans\x\}\right|\\
&=&	\left|-\frac{h^2}{6}\int
\frac{\partial^3}{\partial(\bb\trans\x)^3}
f_{\bb\trans\X}(\bb\trans\x^*) E
\left\{\X_jI(Z_j\ge Z)\mid\bb\trans\x^*\right\}
u^3 K'(u)du\right|\\
&\le&\frac{h^2}{6}\sup_{\bb\trans\x}\left|
\frac{\partial^3}{\partial(\bb\trans\x)^3}
f_{\bb\trans\X}(\bb\trans\x) E
\left\{\X_jI(Z_j\ge Z)\mid\bb\trans\x\right\}\right|
\left\{3\int u^2 K(u)du\right\}\\
\ese
and 
\bse
&&\var\left\{-\frac{1}{n}\sumj \X_j I(Z_j\ge Z)  K_h'(\bb\trans\X_j-\bb\trans\x)\right\}\\
&=&\frac{1}{n}\left[E\left\{\X_jI(Z_j\ge Z)  K_h'(\bb\trans\X_j-\bb\trans\x)\left\}\right\{\X_jI(Z_j\ge Z)  K_h'(\bb\trans\X_j-\bb\trans\x)\right\}\trans\right.\\
&&\left.-\left\{E\X_j I(Z_j\ge Z)  K_h'(\bb\trans\X_j-\bb\trans\x)\right\}\left\{E\X_j I(Z_j\ge Z)  K_h'(\bb\trans\X_j-\bb\trans\x)\right\}\trans\right]\\
&=&\frac{1}{n}\left(\int\frac{1}{h^4}\x_j\x_j\trans I(Z_j\ge Z)K'^2\left\{{(\bb\trans\x_j-\bb\trans\x)}/{h}\right\}\right.\\
&&\left.\times f_c(c,\x_j)f(t,\bb\trans\x_j)f_{\X_l\mid\bb\trans\X}(\x_{jl}\mid\bb\trans\x_j)f_{\bb\trans\X}(\bb\trans\x_j)dcdtd\x_{jl} d\bb\trans\x_j\right.\\
&&-\left.\left[\frac{\partial}{\partial \bb\trans\x}
f_{\bb\trans\X}(\bb\trans\x)E\{\X_j
I(Z_j\ge
Z)\mid\bb\trans\x\}\right]\left[\frac{\partial}{\partial
	\bb\trans\x}f_{\bb\trans\X}(\bb\trans\x)E\{\X_j I(Z_j\ge Z)\mid\bb\trans\x\}\right]\trans+O(h^2)\right)\\
&=&\frac{1}{n}\int\frac{1}{h^3}(\bb\trans\x+hu,\x_{jl}\trans)\trans(\bb\trans\x+hu,\x_{jl}\trans) S_c(Z,\x_{jl},\bb\trans\x+hu)\\
&&\times S(Z,\bb\trans\x+hu)f_{\X_l\mid\bb\trans\X}(\x_{jl}\mid\bb\trans\x+hu)f_{\bb\trans\X}(\bb\trans\x+hu)d\x_j K'^2(u)du\\
&&+O(1/n)\\
&=&\frac{1}{nh^3}\int(\bb\trans\x,\x_{jl}\trans)\trans(\bb\trans\x,\x_{jl}\trans)S_c(Z,\x_{jl},\bb\trans\x)S(Z,\bb\trans\x)f_{\bb\trans\X}(\bb\trans\x)f_{\X_l\mid\bb\trans\X}(\x_{jl}\mid\bb\trans\x)d\x_{jl}
\\
&&\times\int K'^2(u)du\\
&&+\frac{1}{2nh^3}
\frac{\partial^2}{\partial(\bb\trans\x)^2}
\iint(\bb\trans\x^*,\x_{jl}\trans)\trans(\bb\trans\x^*,\x_{jl}\trans)S_c(Z,\x_{jl},\bb\trans\x^*)S(Z,\bb\trans\x^*)\\
&&\times
f_{\bb\trans\X}(\bb\trans\x^*)f_{\X_l\mid\bb\trans\X}(\x_{jl}\mid\bb\trans\x^*)d\x_{jl}
h^2u^2K'^2(u)du+O(1/n)\\
&\le&\frac{1}{nh^3} \sup_{\bb\trans\x}\left|f_{\bb\trans\X}(\bb\trans\x)E\{\X_j\X_j\trans
I(Z_j\ge Z)\mid\bb\trans\x\}\right|\int K'^2(u)du\\
&&+\frac{1}{2nh} \sup_{\bb\trans\x^*}\left|\frac{\partial^2}{\partial(\bb\trans\x)^2}
f_{\bb\trans\X}(\bb\trans\x^*)E\{\X_j\X_j\trans
I(Z_j\ge Z)\mid\bb\trans\x^*\}\right|\int u^2K'^2(u)du+O(1/n).
\ese
So
\be\label{eq:fbetaxi1}
&&\frac{1}{n}\sumj \X_j I(Z_j\ge Z)  K_h'(\bb\trans\X_j-\bb\trans\x)\n\\
&=&-\frac{\partial}{\partial \bb\trans\x} 
f_{\bb\trans\X}(\bb\trans\x)E\{\X_j I(Z_j\ge Z)\mid\bb\trans\x\}+O_p(n^{-1/2}h^{-3/2}+h^2)
\ee
uniformly under conditions \ref{assum:kernel}-\ref{assum:bandwidth}.

We next show
\be\label{eq:fbeta1}
-\frac{1}{n}\sumj K_h'(\bb\trans\X_j-\bb\trans\x)=f_{\bb\trans\X}'(\bb\trans\x)+O_p(n^{-1/2}h^{-3/2}+h^2).
\ee
To see this, we compute the absolute bias of the left hand side as
\bse
&&\left|E\left\{-\frac{1}{n}\sumj K_h'(\bb\trans\X_j-\bb\trans\x)\right\}-f_{\bb\trans\X}(\bb\trans\x)\right|\\
&=&\left|E\left\{-K_h'(\bb\trans\X_j-\bb\trans\x)\right\}-f_{\bb\trans\X}(\bb\trans\x)\right|\\
&=&\left|\int-\frac{1}{h^2}K'\left(\frac{\bb\trans\x_j-\bb\trans\x}{h}\right)f_{\bb\trans\X}(\bb\trans\x_j)d\bb\trans\x_j -f_{\bb\trans\X}(\bb\trans\x)\right|\\
&=&\left|\int -\frac{1}{h}K'(u)f_{\bb\trans\X}(\bb\trans\x+hu)du -f_{\bb\trans\X}(\bb\trans\x)\right|\\
&=&\left|\int-\frac{1}{h}
K'(u)\left\{f_{\bb\trans\X}(\bb\trans\x)+f_{\bb\trans\X}'(\bb\trans\x)hu+\frac{1}{2}f''_{\bb\trans\X}(\bb\trans\x)h^2u^2+\frac{1}{6}f^{(3)}_{\bb\trans\X}(\bb\trans\x^*)h^3u^3\right\}du\right.\\
&&\left. -f_{\bb\trans\X}(\bb\trans\x)\right|\\
&=&\left|f_{\bb\trans\X}'(\bb\trans\x)-\frac{h^2}{6}\int f^{(3)}_{\bb\trans\X}(\bb\trans\x^*)u^3K'(u)du -f_{\bb\trans\X}(\bb\trans\x)\right|\\
&\le&\frac{h^2}{2}\sup_{\bb\trans\x} |f^{(3)}_{\bb\trans\X}(\bb\trans\x)|\int u^2K(u)du,
\ese
and the variance to be
\bse
&&\var\left\{\frac{1}{n}\sumj K_h'(\bb\trans\X_j-\bb\trans\x)\right\}\\
&=&\frac{1}{n}\var K_h'(\bb\trans\X_j-\bb\trans\x)\\
&=&\frac{1}{n}\left[EK_h'^2(\bb\trans\X_j-\bb\trans\x)-\left\{EK_h'(\bb\trans\X_j-\bb\trans\x)\right\}^2\right]\\
&=&\frac{1}{n}\left[\int\frac{1}{h^4}K'^2\{(\bb\trans\x_j-\bb\trans\x)/h\}f_{\bb\trans\x}(\bb\trans\x_j)d\bb\trans\x_j-f_{\bb\trans\x}^2(\bb\trans\x)+O(h^2)\right]\\
&=&\frac{1}{nh^3}\int
K'^2(u)f_{\bb\trans\x}(\bb\trans\x+hu)du-\frac{1}{n}f_{\bb\trans\x}^2(\bb\trans\x)+O(h^2/n)
\\
&\le&\frac{1}{nh^3}f_{\bb\trans\x}(\bb\trans\x)\int K'^2(u)du
+\frac{1}{2nh}\sup_{\bb\trans\x}|f_{\bb\trans\x}''(\bb\trans\x)|\int u^2K'^2(u)du
+O(1/n).
\ese
Therefore, (\ref{eq:fbeta1}) holds
uniformly under conditions \ref{assum:kernel}, \ref{assum:bandwidth} and \ref{assum:fbeta}.

We compute the expectation 
\bse
&&E\left\{-\frac{1}{n}\sumj I(Z_j\ge Z)  K_h'(\bb\trans\X_j-\bb\trans\x)\right\}\\
&=&E\left\{- I(Z_j\ge Z)  K_h'(\bb\trans\X_j-\bb\trans\x)\right\}\\
&=&-\int\frac{1}{h^2} I(c\ge Z)I(t\ge Z)  K'\left\{(\bb\trans\x_j-\bb\trans\x)/h\right\}\\
&&\times f_c(c,\x_j)f(t,\bb\trans\x_j)f_{\X\mid\bb\trans\X}(\x_j\mid\bb\trans\x_j)f_{\bb\trans\X}(\bb\trans\x_j)dcdtd\x_j d\bb\trans\x_j\\
&=&-\int\frac{1}{h^2} K'\left\{(\bb\trans\x_j-\bb\trans\x)/h\right\}
S_c(Z,\x_{jl},\bb\trans\x_j)S(Z,\bb\trans\x_j)f_{\X_{l}\mid\bb\trans\X}(\x_{jl}\mid\bb\trans\x_j)\\
&&\times f_{\bb\trans\X}(\bb\trans\x_j)d\x_{jl} d\bb\trans\x_j\\
&=&-\frac{1}{h}\int K'(u)S_c(Z,\x_{jl},\bb\trans\x+hu)S(Z,\bb\trans\x+hu)f_{\X_l\mid\bb\trans\X}(\x_{jl}\mid\bb\trans\x+hu)\\
&&\times f_{\bb\trans\X}(\bb\trans\x+hu) d\x_{jl}du\\
&=&-\frac{1}{h}\int S_c(Z,\x_{jl},\bb\trans\x)S(Z,\bb\trans\x)f_{\X_l\mid\bb\trans\X}(\x_{jl}\mid\bb\trans\x)f_{\bb\trans\X}(\bb\trans\x)d\x_{jl}
\int K'(u)du\\
&&-\frac{1}{h}\int\frac{\partial}{\partial\bb\trans\x} \left\{ S_c(Z,\x_{jl},\bb\trans\x)S(Z,\bb\trans\x)f_{\X_l\mid\bb\trans\X}(\x_{jl}\mid\bb\trans\x)f_{\bb\trans\X}(\bb\trans\x)\right\}d\x_{jl} \\
&&\times hu K'(u)du\\
&&-\frac{1}{2h}\int\frac{\partial^2}{\partial(\bb\trans\x)^2}\left\{S_c(Z,\x_{jl},\bb\trans\x)S(Z,\bb\trans\x)f_{\X_l\mid\bb\trans\X}(\x_{jl}\mid\bb\trans\x)f_{\bb\trans\X}(\bb\trans\x)\right\}d\x_{jl}\\
&&\times \int h^2u^2 K'(u)du\\
&&-\frac{1}{6h}\iint \frac{\partial^3}{\partial(\bb\trans\x)^3}\left\{
S_c(Z,\x_{jl},\bb\trans\x^*)S(Z,\bb\trans\x^*)f_{\X_l\mid\bb\trans\X}(\x_{jl}\mid\bb\trans\x^*)f_{\bb\trans\X}(\bb\trans\x^*)\right\}d\x_{jl}
\\
&&\times h^3u^3 K'(u)du\\
&=&
\frac{\partial}{\partial\bb\trans\x} \int S_c(Z,\x_{jl},\bb\trans\x)S(Z,\bb\trans\x)f_{\X_l\mid\bb\trans\X}(\x_{jl}\mid\bb\trans\x)f_{\bb\trans\X}(\bb\trans\x)d\x_{jl}\\
&&-\frac{h^2}{6}\int \int\frac{\partial^3}{\partial(\bb\trans\x)^3}
S_c(Z,\x_{jl},\bb\trans\x^*)S(Z,\bb\trans\x^*)f_{\X_l\mid\bb\trans\X}(\x_{jl}\mid\bb\trans\x^*)f_{\bb\trans\X}(\bb\trans\x^*)d\x_{jl}
\\
&&\times u^3 K'(u)du\\
&=&
\frac{\partial}{\partial \bb\trans\x}
S(Z,\bb\trans\x)f_{\bb\trans\X}(\bb\trans\x)E\{S_c(Z,\X_j)\mid\bb\trans\x\}\\
&&-\frac{h^2}{6}\int
\frac{\partial^3}{\partial(\bb\trans\x)^3}
S(Z,\bb\trans\x^*)f_{\bb\trans\X}(\bb\trans\x^*) E
\left\{S_c(Z,\X_j)\mid\bb\trans\x^*\right\}
u^3 K'(u)du\\
&=&
\frac{\partial}{\partial \bb\trans\x}
f_{\bb\trans\X}(\bb\trans\x)E\{I(Z_j\ge Z)\mid\bb\trans\x\}\\
&&-\frac{h^2}{6}\int
\frac{\partial^3}{\partial(\bb\trans\x)^3}
f_{\bb\trans\X}(\bb\trans\x^*) E
\left\{I(Z_j\ge Z)\mid\bb\trans\x^*\right\}
u^3 K'(u)du.
\ese
Hence the absolute bias is
\bse
&&\left|
E\left\{-\frac{1}{n}\sumj I(Z_j\ge Z)
K_h'(\bb\trans\X_j-\bb\trans\x)\right\}
-\frac{\partial}{\partial \bb\trans\x}
f_{\bb\trans\X}(\bb\trans\x)E\{I(Z_j\ge
Z)\mid\bb\trans\x\}\right|\\
&=&	\left|-\frac{h^2}{6}\int
\frac{\partial^3}{\partial(\bb\trans\x)^3}
f_{\bb\trans\X}(\bb\trans\x^*) E
\left\{I(Z_j\ge Z)\mid\bb\trans\x^*\right\}
u^3 K'(u)du\right|\\
&\le&\frac{h^2}{2}\sup_{\bb\trans\x}\left|
\frac{\partial^3}{\partial(\bb\trans\x)^3}
f_{\bb\trans\X}(\bb\trans\x) E
\left\{I(Z_j\ge Z)\mid\bb\trans\x\right\}\right|
\left\{\int u^2 K(u)du\right\}\\
\ese
and 
\bse
&&\var\left\{-\frac{1}{n}\sumj  I(Z_j\ge Z)  K_h'(\bb\trans\X_j-\bb\trans\x)\right\}\\
&=&\frac{1}{n}\left[E\left\{I(Z_j\ge Z)  K_h'(\bb\trans\X_j-\bb\trans\x)\left\}\right\{I(Z_j\ge Z)  K_h'(\bb\trans\X_j-\bb\trans\x)\right\}\trans\right.\\
&&\left.-\left\{EI(Z_j\ge Z)  K_h'(\bb\trans\X_j-\bb\trans\x)\right\}\left\{EI(Z_j\ge Z)  K_h'(\bb\trans\X_j-\bb\trans\x)\right\}\trans\right]\\
&=&\frac{1}{n}\left(\int\frac{1}{h^4}I(Z_j\ge Z)K'^2\left\{{(\bb\trans\x_j-\bb\trans\x)}/{h}\right\}\right.\\
&&\left.\times f_c(c,\x_j)f(t,\bb\trans\x_j)f_{\X_l\mid\bb\trans\X}(\x_{jl}\mid\bb\trans\x_j)f_{\bb\trans\X}(\bb\trans\x_j)dcdtd\x_{jl} d\bb\trans\x_j\right.\\
&&-\left.\left[\frac{\partial}{\partial \bb\trans\x}
f_{\bb\trans\X}(\bb\trans\x)E\{I(Z_j\ge
Z)\mid\bb\trans\x\}\right]\left[\frac{\partial}{\partial
	\bb\trans\x}f_{\bb\trans\X}(\bb\trans\x)E\{I(Z_j\ge Z)\mid\bb\trans\x\}\right]\trans+O(h^2)\right)\\
&=&\frac{1}{n}\int\frac{1}{h^3} S_c(Z,\x_{jl},\bb\trans\x+hu)\\
&&\times S(Z,\bb\trans\x+hu)f_{\X_l\mid\bb\trans\X}(\x_{jl}\mid\bb\trans\x+hu)f_{\bb\trans\X}(\bb\trans\x+hu)d\x_j K'^2(u)du+O(1/n)\\
&=&\frac{1}{nh^3}\int S_c(Z,\x_{jl},\bb\trans\x)S(Z,\bb\trans\x)f_{\bb\trans\X}(\bb\trans\x)f_{\X_l\mid\bb\trans\X}(\x_{jl}\mid\bb\trans\x)d\x_{jl}\int K'^2(u)du\\
&&+\frac{1}{2nh^3}
\frac{\partial^2}{\partial(\bb\trans\x)^2}
\iint S_c(Z,\x_{jl},\bb\trans\x^*)S(Z,\bb\trans\x^*)\\
&&\times
f_{\bb\trans\X}(\bb\trans\x^*)f_{\X_l\mid\bb\trans\X}(\x_{jl}\mid\bb\trans\x^*)d\x_{jl}
h^2u^2K'^2(u)du+O(1/n)\\
&\le&\frac{1}{nh^3} \sup_{\bb\trans\x}\left|f_{\bb\trans\X}(\bb\trans\x)E\{I(Z_j\ge Z)\mid\bb\trans\x\}\right|\int K'^2(u)du\\
&&+\frac{1}{2nh} \sup_{\bb\trans\x^*}\left|\frac{\partial^2}{\partial(\bb\trans\x)^2}
f_{\bb\trans\X}(\bb\trans\x^*)E\{I(Z_j\ge Z)\mid\bb\trans\x^*\}\right|\int u^2K'^2(u)du+O(1/n).
\ese
So
\be\label{eq:fbetai1}
-\frac{1}{n}\sumj I(Z_j\ge Z)  K_h'(\bb\trans\X_j-\bb\trans\x)
=\frac{\partial}{\partial \bb\trans\x} 
f_{\bb\trans\X}(\bb\trans\x)E\{I(Z_j\ge Z)\mid\bb\trans\x\}+O_p(n^{-1/2}h^{-3/2}+h^2)
\notag\\
\ee
uniformly under conditions \ref{assum:kernel}-\ref{assum:bandwidth}.

Following similar derivations, we have
\bse
&&E\left\{\frac{1}{n}\sumj \X_j I(Z_j\ge Z)  K_h(\bb\trans\X_j-\bb\trans\x)\right\}\\
&=&E\left\{\X_j I(Z_j\ge Z)  K_h(\bb\trans\X_j-\bb\trans\x)\right\}\\
&=&\int\frac{1}{h}\x_j I(c\ge Z)I(t\ge Z)  K\left\{(\bb\trans\x_j-\bb\trans\x)/h\right\}\\
&&\times f_c(c,\x_j)f(t,\bb\trans\x_j)f_{\X\mid\bb\trans\X}(\x_j\mid\bb\trans\x_j)f_{\bb\trans\X}(\bb\trans\x_j)dcdtd\x_j d\bb\trans\x_j\\
&=&\int\frac{1}{h}\x_j K\left\{(\bb\trans\x_j-\bb\trans\x)/h\right\}
S_c(Z,\x_{jl},\bb\trans\x_j)S(Z,\bb\trans\x_j)f_{\X_{l}\mid\bb\trans\X}(\x_{jl}\mid\bb\trans\x_j)\\
&&\times f_{\bb\trans\X}(\bb\trans\x_j)d\x_{jl} d\bb\trans\x_j\\
&=&\int (\bb\trans\x+hu,\x_{jl}\trans)\trans
K(u)S_c(Z,\x_{jl},\bb\trans\x+hu)S(Z,\bb\trans\x+hu)f_{\X_l\mid\bb\trans\X}(\x_{jl}\mid\bb\trans\x+hu)\\
&&\times f_{\bb\trans\X}(\bb\trans\x+hu)d\x_{jl} du\\
&=&\int(\bb\trans\x,\x_{jl}\trans)\trans
S_c(Z,\x_{jl},\bb\trans\x)S(Z,\bb\trans\x)f_{\X_l\mid\bb\trans\X}(\x_{jl}\mid\bb\trans\x)f_{\bb\trans\X}(\bb\trans\x)d\x_{jl}
\int K(u)du\\
&&+\int\frac{\partial}{\partial\bb\trans\x} \left\{(\bb\trans\x,\x_{jl}\trans)\trans S_c(Z,\x_{jl},\bb\trans\x)S(Z,\bb\trans\x)f_{\X_l\mid\bb\trans\X}(\x_{jl}\mid\bb\trans\x)f_{\bb\trans\X}(\bb\trans\x)\right\}d\x_{jl} \\
&&\times \int hu K(u)du\\
&&+\frac{1}{2}\iint\frac{\partial^2}{\partial(\bb\trans\x)^2}
\left\{(\bb\trans\x^*,\x_{jl}\trans)\trans
S_c(Z,\x_{jl},\bb\trans\x^*)S(Z,\bb\trans\x^*)f_{\X_l\mid\bb\trans\X}(\x_{jl}\mid\bb\trans\x^*)f_{\bb\trans\X}(\bb\trans\x^*)\right\}d\x_{jl}\\
&&\times  h^2u^2 K(u)du\\
&=&
\int(\bb\trans\x,\x_{jl}\trans)\trans
S_c(Z,\x_{jl},\bb\trans\x)S(Z,\bb\trans\x)f_{\X_l\mid\bb\trans\X}(\x_{jl}\mid\bb\trans\x)f_{\bb\trans\X}(\bb\trans\x)d\x_{jl}\\
&&+\frac{h^2}{2}\iint\frac{\partial^2}{\partial(\bb\trans\x)^2}
\left\{(\bb\trans\x^*,\x_{jl}\trans)\trans
S_c(Z,\x_{jl},\bb\trans\x^*)S(Z,\bb\trans\x^*)f_{\X_l\mid\bb\trans\X}(\x_{jl}\mid\bb\trans\x^*)f_{\bb\trans\X}(\bb\trans\x^*)\right\}d\x_{jl}\\
&&\times u^2 K(u)du\\
&=&
S(Z,\bb\trans\x)f_{\bb\trans\X}(\bb\trans\x)E\{\X_j
S_c(Z,\X_j)\mid\bb\trans\x\}\\
&&+\frac{h^2}{2}\int
\frac{\partial^2}{\partial(\bb\trans\x)^2}
S(Z,\bb\trans\x^*)f_{\bb\trans\X}(\bb\trans\x^*) E
\left\{\X_j
S_c(Z,\X_j)\mid\bb\trans\x^*\right\}
u^2 K(u)du\\
&=&
f_{\bb\trans\X}(\bb\trans\x)E\{\X_jI(Z_j\ge Z)\mid\bb\trans\x\}\\
&&+\frac{h^2}{2}\int
\frac{\partial^2}{\partial(\bb\trans\x)^2}
f_{\bb\trans\X}(\bb\trans\x^*) E
\left\{\X_jI(Z_j\ge Z)\mid\bb\trans\x^*\right\}
u^2 K(u)du.
\ese
Hence the absolute bias is
\bse
&&\left|
E\left\{\frac{1}{n}\sumj \X_j I(Z_j\ge Z)
K_h(\bb\trans\X_j-\bb\trans\x)\right\}
-f_{\bb\trans\X}(\bb\trans\x)E\{\X_j(Z_j\ge
Z)\mid\bb\trans\x\}\right|\\
&=&	\left|\frac{h^2}{2}\int
\frac{\partial^2}{\partial(\bb\trans\x)^2}
f_{\bb\trans\X}(\bb\trans\x^*) E
\left\{\X_jI(Z_j\ge Z)\mid\bb\trans\x^*\right\}
u^2 K(u)du\right|\\
&\le&\frac{h^2}{2}\sup_{\bb\trans\x}\left|
\frac{\partial^2}{\partial(\bb\trans\x)^2}
f_{\bb\trans\X}(\bb\trans\x) E
\left\{\X_jI(Z_j\ge Z)\mid\bb\trans\x\right\}\right|
\left\{\int u^2 K(u)du\right\}\\
\ese
and 
\bse
&&\var\left\{-\frac{1}{n}\sumj \X_j I(Z_j\ge Z)  K_h(\bb\trans\X_j-\bb\trans\x)\right\}\\
&=&\frac{1}{n}\left[E\left\{\X_jI(Z_j\ge Z)  K_h(\bb\trans\X_j-\bb\trans\x)\left\}\right\{\X_jI(Z_j\ge Z)  K_h(\bb\trans\X_j-\bb\trans\x)\right\}\trans\right.\\
&&\left.-\left\{E\X_j I(Z_j\ge Z)  K_h(\bb\trans\X_j-\bb\trans\x)\right\}\left\{E\X_j I(Z_j\ge Z)  K_h(\bb\trans\X_j-\bb\trans\x)\right\}\trans\right]\\
&=&\frac{1}{n}\left(\int\frac{1}{h^2}\x_j\x_j\trans I(Z_j\ge Z)K^2\left\{{(\bb\trans\x_j-\bb\trans\x)}/{h}\right\}\right.\\
&&\left.\times f_c(c,\x_j)f(t,\bb\trans\x_j)f_{\X_l\mid\bb\trans\X}(\x_{jl}\mid\bb\trans\x_j)f_{\bb\trans\X}(\bb\trans\x_j)dcdtd\x_{jl} d\bb\trans\x_j\right.\\
&&-\left.\left[
f_{\bb\trans\X}(\bb\trans\x)E\{\X_j
I(Z_j\ge
Z)\mid\bb\trans\x\}\right]\left[f_{\bb\trans\X}(\bb\trans\x)E\{\X_j
I(Z_j\ge Z)\mid\bb\trans\x\}\right]\trans\right.\\
&&\left.+O(h^2)\right)\\
&=&\frac{1}{n}\int\frac{1}{h}(\bb\trans\x+hu,\x_{jl}\trans)\trans(\bb\trans\x+hu,\x_{jl}\trans) S_c(Z,\x_{jl},\bb\trans\x+hu)\\
&&\times S(Z,\bb\trans\x+hu)f_{\X_l\mid\bb\trans\X}(\x_{jl}\mid\bb\trans\x+hu)f_{\bb\trans\X}(\bb\trans\x+hu)d\x_j K^2(u)du\\
&&+O(1/n)\\
&=&\frac{1}{nh}\int(\bb\trans\x,\x_{jl}\trans)\trans(\bb\trans\x,\x_{jl}\trans)S_c(Z,\x_{jl},\bb\trans\x)S(Z,\bb\trans\x)f_{\bb\trans\X}(\bb\trans\x)f_{\X_l\mid\bb\trans\X}(\x_{jl}\mid\bb\trans\x)d\x_{jl}
\\
&&\times\int K^2(u)du\\
&&+\frac{1}{2nh}
\frac{\partial^2}{\partial(\bb\trans\x)^2}
\iint(\bb\trans\x^*,\x_{jl}\trans)\trans(\bb\trans\x^*,\x_{jl}\trans)S_c(Z,\x_{jl},\bb\trans\x^*)S(Z,\bb\trans\x^*)\\
&&\times
f_{\bb\trans\X}(\bb\trans\x^*)f_{\X_l\mid\bb\trans\X}(\x_{jl}\mid\bb\trans\x^*)d\x_{jl}
h^2u^2K^2(u)du+O(1/n)\\
&\le&\frac{1}{nh} \sup_{\bb\trans\x}\left|f_{\bb\trans\X}(\bb\trans\x)E\{\X_j\X_j\trans
I(Z_j\ge Z)\mid\bb\trans\x\}\right|\int K^2(u)du\\
&&+\frac{h}{2n} \sup_{\bb\trans\x^*}\left|\frac{\partial^2}{\partial(\bb\trans\x)^2}
f_{\bb\trans\X}(\bb\trans\x^*)E\{\X_j\X_j\trans
I(Z_j\ge Z)\mid\bb\trans\x^*\}\right|\int u^2K^2(u)du+O(1/n).
\ese
So
\be\label{eq:fbetaxi}
\frac{1}{n}\sumj \X_j I(Z_j\ge Z)  K_h(\bb\trans\X_j-\bb\trans\x)
=	f_{\bb\trans\X}(\bb\trans\x)E\{\X_j I(Z_j\ge Z)\mid\bb\trans\x\}+O_p(n^{-1/2}h^{-1/2}+h^2)
\notag\\
\ee
uniformly under conditions \ref{assum:kernel}-\ref{assum:bandwidth}.

\bse
&&E\left\{\frac{1}{n}\sumj I(Z_j\ge Z)  K_h(\bb\trans\X_j-\bb\trans\x)\right\}\\
&=&E\left\{I(Z_j\ge Z)  K_h(\bb\trans\X_j-\bb\trans\x)\right\}\\
&=&\int\frac{1}{h}I(c\ge Z)I(t\ge Z)  K\left\{(\bb\trans\x_j-\bb\trans\x)/h\right\}\\
&&\times f_c(c,\x_j)f(t,\bb\trans\x_j)f_{\X\mid\bb\trans\X}(\x_j\mid\bb\trans\x_j)f_{\bb\trans\X}(\bb\trans\x_j)dcdtd\x_j d\bb\trans\x_j\\
&=&\int\frac{1}{h}K\left\{(\bb\trans\x_j-\bb\trans\x)/h\right\}
S_c(Z,\x_{jl},\bb\trans\x_j)S(Z,\bb\trans\x_j)f_{\X_{l}\mid\bb\trans\X}(\x_{jl}\mid\bb\trans\x_j)\\
&&\times f_{\bb\trans\X}(\bb\trans\x_j)d\x_{jl} d\bb\trans\x_j\\
&=&\int
K(u)S_c(Z,\x_{jl},\bb\trans\x+hu)S(Z,\bb\trans\x+hu)f_{\X_l\mid\bb\trans\X}(\x_{jl}\mid\bb\trans\x+hu)\\
&&\times f_{\bb\trans\X}(\bb\trans\x+hu)d\x_{jl} du\\
&=&\int
S_c(Z,\x_{jl},\bb\trans\x)S(Z,\bb\trans\x)f_{\X_l\mid\bb\trans\X}(\x_{jl}\mid\bb\trans\x)f_{\bb\trans\X}(\bb\trans\x)d\x_{jl}
\int K(u)du\\
&&+\int\frac{\partial}{\partial\bb\trans\x} \left\{ S_c(Z,\x_{jl},\bb\trans\x)S(Z,\bb\trans\x)f_{\X_l\mid\bb\trans\X}(\x_{jl}\mid\bb\trans\x)f_{\bb\trans\X}(\bb\trans\x)\right\}d\x_{jl} \\
&&\times \int hu K(u)du\\
&&+\frac{1}{2}\iint\frac{\partial^2}{\partial(\bb\trans\x)^2}
\left\{
S_c(Z,\x_{jl},\bb\trans\x^*)S(Z,\bb\trans\x^*)f_{\X_l\mid\bb\trans\X}(\x_{jl}\mid\bb\trans\x^*)f_{\bb\trans\X}(\bb\trans\x^*)\right\}d\x_{jl}\\
&&\times  h^2u^2 K(u)du\\
&=&
\int S_c(Z,\x_{jl},\bb\trans\x)S(Z,\bb\trans\x)f_{\X_l\mid\bb\trans\X}(\x_{jl}\mid\bb\trans\x)f_{\bb\trans\X}(\bb\trans\x)d\x_{jl}\\
&&+\frac{h^2}{2}\iint\frac{\partial^2}{\partial(\bb\trans\x)^2}
\left\{
S_c(Z,\x_{jl},\bb\trans\x^*)S(Z,\bb\trans\x^*)f_{\X_l\mid\bb\trans\X}(\x_{jl}\mid\bb\trans\x^*)f_{\bb\trans\X}(\bb\trans\x^*)\right\}d\x_{jl}\\
&&\times u^2 K(u)du\\
&=&
S(Z,\bb\trans\x)f_{\bb\trans\X}(\bb\trans\x)E\{
S_c(Z,\X_j)\mid\bb\trans\x\}\\
&&+\frac{h^2}{2}\int
\frac{\partial^2}{\partial(\bb\trans\x)^2}
S(Z,\bb\trans\x^*)f_{\bb\trans\X}(\bb\trans\x^*) E
\left\{S_c(Z,\X_j)\mid\bb\trans\x^*\right\}
u^2 K(u)du\\
&=&
f_{\bb\trans\X}(\bb\trans\x)E\{I(Z_j\ge Z)\mid\bb\trans\x\}\\
&&+\frac{h^2}{2}\int
\frac{\partial^2}{\partial(\bb\trans\x)^2}
f_{\bb\trans\X}(\bb\trans\x^*) E
\left\{I(Z_j\ge Z)\mid\bb\trans\x^*\right\}
u^2 K(u)du.
\ese
Hence the absolute bias is
\bse
&&\left|
E\left\{\frac{1}{n}\sumj I(Z_j\ge Z)
K_h(\bb\trans\X_j-\bb\trans\x)\right\}
-f_{\bb\trans\X}(\bb\trans\x)E\{I(Z_j\ge
Z)\mid\bb\trans\x\}\right|\\
&=&	\left|\frac{h^2}{2}\int
\frac{\partial^2}{\partial(\bb\trans\x)^2}
f_{\bb\trans\X}(\bb\trans\x^*) E
\left\{I(Z_j\ge Z)\mid\bb\trans\x^*\right\}
u^2 K(u)du\right|\\
&\le&\frac{h^2}{2}\sup_{\bb\trans\x}\left|
\frac{\partial^2}{\partial(\bb\trans\x)^2}
f_{\bb\trans\X}(\bb\trans\x) E
\left\{I(Z_j\ge Z)\mid\bb\trans\x\right\}\right|
\left\{\int u^2 K(u)du\right\}\\
\ese
and 
\bse
&&\var\left\{\frac{1}{n}\sumj I(Z_j\ge Z)  K_h(\bb\trans\X_j-\bb\trans\x)\right\}\\
&=&\frac{1}{n}\left[E\left\{I(Z_j\ge Z)  K_h(\bb\trans\X_j-\bb\trans\x)\left\}\right\{\X_jI(Z_j\ge Z)  K_h(\bb\trans\X_j-\bb\trans\x)\right\}\trans\right.\\
&&\left.-\left\{E I(Z_j\ge Z)  K_h(\bb\trans\X_j-\bb\trans\x)\right\}\left\{E I(Z_j\ge Z)  K_h(\bb\trans\X_j-\bb\trans\x)\right\}\trans\right]\\
&=&\frac{1}{n}\left(\int\frac{1}{h^2}I(Z_j\ge Z)K^2\left\{{(\bb\trans\x_j-\bb\trans\x)}/{h}\right\}\right.\\
&&\left.\times f_c(c,\x_j)f(t,\bb\trans\x_j)f_{\X_l\mid\bb\trans\X}(\x_{jl}\mid\bb\trans\x_j)f_{\bb\trans\X}(\bb\trans\x_j)dcdtd\x_{jl} d\bb\trans\x_j\right.\\
&&-\left.\left[
f_{\bb\trans\X}(\bb\trans\x)E\{I(Z_j\ge
Z)\mid\bb\trans\x\}\right]\left[f_{\bb\trans\X}(\bb\trans\x)E\{I(Z_j\ge Z)\mid\bb\trans\x\}\right]\trans+O(h^2)\right)\\
&=&\frac{1}{n}\int\frac{1}{h}S_c(Z,\x_{jl},\bb\trans\x+hu)\\
&&\times S(Z,\bb\trans\x+hu)f_{\X_l\mid\bb\trans\X}(\x_{jl}\mid\bb\trans\x+hu)f_{\bb\trans\X}(\bb\trans\x+hu)d\x_j K^2(u)du+O(1/n)\\
&=&\frac{1}{nh}\int S_c(Z,\x_{jl},\bb\trans\x)S(Z,\bb\trans\x)f_{\bb\trans\X}(\bb\trans\x)f_{\X_l\mid\bb\trans\X}(\x_{jl}\mid\bb\trans\x)d\x_{jl}\int K^2(u)du\\
&&+\frac{1}{2nh}
\frac{\partial^2}{\partial(\bb\trans\x)^2}
\iint S_c(Z,\x_{jl},\bb\trans\x^*)S(Z,\bb\trans\x^*)\\
&&\times
f_{\bb\trans\X}(\bb\trans\x^*)f_{\X_l\mid\bb\trans\X}(\x_{jl}\mid\bb\trans\x^*)d\x_{jl}
h^2u^2K^2(u)du+O(1/n)\\
&\le&\frac{1}{nh} \sup_{\bb\trans\x}\left|f_{\bb\trans\X}(\bb\trans\x)E\{I(Z_j\ge Z)\mid\bb\trans\x\}\right|\int K^2(u)du\\
&&+\frac{h}{2n} \sup_{\bb\trans\x^*}\left|\frac{\partial^2}{\partial(\bb\trans\x)^2}
f_{\bb\trans\X}(\bb\trans\x^*)E\{I(Z_j\ge Z)\mid\bb\trans\x^*\}\right|\int u^2K^2(u)du+O(1/n).
\ese

So
\be\label{eq:fbetai}
\frac{1}{n}\sumj I(Z_j\ge Z)  K_h(\bb\trans\X_j-\bb\trans\x)
=	f_{\bb\trans\X}(\bb\trans\x)E\{I(Z_j\ge Z)\mid\bb\trans\x\}+O_p(n^{-1/2}h^{-1/2}+h^2)
\notag\\\ee
uniformly under conditions \ref{assum:kernel}-\ref{assum:bandwidth}.
\bse
&&E\left\{\frac{1}{n}\sumj I(Z_j\ge Z)  K_h''(\bb\trans\X_j-\bb\trans\x)\right\}\\
&=&E\left\{I(Z_j\ge Z)  K_h''(\bb\trans\X_j-\bb\trans\x)\right\}\\
&=&\int\frac{1}{h^3}I(c\ge Z)I(t\ge Z)  K''\left\{(\bb\trans\x_j-\bb\trans\x)/h\right\}\\
&&\times f_c(c,\x_j)f(t,\bb\trans\x_j)f_{\X\mid\bb\trans\X}(\x_j\mid\bb\trans\x_j)f_{\bb\trans\X}(\bb\trans\x_j)dcdtd\x_j d\bb\trans\x_j\\
&=&\int\frac{1}{h^3}K''\left\{(\bb\trans\x_j-\bb\trans\x)/h\right\}
S_c(Z,\x_{jl},\bb\trans\x_j)S(Z,\bb\trans\x_j)f_{\X_{l}\mid\bb\trans\X}(\x_{jl}\mid\bb\trans\x_j)\\
&&\times f_{\bb\trans\X}(\bb\trans\x_j)d\x_{jl} d\bb\trans\x_j\\
&=&\int\frac{1}{h^2}
K''(u)S_c(Z,\x_{jl},\bb\trans\x+hu)S(Z,\bb\trans\x+hu)f_{\X_l\mid\bb\trans\X}(\x_{jl}\mid\bb\trans\x+hu)\\
&&\times f_{\bb\trans\X}(\bb\trans\x+hu)d\x_{jl} du\\
&=&\frac{1}{2h^2}\int\frac{\partial^2}{\partial(\bb\trans\x)^2}
\left\{
S_c(Z,\x_{jl},\bb\trans\x)S(Z,\bb\trans\x)f_{\X_l\mid\bb\trans\X}(\x_{jl}\mid\bb\trans\x)f_{\bb\trans\X}(\bb\trans\x)\right\}d\x_{jl}\\
&&\times  h^2\int u^2 K''(u)du\\
&&+\frac{h^2}{24}\iint\frac{\partial^4}{\partial(\bb\trans\x)^4}
\left\{
S_c(Z,\x_{jl},\bb\trans\x^*)S(Z,\bb\trans\x^*)f_{\X_l\mid\bb\trans\X}(\x_{jl}\mid\bb\trans\x^*)f_{\bb\trans\X}(\bb\trans\x^*)\right\}d\x_{jl}\\
&&\times  u^4 K''(u)du\\
&=&
\int\frac{\partial^2}{\partial(\bb\trans\x)^2}
\left\{
S_c(Z,\x_{jl},\bb\trans\x)S(Z,\bb\trans\x)f_{\X_l\mid\bb\trans\X}(\x_{jl}\mid\bb\trans\x)f_{\bb\trans\X}(\bb\trans\x)\right\}d\x_{jl}\\
&&+\frac{h^2}{24}\iint\frac{\partial^4}{\partial(\bb\trans\x)^4}
\left\{
S_c(Z,\x_{jl},\bb\trans\x^*)S(Z,\bb\trans\x^*)f_{\X_l\mid\bb\trans\X}(\x_{jl}\mid\bb\trans\x^*)f_{\bb\trans\X}(\bb\trans\x^*)\right\}d\x_{jl}\\
&&\times  u^4 K''(u)du\\
&=&
\frac{\partial^2}{\partial(\bb\trans\x)^2}S(Z,\bb\trans\x)f_{\bb\trans\X}(\bb\trans\x)E\{
S_c(Z,\X_j)\mid\bb\trans\x\}\\
&&+\frac{h^2}{24}\iint\frac{\partial^4}{\partial(\bb\trans\x)^4}
\left\{
S_c(Z,\x_{jl},\bb\trans\x^*)S(Z,\bb\trans\x^*)f_{\X_l\mid\bb\trans\X}(\x_{jl}\mid\bb\trans\x^*)f_{\bb\trans\X}(\bb\trans\x^*)\right\}d\x_{jl}\\
&&\times  u^4 K''(u)du\\
&=&
\frac{\partial^2}{\partial(\bb\trans\x)^2}f_{\bb\trans\X}(\bb\trans\x)E\{I(Z_j\ge Z)\mid\bb\trans\x\}\\
&&+\frac{h^2}{24}\iint\frac{\partial^4}{\partial(\bb\trans\x)^4}
\left\{
S_c(Z,\x_{jl},\bb\trans\x^*)S(Z,\bb\trans\x^*)f_{\X_l\mid\bb\trans\X}(\x_{jl}\mid\bb\trans\x^*)f_{\bb\trans\X}(\bb\trans\x^*)\right\}d\x_{jl}\\
&&\times  u^4 K''(u)du.
\ese
Hence the absolute bias is
\bse
&&\left|
E\left\{\frac{1}{n}\sumj I(Z_j\ge Z)
K_h''(\bb\trans\X_j-\bb\trans\x)\right\}
-\frac{\partial^2}{\partial(\bb\trans\x)^2}f_{\bb\trans\X}(\bb\trans\x)E\{I(Z_j\ge Z)\mid\bb\trans\x\}\right|\\
&=&\left|\frac{h^2}{24}\iint\frac{\partial^4}{\partial(\bb\trans\x)^4}
\left\{
S_c(Z,\x_{jl},\bb\trans\x^*)S(Z,\bb\trans\x^*)f_{\X_l\mid\bb\trans\X}(\x_{jl}\mid\bb\trans\x^*)f_{\bb\trans\X}(\bb\trans\x^*)\right\}d\x_{jl}u^4 K''(u)du\right|\\
&\le&\frac{h^2}{2}\sup_{\bb\trans\x^*}\left|
\frac{\partial^4}{\partial(\bb\trans\x)^4}
f_{\bb\trans\X}(\bb\trans\x) E
\left\{I(Z_j\ge Z)\mid\bb\trans\x\right\}\right|
\left\{\int u^2 K(u)du\right\}.
\ese
To analyze the variance, 
\bse
&&\var\left\{\frac{1}{n}\sumj I(Z_j\ge Z)  K_h''(\bb\trans\X_j-\bb\trans\x)\right\}\\
&=&\frac{1}{n}\left[E\left\{I(Z_j\ge Z)  K_h''(\bb\trans\X_j-\bb\trans\x)\left\}\right\{\X_jI(Z_j\ge Z)  K_h''(\bb\trans\X_j-\bb\trans\x)\right\}\trans\right.\\
&&\left.-\left\{E I(Z_j\ge Z)  K_h''(\bb\trans\X_j-\bb\trans\x)\right\}\left\{E I(Z_j\ge Z)  K_h''(\bb\trans\X_j-\bb\trans\x)\right\}\trans\right]\\
&=&\frac{1}{n}\left(\int\frac{1}{h^6}I(Z_j\ge Z)K''^2\left\{{(\bb\trans\x_j-\bb\trans\x)}/{h}\right\}\right.\\
&&\left.\times f_c(c,\x_j)f(t,\bb\trans\x_j)f_{\X_l\mid\bb\trans\X}(\x_{jl}\mid\bb\trans\x_j)f_{\bb\trans\X}(\bb\trans\x_j)dcdtd\x_{jl} d\bb\trans\x_j\right.\\
&&-\left.\left[\frac{\partial^2}{\partial (\bb\trans\x)^2}
f_{\bb\trans\X}(\bb\trans\x)E\{I(Z_j\ge
Z)\mid\bb\trans\x\}\right]\left[\frac{\partial^2}{(\partial
	\bb\trans\x)^2}f_{\bb\trans\X}(\bb\trans\x)E\{I(Z_j\ge Z)\mid\bb\trans\x\}\right]\trans+O(h^2)\right)\\
&=&\frac{1}{n}\int\frac{1}{h^5}S_c(Z,\x_{jl},\bb\trans\x+hu)\\
&&\times S(Z,\bb\trans\x+hu)f_{\X_l\mid\bb\trans\X}(\x_{jl}\mid\bb\trans\x+hu)f_{\bb\trans\X}(\bb\trans\x+hu)d\x_j K''^2(u)du+O(1/n)\\
&=&\frac{1}{nh^5}\int S_c(Z,\x_{jl},\bb\trans\x)S(Z,\bb\trans\x)f_{\bb\trans\X}(\bb\trans\x)f_{\X_l\mid\bb\trans\X}(\x_{jl}\mid\bb\trans\x)d\x_{jl}\int K''^2(u)du\\
&&+\frac{1}{2nh^5}
\frac{\partial^2}{\partial(\bb\trans\x)^2}
\iint S_c(Z,\x_{jl},\bb\trans\x^*)S(Z,\bb\trans\x^*)\\
&&\times
f_{\bb\trans\X}(\bb\trans\x^*)f_{\X_l\mid\bb\trans\X}(\x_{jl}\mid\bb\trans\x^*)d\x_{jl}
h^2u^2K^2(u)du+O(1/n)\\
&\le&\frac{1}{nh^5} \sup_{\bb\trans\x}\left|f_{\bb\trans\X}(\bb\trans\x)E\{I(Z_j\ge Z)\mid\bb\trans\x\}\right|\int K''^2(u)du\\
&&+\frac{1}{2nh^3} \sup_{\bb\trans\x^*}\left|\frac{\partial^2}{\partial(\bb\trans\x)^2}
f_{\bb\trans\X}(\bb\trans\x^*)E\{I(Z_j\ge Z)\mid\bb\trans\x^*\}\right|\int u^2K''^2(u)du+O(1/n).
\ese

So
\be\label{eq:fbetai2}
\frac{1}{n}\sumj I(Z_j\ge Z)  K_h(\bb\trans\X_j-\bb\trans\x)
=\frac{\partial^2}{\partial(\bb\trans\x)^2}	f_{\bb\trans\X}(\bb\trans\x)E\{I(Z_j\ge Z)\mid\bb\trans\x\}+O_p(n^{-1/2}h^{-5/2}+h^2)
\notag\\\ee
uniformly under conditions \ref{assum:kernel}-\ref{assum:bandwidth}.

Combining the results of (\ref{eq:fbeta}) and (\ref{eq:fbetai}), we have
\be
\wh E\left\{Y(Z)\mid\bb\trans\X\right\}
&=&\frac{\sumj I(Z_j\ge Z) K_h(\bb\trans\X_j-\bb\trans\X)}
{\sumj K_h(\bb\trans\X_j-\bb\trans\X)}\n\\
&=&\frac{f_{\bb\trans\X}(\bb\trans\X)E\{I(Z_j\ge Z)\mid\bb\trans\X\}+O_p(n^{-1/2}h^{-1/2}+h^2)}
{f_{\bb\trans\X}(\bb\trans\X)+O_p(n^{-1/2}h^{-1/2}+h^2)}\n\\
&=&\frac{f_{\bb\trans\X}(\bb\trans\X)E\{Y(Z)\mid\bb\trans\X\}}
{f_{\bb\trans\X}(\bb\trans\X)}+O_p(n^{-1/2}h^{-1/2}+h^2)\n\\
&=&E\{Y(Z)\mid\bb\trans\X\}+O_p\{(nh)^{-1/2}+h^2\},\n
\ee
uniformly under conditions \ref{assum:kernel}-\ref{assum:bandwidth}. Combining the results of (\ref{eq:fbeta}) and (\ref{eq:fbetaxi}), we have
\be
\wh E\left\{\X Y(Z)\mid\bb\trans\X\right\}
&=&\frac{\sumj \X_j I(Z_j\ge Z)  K_h(\bb\trans\X_j-\bb\trans\X)}
{\sumj K_h(\bb\trans\X_j-\bb\trans\X)}\n\\
&=&\frac{f_{\bb\trans\X}(\bb\trans\X)E\{\X_j I(Z_j\ge Z)\mid\bb\trans\X\}+O_p(n^{-1/2}h^{-1/2}+h^2)}
{f_{\bb\trans\X}(\bb\trans\X)+O_p(n^{-1/2}h^{-1/2}+h^2)}\n\\
&=&\frac{f_{\bb\trans\X}(\bb\trans\X)E\{\X Y(Z)\mid\bb\trans\X\}}
{f_{\bb\trans\X}(\bb\trans\X)}+O_p(n^{-1/2}h^{-1/2}+h^2)\n\\
&=&E\{\X Y(Z)\mid\bb\trans\X\}+O_p\{(nh)^{-1/2}+h^2\} \n
\ee
uniformly under conditions \ref{assum:kernel}-\ref{assum:bandwidth}. Combining the results of 
(\ref{eq:fbeta}), (\ref{eq:fbeta1}), (\ref{eq:fbetai}) and 
(\ref{eq:fbetai1}), we have
\be
&&\frac{\partial}{\partial\bb\trans\X}\wh E\left\{Y(Z)\mid\bb\trans\X\right\}\n\\
&=&-\frac{\sumj I(Z_j\ge Z)  K_h'(\bb\trans\X_j-\bb\trans\X)}
{\sumj K_h(\bb\trans\X_j-\bb\trans\X)}\n\\
&&+\frac{\{\sumj I(Z_j\ge Z)K_h(\bb\trans\X_j-\bb\trans\X)\}
	\{\sumj K_h'(\bb\trans\X_j-\bb\trans\X)\}}
{\{\sumj K_h(\bb\trans\X_j-\bb\trans\X)\}^2}\n\\
&=&\frac{\partial 
	f_{\bb\trans\X}(\bb\trans\X)E\{I(Z_j\ge Z)\mid\bb\trans\X\}/\partial \bb\trans\X+O_p(n^{-1/2}h^{-3/2}+h^2)}
{f_{\bb\trans\X}(\bb\trans\X)+O_p(n^{-1/2}h^{-1/2}+h^2)}\n\\
&&+\frac{\left[f_{\bb\trans\X}(\bb\trans\X)E\{I(Z_j\ge Z)\mid\bb\trans\X\}+O_p(n^{-1/2}h^{-1/2}+h^2)\right]
	\left[-f_{\bb\trans\X}'(\bb\trans\X)+O_p(n^{-1/2}h^{-3/2}+h^2)\right]}
{f_{\bb\trans\X}^2(\bb\trans\X)+O(n^{-1/2}h^{-1/2}+h^2)}\n\\
&=&\frac{f_{\bb\trans\X}'(\bb\trans\X)E\{I(Z_j\ge Z)\mid\bb\trans\X\}}{f_{\bb\trans\X}(\bb\trans\X)}+\frac{\partial E\{I(Z_j\ge Z)\mid\bb\trans\X\}}{\partial\bb\trans\X}\n\\
&&-\frac{f_{\bb\trans\X}'(\bb\trans\X)E\{I(Z_j\ge Z)\mid\bb\trans\X\}}
{f_{\bb\trans\X}^2(\bb\trans\X)}+O_p(n^{-1/2}h^{-3/2}+h^2)\n\\
&=&\frac{\partial}{\partial\bb\trans\X} E\{Y(Z)\mid\bb\trans\X\}+O_p\{(nh^3)^{-1/2}+h^2\}\n
\ee
uniformly under conditions \ref{assum:kernel}-\ref{assum:bandwidth}. 
Finally, combining the results of 
(\ref{eq:fbeta}), (\ref{eq:fbeta1}), (\ref{eq:fbetaxi}) and 
(\ref{eq:fbetaxi1}), we have
\be
&&\frac{\partial}{\partial\bb\trans\X}\wh E\left\{\X Y(Z)\mid\bb\trans\X\right\}\n\\
&=&-\frac{\sumj \X_j I(Z_j\ge Z)  K_h'(\bb\trans\X_j-\bb\trans\X)}
{\sumj K_h(\bb\trans\X_j-\bb\trans\X)}\n\\
&&+\frac{\{\sumj \X_j I(Z_j\ge Z)K_h(\bb\trans\X_j-\bb\trans\X)\}
	\{\sumj K_h'(\bb\trans\X_j-\bb\trans\X)\}}
{\{\sumj K_h(\bb\trans\X_j-\bb\trans\X)\}^2}\n\\
&=&\frac{\partial 
	f_{\bb\trans\X}(\bb\trans\X)E\{\X_j I(Z_j\ge Z)\mid\bb\trans\X\}/\partial \bb\trans\X+O_p(n^{-1/2}h^{-3/2}+h^2)}
{f_{\bb\trans\X}(\bb\trans\X)+O_p(n^{-1/2}h^{-1/2}+h^2)}\n\\
&&+\frac{\left[f_{\bb\trans\X}(\bb\trans\X)E\{\X_j I(Z_j\ge Z)\mid\bb\trans\X\}+O_p(n^{-1/2}h^{-1/2}+h^2)\right]
	\left[-f_{\bb\trans\X}'(\bb\trans\X)+O_p(n^{-1/2}h^{-3/2}+h^2)\right]}
{f_{\bb\trans\X}^2(\bb\trans\X)+O(n^{-1/2}h^{-1/2}+h^2)}\n\\
&=&\frac{f_{\bb\trans\X}'(\bb\trans\X)E\{\X Y(Z)\mid\bb\trans\X\}}
{f_{\bb\trans\X}(\bb\trans\X)}+\frac{f_{\bb\trans\X}(\bb\trans\X)\partial E\{\X Y(Z)\mid\bb\trans\X\}/\partial\bb\trans\X}
{f_{\bb\trans\X}(\bb\trans\X)}\n\\
&&-\frac{f_{\bb\trans\X}'(\bb\trans\X)E\{\X Y(Z)\mid\bb\trans\X\}}
{f_{\bb\trans\X}(\bb\trans\X)}+O_p(n^{-1/2}h^{-3/2}+h^2)\n\\
&=&\frac{\partial}{\partial\bb\trans\X}E\{\X Y(Z)\mid\bb\trans\X\}+O_p\{(nh^3)^{-1/2}+h^2\} \n
\ee
uniformly under conditions \ref{assum:kernel}-\ref{assum:bandwidth}.

Now we inspect the consistency of the Kaplan Meier estimator on the
hazard function and its derivatives. 
Let 
$A=n^{-1}\sumj I(Z_j\ge
Z_i)K_h(\bb\trans\X_j-\bb\trans\x)-f_{\bb\trans\X}(\bb\trans\x)E\{I(Z\ge 
Z_i)\mid\bb\trans\x\}$.
\bse
\widehat{\lambda}(z,\bb\trans\x) &=& \int_0^{\infty} K_b(t-z) d\widehat{\Lambda}(t|\bb\trans\x)\\
&=& \sum_{i=1}^nK_b(Z_i-z)
\frac{\Delta_iK_h(\bb\trans\X_i-\bb\trans\x)}{\sumj I(Z_j\ge
	Z_i)K_h(\bb\trans\X_j-\bb\trans\x)}\\
&=&\frac{1}{n}\sum_{i=1}^nK_b(Z_i-z)
\frac{\Delta_iK_h(\bb\trans\X_i-\bb\trans\x)}{f_{\bb\trans\X}(\bb\trans\x)E\{I(Z\ge
	Z_i)\mid\bb\trans\x\}+A}\\
&=&\frac{1}{n}\sum_{i=1}^n 
\frac{K_b(Z_i-z)\Delta_iK_h(\bb\trans\X_i-\bb\trans\x)}{f_{\bb\trans\X}(\bb\trans\x)E\{I(Z\ge 
	Z_i)\mid\bb\trans\x\}}\{1+O_p(A)\},
\ese
We first inspect
\bse
\frac{1}{n}\sum_{i=1}^n 
\frac{K_b(Z_i-z)\Delta_iK_h(\bb\trans\X_i-\bb\trans\x)}{f_{\bb\trans\X}(\bb\trans\x)E\{I(Z\ge 
	Z_i)\mid\bb\trans\x\}}.
\ese
First,
\bse
&& E \left[\frac{1}{n}\sum_{i=1}^n 
\frac{K_b(Z_i-z)\Delta_iK_h(\bb\trans\X_i-\bb\trans\x)}{f_{\bb\trans\X}(\bb\trans\x)E\{I(Z\ge 
	Z_i)\mid\bb\trans\x\}}\right]\\
&=&E\left[
\frac{K_b(Z_i-z)\Delta_iK_h(\bb\trans\X_i-\bb\trans\x)}{f_{\bb\trans\X}(\bb\trans\x)
	S(Z_i,\bb\trans\x)E\{S_c( 
	Z_i,\X_i)\mid\bb\trans\x\}}\right]\\
&=&\iint
\frac{K_b(z_i-z)K_h(\bb\trans\x_i-\bb\trans\x)}{f_{\bb\trans\X}(\bb\trans\x)
	S(z_i,\bb\trans\x)E\{S_c( 
	z_i,\X_i)\mid\bb\trans\x\}}
f(z_i,\bb\trans\x_i)S_c(z_i,\x_i)f_{\X\mid\bb\trans\X}(\x_i,\bb\trans\x_i)\\
&&\times f_{\bb\trans\X}(\bb\trans\x_i)dz_id\x_id\bb\trans\x_i\\
&=&\iint
\frac{K_b(z_i-z)K_h(\bb\trans\x_i-\bb\trans\x)E\{S_c( 
	z_i,\X_i)\mid\bb\trans\x_i\}}{f_{\bb\trans\X}(\bb\trans\x)
	S(z_i,\bb\trans\x)E\{S_c( 
	z_i,\X_i)\mid\bb\trans\x\}}
f(z_i,\bb\trans\x_i)\\
&&\times f_{\bb\trans\X}(\bb\trans\x_i)dz_id\bb\trans\x_i\\
&=&\iint
\frac{K(v)K(u)E\{S_c( 
	z+bv,\X_i)\mid\bb\trans\x+hu\}}{f_{\bb\trans\X}(\bb\trans\x)
	S(z+bv,\bb\trans\x)E\{S_c( 
	z+bv,\X_i)\mid\bb\trans\x\}}
f(z+bv,\bb\trans\x+hu)\\
&&\times f_{\bb\trans\X}(\bb\trans\x+hu)dvdu\\
&=&\iint  K(v)K(u)
\lambda(z,\bb\trans\x)dvdu +\frac{b^2\partial^2}{2\partial z^2}\iint  
\frac{K(v)K(u)}{
	S(z^*,\bb\trans\x)}
f(z^*,\bb\trans\x)v^2dvdu\\
&&+
\frac{h^2\partial^2}{2\partial (\bb\trans\x)^2}\iint \frac{ K(v)K(u) E\{S_c( 
	z,\X_i)\mid\bb\trans\x^*\}
	f(z,\bb\trans\x^*) f_{\bb\trans\X}(\bb\trans\x^*)}{f_{\bb\trans\X}(\bb\trans\x)
	S(z,\bb\trans\x)
	E\{S_c( 
	z,\X_i)\mid\bb\trans\x\} }u^2dvdu\\
&=&
\lambda(z,\bb\trans\x) +\frac{b^2\partial^2}{2\partial z^2}\iint  
\frac{K(v)K(u)}{
	S(z^*,\bb\trans\x)}
f(z^*,\bb\trans\x)v^2dvdu\\
&&+
\frac{h^2\partial^2}{2\partial (\bb\trans\x)^2}\iint \frac{ K(v)K(u) E\{S_c( 
	z,\X_i)\mid\bb\trans\x^*\}
	f(z,\bb\trans\x^*) f_{\bb\trans\X}(\bb\trans\x^*)}{f_{\bb\trans\X}(\bb\trans\x)
	S(z,\bb\trans\x)
	E\{S_c( 
	z,\X_i)\mid\bb\trans\x\} }u^2dvdu,
\ese
where throughout the text, $z^*$ is on the line connecting $z$ and
$z+bv$. Thus, the absolute bias 
is
\bse
&&\left|E \left\{\frac{1}{n}\sum_{i=1}^n 
\frac{K_b(Z_i-z)\Delta_iK_h(\bb\trans\X_i-\bb\trans\x)}{f_{\bb\trans\X}(\bb\trans\x)E\{I(Z\ge 
	Z_i)\mid\bb\trans\x\}}\right\}
-\lambda(z,\bb\trans\x)
\right|\\
&\le&b^2\sup_{z^*,\bb\trans\x}\left|\frac{\partial^2}{2\partial z^2}
\frac{f(z^*,\bb\trans\x)}{
	S(z^*,\bb\trans\x)}
\right|\int v^2K(v)dv\\
&&+h^2\sup_{z,\bb\trans\x, \bb\trans\x^*}
\left|
\frac{\partial^2}{2\partial (\bb\trans\x)^2} \frac{E\{S_c( 
	z,\X_i)\mid\bb\trans\x^*\}
	f(z,\bb\trans\x^*) f_{\bb\trans\X}(\bb\trans\x^*)}{f_{\bb\trans\X}(\bb\trans\x)
	S(z,\bb\trans\x)
	E\{S_c( 
	z,\X_i)\mid\bb\trans\x\} }\right|
\int u^2K(u)du\\
&=&O(h^2+b^2)
\ese
under conditions \ref{assum:kernel}-\ref{assum:survivalfunction}. Following the same procedure, noting that
$A=O_p\{(nh)^{-1/2}+h^2\}$ uniformly, 
we can show that 
\bse
\frac{1}{n}\sum_{i=1}^n 
\frac{K_b(Z_i-z)\Delta_iK_h(\bb\trans\X_i-\bb\trans\x)}{f_{\bb\trans\X}(\bb\trans\x)E\{I(Z\ge 
	Z_i)\mid\bb\trans\x\}}O_p(A)
=O_p\{h^2+(nh)^{-1/2}\},
\ese
hence
the bias of $\wh\lambda(z,\bb\trans\x)$ is of order
$O_p\{(nh)^{-1/2}+h^2+b^2\}$ uniformly.
On the other hand, the variance of $\wh{\lambda}(z,\bb\trans\x)$ is
\bse
&&\var\left\{\wh{\lambda}(z,\bb\trans\x)\right\}\\
&=&\var\left[\frac{1}{n}\sumi K_b(Z_i-z) \frac{\Delta_iK_h(\bb\trans\X_i-\bb\trans\x)}{f_{\bb\trans\X}(\bb\trans\x)E\{I(Z\ge Z_i)\mid\bb\trans\x\}+A}\right]\\
&=&\var\left[\frac{1}{n}\sumi \frac{K_b(Z_i-z)\Delta_iK_h(\bb\trans\X_i-\bb\trans\x)}{f_{\bb\trans\X}(\bb\trans\x)E\{I(Z\ge Z_i)\mid\bb\trans\x\}}\{1+O_p(A)\}\right]\\
&\le&2\var\left[\frac{1}{n}\sumi\frac{K_b(Z_i-z)\Delta_iK_h(\bb\trans\X_i-\bb\trans\x)}{f_{\bb\trans\X}(\bb\trans\x)E\{I(Z\ge Z_i)\mid\bb\trans\x\}}\right]
+2\var\left[\frac{1}{n}\sumi\frac{K_b(Z_i-z)\Delta_iK_h(\bb\trans\X_i-\bb\trans\x)}{f_{\bb\trans\X}(\bb\trans\x)E\{I(Z\ge Z_i)\mid\bb\trans\x\}}O_p(A)\right].
\ese
We inspect the first term first.
\bse
&&2\var\left[\frac{1}{n}\sumi\frac{K_b(Z_i-z)\Delta_iK_h(\bb\trans\X_i-\bb\trans\x)}{f_{\bb\trans\X}(\bb\trans\x)E\{I(Z\ge Z_i)\mid\bb\trans\x\}}\right]\\
&=&\frac{2}{n}\var\left[ \frac{K_b(Z_i-z)\Delta_iK_h(\bb\trans\X_i-\bb\trans\x)}{f_{\bb\trans\X}(\bb\trans\x)E\{I(Z\ge Z_i)\mid\bb\trans\x\}}\right]\\
&=&\frac{2}{n}\left(E\left[\frac{K_b(Z_i-z)\Delta_iK_h(\bb\trans\X_i-\bb\trans\x)}{f_{\bb\trans\X}(\bb\trans\x)E\{I(Z\ge
	Z_i)\mid\bb\trans\x\}}\right]^2-\lambda^2(z,\bb\trans\x)+O(n^{-1/2}h^{-1/2}+h^2+b^2)\right)\\
&=&\frac{2}{n}\left(E\left[\frac{K_b(Z_i-z)\Delta_iK_h(\bb\trans\X_i-\bb\trans\x)}{f_{\bb\trans\X}(\bb\trans\x)E\{I(Z\ge Z_i)\mid\bb\trans\x\}}\right]^2\right)+O(1/n)\\
&=&\frac{2}{b^2h^2n}\iint\frac{K^2\{(Z_i-z)/b\}K^2\{(\bb\trans\x_i-\bb\trans\x)/h\}}{f_{\bb\trans\X}^2(\bb\trans\x)S^2(z_i,\bb\trans\x)E^2\{S_c(z_i,\X_i)\mid\bb\trans\x\}}\\
&&\times f(z_i,\bb\trans\x_i)S_c(z_i,\x_i)f_{\X\mid\bb\trans\X}(\x_i,\bb\trans\x_i)f_{\bb\trans\X}(\bb\trans\x_i)dz_id\x_id\bb\trans\x_i+O(1/n)\\
&=&\frac{2}{bhn}\iint\frac{K^2(v)K^2(u)}{f_{\bb\trans\X}^2(\bb\trans\x)S^2(z+bv,\bb\trans\x)E^2\{S_c(z+bv,\X_i)\mid\bb\trans\x\}}\\
&&\times f(z+bv,\bb\trans\x+hu)E\{S_c(z+bv,\X_i)\mid\bb\trans\x+hu\}f_{\bb\trans\X}(\bb\trans\x+hu)dvdu+O(1/n)\\
&=&\frac{2}{bhn}\iint\frac{f(z,\bb\trans\x)K^2(v)K^2(u)}{f_{\bb\trans\X}(\bb\trans\x)S^2(z,\bb\trans\x)E\{S_c(z,\X_i)\mid\bb\trans\x\}}dvdu\\
&&+\frac{b\partial^2}{nh\partial z^2}\iint\frac{f(z^*,\bb\trans\x)K^2(v)K^2(u)}{f_{\bb\trans\X}(\bb\trans\x)S^2(z^*,\bb\trans\x)E\{S_c(z^*,\X_i)\mid\bb\trans\x\}}v^2dvdu\\
&&+\frac{h\partial^2}{nb\partial
	(\bb\trans\x)^2}\iint\frac{f_{\bb\trans\X}(\bb\trans\x^*)f(z,\bb\trans\x^*)E\{S_c(z,\X_i)\mid\bb\trans\x^*\}K^2(v)K^2(u)}{f_{\bb\trans\X}^2(\bb\trans\x)S^2(z,\bb\trans\x)E\{S_c(z,\X_i)^2\mid\bb\trans\x\}}u^2dvdu+O(1/n)\\
&\le&\frac{2}{bhn}\frac{f(z,\bb\trans\x)}{f_{\bb\trans\X}(\bb\trans\x)S^2(z,\bb\trans\x)E\{S_c(z,\X_i)\mid\bb\trans\x\}}\left\{\int
K^2(u)du\right\}^2\\
&&+\frac{b}{nh}\sup_{z^*,\bb\trans\x}\left|\frac{\partial^2}{\partial z^2}\frac{f(z^*,\bb\trans\x)}{f_{\bb\trans\X}(\bb\trans\x)S^2(z^*,\bb\trans\x)E\{S_c(z^*,\X_i)\mid\bb\trans\x\}}\right|
\left\{\int K^2(u)u^2 du\right\}\left\{\int K^2(u) du\right\}\\
&&+\frac{h}{nb}\sup_{z,\bb\trans\x,\bb\trans\x^*}\left|\frac{\partial^2}{\partial
	(\bb\trans\x)^2}\frac{f_{\bb\trans\X}(\bb\trans\x^*)f(z,\bb\trans\x^*)E\{S_c(z,\X_i)\mid\bb\trans\x^*\}}{f_{\bb\trans\X}^2(\bb\trans\x)S^2(z,\bb\trans\x)E\{S_c(z,\X_i)^2\mid\bb\trans\x\}}\right|\\
&&\times\left\{\int K^2(u)u^2 du\right\}\left\{\int K^2(u) du\right\}
+O(1/n)\\
&=&O\{1/(nhb)+h/(nb)+b/(nh)+1/n\}\\
&=&O\{1/(nhb)\}
\ese
uniformly
under conditions \ref{assum:kernel}-\ref{assum:survivalfunction} and $\bb\trans\x^*$ is on the line connecting $\bb\trans\x$ and $\bb\trans\x+hu$.
For the second term
\bse
&&2\var\left[\frac{1}{n}\sumi\frac{K_b(Z_i-z)\Delta_iK_h(\bb\trans\X_i-\bb\trans\x)}{f_{\bb\trans\X}(\bb\trans\x)E\{I(Z\ge Z_i)\mid\bb\trans\x\}}O_p(A)\right]\\
&\le&2E\left[\frac{1}{n}\sumi\frac{K_b(Z_i-z)\Delta_iK_h(\bb\trans\X_i-\bb\trans\x)}{f_{\bb\trans\X}(\bb\trans\x)E\{I(Z\ge Z_i)\mid\bb\trans\x\}}O_p(A)\right]^2\\
&\le&2E\left(\frac{1}{n}\sumi\left[\frac{K_b(Z_i-z)\Delta_iK_h(\bb\trans\X_i-\bb\trans\x)}{f_{\bb\trans\X}(\bb\trans\x)E\{I(Z\ge Z_i)\mid\bb\trans\x\}}\right]^2\sup\{O_p^2(A)\}\right)\\
&=&2E\left(\frac{1}{n}\sumi\left[\frac{K_b(Z_i-z)\Delta_iK_h(\bb\trans\X_i-\bb\trans\x)}{f_{\bb\trans\X}(\bb\trans\x)E\{I(Z\ge Z_i)\mid\bb\trans\x\}}\right]^2\left[O_p\{(nh)^{-1}+h^4\}\right]\right)\\
&=&2E\left[\frac{K_b(Z_i-z)\Delta_iK_h(\bb\trans\X_i-\bb\trans\x)}{f_{\bb\trans\X}(\bb\trans\x)E\{I(Z\ge Z_i)\mid\bb\trans\x\}}\right]^2 O\{(nh)^{-1}+h^4\}\\
&=&\left[\frac{2}{bhn}\frac{f(z,\bb\trans\x)}{f_{\bb\trans\X}(\bb\trans\x)S^2(z,\bb\trans\x)E\{S_c(z,\X_i)\mid\bb\trans\x\}}\left\{\int K^2(u)du\right\}^2+O(n^{-1}b^{-1}h+n^{-1}h^{-1}b)\right]\\
&&\times O\{(nh)^{-1}+h^4\}\\
&=&O\{(nh)^{-2}b^{-1}+n^{-1}h^3b^{-1}\}
\ese
under conditions \ref{assum:kernel}-\ref{assum:survivalfunction} uniformly. Summarizing the above results, the
variance of $\wh\lambda(z,\bb\trans\x)$ is of order $1/(nhb)$
uniformly. 
Hence we have the consistency of estimator
$\widehat{\lambda}(z,\bb\trans\x)$, specifically 
\be
\widehat{\lambda}(z,\bb\trans\x)=\lambda(z,\bb\trans\x)+O_p\{(nhb)^{-1/2}+h^2+b^2\}\n
\ee
uniformly under condition \ref{assum:kernel}-\ref{assum:survivalfunction}.
Next we inspect the estimator for the first
derivative of hazard function $\lambda(z,\bb\trans\x)$.
Let 
\bse
\wh\blam_{11}&=&
-\sum_{i=1}^nK_b(Z_i-z) \frac{\Delta_iK_h'(\bb\trans\X_i-\bb\trans\x)}{\sum_j I(Z_j\ge
	z)K_h(\bb\trans\X_j-\bb\trans\x)}\\
\wh\blam_{12}&=&\sum_{i=1}^nK_b(Z_i-z)\Delta_iK_h(\bb\trans\X_i-\bb\trans\x) \frac{\sumj I(Z_j\ge z)K_h'(\bb\trans\X_j-\bb\trans\x)}{\{\sumj I(Z_j\ge
	z)K_h(\bb\trans\X_j-\bb\trans\x)\}^2}.
\ese
Then 
$\wh\blam_1(x,\bb\trans\x)=\wh\blam_{11}+\wh\blam_{12}$.
To analyze $\wh\blam_{11}$,
\bse
E\wh\blam_{11}&=&E\left\{-\sumi K_b(Z_i-z) \frac{\Delta_iK_h'(\bb\trans\X_i-\bb\trans\x)}{\sumj I(Z_j\ge
	z)K_h(\bb\trans\X_j-\bb\trans\x)}\right\}\\
&=&E\left[\frac{1}{n}\sumi-K_b(Z_i-z) \frac{\Delta_iK_h'(\bb\trans\X_i-\bb\trans\x)}{f_{\bb\trans\X}(\bb\trans\x)E\{I(Z\ge Z_i)\mid\bb\trans\x\}+A}\right]\\
&=&E\left[\frac{1}{n}\sumi \frac{-K_b(Z_i-z)\Delta_iK_h'(\bb\trans\X_i-\bb\trans\x)}{f_{\bb\trans\X}(\bb\trans\x)S(Z_i,\bb\trans\x)E\{S_c(Z_i,\X_i)\mid\bb\trans\x\}}\right]\\
&&+E\left[ \frac{1}{n}\sumi\frac{-K_b(Z_i-z)\Delta_iK_h'(\bb\trans\X_i-\bb\trans\x)}{f_{\bb\trans\X}(\bb\trans\x)S(Z_i,\bb\trans\x)E\{S_c(Z_i,\X_i)\mid\bb\trans\x\}}O_p(A)\right]\\
&=&E\left[ \frac{-K_b(Z_i-z)\Delta_iK_h'(\bb\trans\X_i-\bb\trans\x)}{f_{\bb\trans\X}(\bb\trans\x)S(Z_i,\bb\trans\x)E\{S_c(Z_i,\X_i)\mid\bb\trans\x\}}\right]\\
&&+E\left[ \frac{1}{n}\sumi\frac{-K_b(Z_i-z)\Delta_iK_h'(\bb\trans\X_i-\bb\trans\x)}{f_{\bb\trans\X}(\bb\trans\x)S(Z_i,\bb\trans\x)E\{S_c(Z_i,\X_i)\mid\bb\trans\x\}}O_p(A)\right].
\ese
We inspect the first term to obtain
\bse
&&E\left[ \frac{-K_b(Z_i-z)\Delta_iK_h'(\bb\trans\X_i-\bb\trans\x)}{f_{\bb\trans\X}(\bb\trans\x)S(Z_i,\bb\trans\x)E\{S_c(Z_i,\X_i)\mid\bb\trans\x\}}\right]\\
&=&-\iint
\frac{K_b(z_i-z)K_h'(\bb\trans\x_i-\bb\trans\x)}{f_{\bb\trans\X}(\bb\trans\x)
	S(z_i,\bb\trans\x)E\{S_c( 
	z_i,\X_i)\mid\bb\trans\x\}}
f(z_i,\bb\trans\x_i)S_c(z_i,\x_i)f_{\X\mid\bb\trans\X}(\x_i,\bb\trans\x_i)\\
&&\times f_{\bb\trans\X}(\bb\trans\x_i)dz_id\x_id\bb\trans\x_i\\
&=&-\iint
\frac{1}{h}\frac{K(v)K'(u)f(z+bv,\bb\trans\x+hu)E\{S_c(z+bv,\X_i)\mid\bb\trans\x+hu\}}{f_{\bb\trans\X}(\bb\trans\x)
	S(z+bv,\bb\trans\x)E\{S_c( 
	z+bv,\X_i)\mid\bb\trans\x\}}f_{\bb\trans\X}(\bb\trans\x+hu)dvdu\\
&=&-\iint
\frac{\partial\left[f(z,\bb\trans\x)E\{S_c(z,\X_i)\mid\bb\trans\x\}f_{\bb\trans\X}(\bb\trans\x)\right]/\partial\bb\trans\x}{f_{\bb\trans\X}(\bb\trans\x)
	S(z,\bb\trans\x)E\{S_c( 
	z,\X_i)\mid\bb\trans\x\}}K(v)uK'(u)dvdu\\
&&-\frac{b^2\partial^3}{2\partial z^2\partial\bb\trans\x}\iint
\frac{f_{\bb\trans\X}(\bb\trans\x^*)f(z^*,\bb\trans\x^*)E\{S_c(z^*,\X_i)\mid\bb\trans\x^*\}}{f_{\bb\trans\X}(\bb\trans\x)
	S(z^*,\bb\trans\x)E\{S_c( 
	z^*,\X_i)\mid\bb\trans\x\}}v^2K(v)uK'(u)dvdu\\
&&-\frac{h^2\partial^3}{6\partial(\bb\trans\x)^3}\iint
\frac{f_{\bb\trans\X}(\bb\trans\x^*)f(z,\bb\trans\x^*)E\{S_c(z,\X_i)\mid\bb\trans\x^*\}}{f_{\bb\trans\X}(\bb\trans\x)
	S(z,\bb\trans\x)E\{S_c( 
	z,\X_i)\mid\bb\trans\x\}}K(v)u^3K'(u)dvdu\\
&=&\frac{\partial\left[f(z,\bb\trans\x)E\{S_c(z,\X_i)\mid\bb\trans\x\}f_{\bb\trans\X}(\bb\trans\x)\right]/\partial\bb\trans\x}{f_{\bb\trans\X}(\bb\trans\x)
	S(z,\bb\trans\x)E\{S_c( 
	z,\X_i)\mid\bb\trans\x\}}\\
&&-\frac{b^2\partial^3}{2\partial z^2\partial\bb\trans\x}\iint
\frac{f_{\bb\trans\X}(\bb\trans\x^*)f(z^*,\bb\trans\x^*)E\{S_c(z^*,\X_i)\mid\bb\trans\x^*\}}{f_{\bb\trans\X}(\bb\trans\x)
	S(z^*,\bb\trans\x)E\{S_c( 
	z^*,\X_i)\mid\bb\trans\x\}}v^2K(v)uK'(u)dvdu\\
&&-\frac{h^2\partial^3}{6\partial(\bb\trans\x)^3}\iint
\frac{f_{\bb\trans\X}(\bb\trans\x^*)f(z,\bb\trans\x^*)E\{S_c(z,\X_i)\mid\bb\trans\x^*\}}{f_{\bb\trans\X}(\bb\trans\x)
	S(z,\bb\trans\x)E\{S_c( 
	z,\X_i)\mid\bb\trans\x\}}K(v)u^3K'(u)dvdu.
\ese
Hence the absolute bias is given by
\bse
&&\left|E\left[ \frac{-K_b(Z_i-z)\Delta_iK_h'(\bb\trans\X_i-\bb\trans\x)}{f_{\bb\trans\X}(\bb\trans\x)S(Z_i,\bb\trans\x)E\{S_c(Z_i,\X_i)\mid\bb\trans\x\}}\right]\right.\\
&&\left.-\frac{\partial\left[f(z,\bb\trans\x)E\{S_c(z,\X_i)\mid\bb\trans\x\}f_{\bb\trans\X}(\bb\trans\x)\right]/\partial\bb\trans\x}{f_{\bb\trans\X}(\bb\trans\x)
	S(z,\bb\trans\x)E\{S_c( 
	z,\X_i)\mid\bb\trans\x\}}\right|\\
&\le&b^2\sup_{z^*,\bb\trans\x,\bb\trans\x^*}\left|\frac{\partial^3}{2\partial z^2\partial\bb\trans\x}\frac{f_{\bb\trans\X}(\bb\trans\x^*)f(z^*,\bb\trans\x^*)E\{S_c(z^*,\X_i)\mid\bb\trans\x^*\}}{f_{\bb\trans\X}(\bb\trans\x)
	S(z^*,\bb\trans\x)E\{S_c( 
	z^*,\X_i)\mid\bb\trans\x\}}\right|\int v^2K(v)dv\\
&&+h^2\sup_{z,\bb\trans\x,\bb\trans\x^*}\left|\frac{\partial^3}{2\partial(\bb\trans\x)^3}\frac{f_{\bb\trans\X}(\bb\trans\x^*)f(z,\bb\trans\x^*)E\{S_c(z,\X_i)\mid\bb\trans\x^*\}}{f_{\bb\trans\X}(\bb\trans\x)
	S(z,\bb\trans\x)E\{S_c( 
	z,\X_i)\mid\bb\trans\x\}}\right|\int u^2K(u)du\\
&=&O(b^2+h^2)
\ese
uniformly under condition \ref{assum:kernel}-\ref{assum:survivalfunction}.

Following the same procedure, we conclude that
\bse
\frac{1}{n}\sumi\frac{-K_b(Z_i-z)\Delta_iK_h'(\bb\trans\X_i-\bb\trans\x)}{f_{\bb\trans\X}(\bb\trans\x)S(Z_i,\bb\trans\x)E\{S_c(Z_i,\X_i)\mid\bb\trans\x\}}O_p(A)=O_p\{h^2+(nh)^{-1/2}\}
\ese
uniformly under conditions \ref{assum:kernel}-\ref{assum:survivalfunction} due to $A=O_p\{h^2+(nh)^{-1/2}\}$. Therefore, we have
\bse
E\wh\blam_{11}&=&\frac{\partial\left[f(z,\bb\trans\x)E\{S_c(z,\X_i)\mid\bb\trans\x\}f_{\bb\trans\X}(\bb\trans\x)\right]/\partial\bb\trans\x}{f_{\bb\trans\X}(\bb\trans\x)
	S(z,\bb\trans\x)E\{S_c( 
	z,\X_i)\mid\bb\trans\x\}}+O\{(nh)^{-1/2}+b^2+h^2\}
\ese
For $\wh\blam_{12}$, let 
$B=	-1/n\sumj I(Z_j\ge z)  K_h'(\bb\trans\X_j-\bb\trans\x)
-\partial
f_{\bb\trans\X}(\bb\trans\x)E\{I(Z_j\ge z)\mid\bb\trans\x\}/\partial \bb\trans\x$,
then
\bse
\wh\blam_{12}&=&K_b(Z_i-z)\Delta_iK_h(\bb\trans\X_i-\bb\trans\x) \frac{\sumj I(Z_j\ge
	z)K_h'(\bb\trans\X_j-\bb\trans\x)}{\{\sumj I(Z_j\ge
	z)K_h(\bb\trans\X_j-\bb\trans\x)\}^2}\\
&=&-\frac{1}{n}\sumi K_b(Z_i-z)\Delta_iK_h(\bb\trans\X_i-\bb\trans\x) \frac{\partial\left[
	f_{\bb\trans\X}(\bb\trans\x)E\{I(Z\ge Z_i)\mid\bb\trans\x\}\right]/ \partial \bb\trans\x+B}{[f_{\bb\trans\X}(\bb\trans\X_i)E\{I(Z\ge Z_i)\mid\bb\trans\x\}+A]^2}\\
&=&-\frac{1}{n}\sumi K_b(Z_i-z)\Delta_iK_h(\bb\trans\X_i-\bb\trans\x) \frac{\partial\left[
	f_{\bb\trans\X}(\bb\trans\z)E\{I(Z\ge Z_i)\mid\bb\trans\x\}\right]/ \partial \bb\trans\x}{f_{\bb\trans\X}^2(\bb\trans\x)E^2\{I(Z\ge Z_i)\mid\bb\trans\x\}}\\
&&\times\left\{1+O_p(B)+O_p(A)\right\}.
\ese
We have
\bse
&&E\left[-\frac{1}{n}\sumi K_b(Z_i-z)\Delta_iK_h(\bb\trans\X_i-\bb\trans\x) \frac{\partial\left[
	f_{\bb\trans\X}(\bb\trans\z)E\{I(Z\ge Z_i)\mid\bb\trans\x\}\right]/ \partial \bb\trans\x}{f_{\bb\trans\X}^2(\bb\trans\x)E^2\{I(Z\ge Z_i)\mid\bb\trans\x\}}\right]\\
&=&-E\left[ K_b(Z_i-z)\Delta_iK_h(\bb\trans\X_i-\bb\trans\x) \frac{\partial\left[
	f_{\bb\trans\X}(\bb\trans\z)E\{I(Z\ge Z_i)\mid\bb\trans\x\}\right]/ \partial \bb\trans\x}{f_{\bb\trans\X}^2(\bb\trans\x)E^2\{I(Z\ge Z_i)\mid\bb\trans\x\}}\right]\\
&=&-\iint K_b(z_i-z)K_h(\bb\trans\x_i-\bb\trans\x) \frac{\partial
	\left[f_{\bb\trans\X}(\bb\trans\x)S(z_i,\bb\trans\x)E\{S_c( 
	z_i,\X_i)\mid\bb\trans\x\}\right]/ \partial \bb\trans\x}{f_{\bb\trans\X}^2(\bb\trans\x)S^2(z_i,\bb\trans\x)E^2\{S_c( 
	z_i,\X_i)\mid\bb\trans\x\}}\\
&&\times f(z_i,\bb\trans\x_i)S_c(z_i,\x_i)f_{\X\mid\bb\trans\X}(\x_i,\bb\trans\x_i)f_{\bb\trans\X}(\bb\trans\x_i)dz_id\x_id\bb\trans\x_i\\
&=&-\iint K_b(z_i-z)K_h(\bb\trans\x_i-\bb\trans\x) \frac{\partial
	\left[f_{\bb\trans\X}(\bb\trans\x)S(z_i,\bb\trans\x)E\{S_c(z_i,\X_i)\mid\bb\trans\x\}\right]/ \partial \bb\trans\x}{f_{\bb\trans\X}^2(\bb\trans\x)S^2(z_i,\bb\trans\x)E^2\{S_c( 
	z_i,\X_i)\mid\bb\trans\x\}}\\
&&\times f(z_i,\bb\trans\x_i)E\{S_c(z_i,\X_i)\mid\bb\trans\x_i\}f_{\bb\trans\X}(\bb\trans\x_i)dz_id\bb\trans\x_i\\
&=&-\iint K(v)K(u) \frac{\partial
	\left[f_{\bb\trans\X}(\bb\trans\x)S(z+bv,\bb\trans\x)E\{S_c( 
	z+bv,\X_i)\mid\bb\trans\x\}\right]/ \partial \bb\trans\x}{f_{\bb\trans\X}^2(\bb\trans\x)S^2(z+bv,\bb\trans\x)E^2\{S_c( 
	z+bv,\X_i)\mid\bb\trans\x\}}\\
&&\times f(z+bv,\bb\trans\x+hu)E\{S_c(z+bv,\X_i)\mid\bb\trans\x+hu\}f_{\bb\trans\X}(\bb\trans\x+hu)dvdu\\
&=&-\iint K(v)K(u) \frac{\partial
	\left[f_{\bb\trans\X}(\bb\trans\x)S(z,\bb\trans\x)E\{S_c( 
	z,\X_i)\mid\bb\trans\x\}\right]/ \partial \bb\trans\x}{f_{\bb\trans\X}(\bb\trans\x)S^2(z,\bb\trans\x)E\{S_c( 
	z,\X_i)\mid\bb\trans\x\}}f(z,\bb\trans\x)dvdu\\
&&-\frac{b^2\partial^2}{2\partial z^2}\iint K(v)K(u) \frac{\partial
	\left[f_{\bb\trans\X}(\bb\trans\x)S(z^*,\bb\trans\x)E\{S_c( 
	z^*,\X_i)\mid\bb\trans\x\}\right]/ \partial \bb\trans\x}{f_{\bb\trans\X}(\bb\trans\x)S^2(z^*,\bb\trans\x)E\{S_c( 
	z^*,\X_i)\mid\bb\trans\x\}}f(z^*,\bb\trans\x)v^2dvdu\\
&&-\frac{h^2\partial^2}{2\partial (\bb\trans\x)^2}\iint K(v)K(u) \frac{\partial
	\left[f_{\bb\trans\X}(\bb\trans\x)S(z,\bb\trans\x)E\{S_c( 
	z,\X_i)\mid\bb\trans\x\}\right]/ \partial \bb\trans\x}{f_{\bb\trans\X}^2(\bb\trans\x)S^2(z,\bb\trans\x)E^2\{S_c( 
	z,\X_i)\mid\bb\trans\x\}}\\
&&\times f(z,\bb\trans\x^*)E\{S_c(z,\X_i)\mid\bb\trans\x^*\}f_{\bb\trans\X}(\bb\trans\x^*)u^2dvdu\\
&=&-\frac{\partial
	\left[f_{\bb\trans\X}(\bb\trans\x)S(z,\bb\trans\x)E\{S_c( 
	z,\X_i)\mid\bb\trans\x\}\right]/ \partial \bb\trans\x}{f_{\bb\trans\X}(\bb\trans\x)S^2(z,\bb\trans\x)E\{S_c( 
	z,\X_i)\mid\bb\trans\x\}}f(z,\bb\trans\x)\\
&&-\frac{b^2\partial^2}{2\partial z^2}\iint  \frac{\partial
	\left[f_{\bb\trans\X}(\bb\trans\x)S(z^*,\bb\trans\x)E\{S_c( 
	z^*,\X_i)\mid\bb\trans\x\}\right]/ \partial \bb\trans\x}{f_{\bb\trans\X}(\bb\trans\x)S^2(z^*,\bb\trans\x)E\{S_c( 
	z^*,\X_i)\mid\bb\trans\x\}}f(z^*,\bb\trans\x)v^2K(v)K(u)dvdu\\
&&-\frac{h^2\partial^2}{2\partial (\bb\trans\x)^2}\iint  \frac{\partial
	\left[f_{\bb\trans\X}(\bb\trans\x)S(z,\bb\trans\x)E\{S_c( 
	z,\X_i)\mid\bb\trans\x\}\right]/ \partial \bb\trans\x}{f_{\bb\trans\X}^2(\bb\trans\x)S^2(z,\bb\trans\x)E^2\{S_c( 
	z,\X_i)\mid\bb\trans\x\}}\\
&&\times f(z,\bb\trans\x^*)E\{S_c(z,\X_i)\mid\bb\trans\x^*\}f_{\bb\trans\X}(\bb\trans\x^*)u^2K(v)K(u)dvdu.
\ese
Therefore
\bse
&&\left|E\left[-\frac{1}{n}\sumi K_b(Z_i-z)\Delta_iK_h(\bb\trans\X_i-\bb\trans\x) \frac{\partial\left[
	f_{\bb\trans\X}(\bb\trans\z)E\{I(Z\ge Z_i)\mid\bb\trans\x\}\right]/ \partial \bb\trans\x}{f_{\bb\trans\X}^2(\bb\trans\x)E^2\{I(Z\ge Z_i)\mid\bb\trans\x\}}\right]\right.\\
&&\left.+\frac{\partial
	\left[f_{\bb\trans\X}(\bb\trans\x)S(z,\bb\trans\x)E\{S_c( 
	z,\X_i)\mid\bb\trans\x\}\right]/ \partial \bb\trans\x}{f_{\bb\trans\X}(\bb\trans\x)S^2(z,\bb\trans\x)E\{S_c( 
	z,\X_i)\mid\bb\trans\x\}}f(z,\bb\trans\x)\right|\\
&\le&b^2\sup_{z^*,\bb\trans\x}\left|\frac{\partial^2}{2\partial z^2}\frac{\partial
	\left[f_{\bb\trans\X}(\bb\trans\x)S(z^*,\bb\trans\x)E\{S_c( 
	z^*,\X_i)\mid\bb\trans\x\}\right]/ \partial \bb\trans\x}{f_{\bb\trans\X}(\bb\trans\x)S^2(z^*,\bb\trans\x)E\{S_c( 
	z^*,\X_i)\mid\bb\trans\x\}}f(z^*,\bb\trans\x)\right|\\
&&\times \left\{\int v^2K(v)dv\right\}\\
&&+h^2\sup_{z,\bb\trans\x,\bb\trans\x^*}\left|\frac{\partial^2}{2\partial (\bb\trans\x)^2} \frac{\partial
	\left[f_{\bb\trans\X}(\bb\trans\x)S(z,\bb\trans\x)E\{S_c( 
	z,\X_i)\mid\bb\trans\x\}\right]/ \partial \bb\trans\x}{f_{\bb\trans\X}^2(\bb\trans\x)S^2(z,\bb\trans\x)E^2\{S_c( 
	z,\X_i)\mid\bb\trans\x\}}\right.\\
&&\times
f(z,\bb\trans\x^*)E\{S_c(,\X_i)\mid\bb\trans\x^*\}f_{\bb\trans\X}(\bb\trans\x^*)\Bigg|\left\{\int
u^2K(u)du\right\}\\
&=&O(b^2+h^2)
\ese
uniformly under conditions \ref{assum:kernel}-\ref{assum:survivalfunction}. 
Noting that $B=O_p(n^{-1/2}h^{-3/2}+h^2)$, based on similar 
procedure, we have 
\bse
&&-\frac{1}{n}\sumi K_b(Z_i-z)\Delta_iK_h(\bb\trans\X_i-\bb\trans\x) \frac{\partial\left[
	f_{\bb\trans\X}(\bb\trans\z)E\{I(Z\ge
	Z_i)\mid\bb\trans\x\}\right]/ \partial
	\bb\trans\x}{f_{\bb\trans\X}^2(\bb\trans\x)E^2\{I(Z\ge
	Z_i)\mid\bb\trans\x\}}O_p(B)\\
&=&O_p(n^{-1/2}h^{-3/2}+h^2),\\
&&-\frac{1}{n}\sumi K_b(Z_i-z)\Delta_iK_h(\bb\trans\X_i-\bb\trans\x) \frac{\partial\left[
	f_{\bb\trans\X}(\bb\trans\z)E\{I(Z\ge
	Z_i)\mid\bb\trans\x\}\right]/ \partial
	\bb\trans\x}{f_{\bb\trans\X}^2(\bb\trans\x)E^2\{I(Z\ge
	Z_i)\mid\bb\trans\x\}}O_p(A)\\
&=&O_p(n^{-1/2}h^{-1/2}+h^2)
\ese
uniformly under condition \ref{assum:kernel}-\ref{assum:survivalfunction}. Therefore, we can conclude that
\bse
E\wh\blam_{12}=-\frac{\partial
	\left[f_{\bb\trans\X}(\bb\trans\x)S(z,\bb\trans\x)E\{S_c( 
	z,\X_i)\mid\bb\trans\x\}\right]/ \partial \bb\trans\x}{f_{\bb\trans\X}(\bb\trans\x)S^2(z,\bb\trans\x)E\{S_c( 
	z,\X_i)\mid\bb\trans\x\}}f(z,\bb\trans\x)+O(n^{-1/2}h^{-3/2}+b^2+h^2)
\ese
In addition, we have
\bse
&&\frac{\partial\left[f(z,\bb\trans\x)E\{S_c(z,\X_i)\mid\bb\trans\x\}f_{\bb\trans\X}(\bb\trans\x)\right]/\partial\bb\trans\x}{f_{\bb\trans\X}(\bb\trans\x)
	S(z,\bb\trans\x)E\{S_c( 
	z,\X_i)\mid\bb\trans\x\}}\\
&&-\frac{\partial
	\left[f_{\bb\trans\X}(\bb\trans\x)S(z,\bb\trans\x)E\{S_c( 
	z,\X_i)\mid\bb\trans\x\}\right]/ \partial \bb\trans\x}{f_{\bb\trans\X}(\bb\trans\x)S^2(z,\bb\trans\x)E\{S_c( 
	z,\X_i)\mid\bb\trans\x\}}f(z,\bb\trans\x)\\
&=&\frac{\partial f(z,\bb\trans\x)/\partial\bb\trans\x}{S(z,\bb\trans\x)}
+\frac{f(z,\bb\trans\x)\partial E\{S_c(z,\X_i)\mid\bb\trans\x\}/\partial\bb\trans\x}{	S(z,\bb\trans\x)E\{S_c( 
	z,\X_i)\mid\bb\trans\x\}}
+\frac{f_{\bb\trans\X}'(\bb\trans\x)f(z,\bb\trans\x)}{f_{\bb\trans\X}(\bb\trans\x)
	S(z,\bb\trans\x)}\\
&&-\frac{f_{\bb\trans\X}'(\bb\trans\x)f(z,\bb\trans\x)}{f_{\bb\trans\X}(\bb\trans\x)S(z,\bb\trans\x)}
-\frac{f(z,\bb\trans\x)\partial S(z,\bb\trans\x)/\partial \bb\trans\x}{S^2(z,\bb\trans\x)}
-\frac{f(z,\bb\trans\x)\partial E\{S_c( 
	z,\X_i)\mid\bb\trans\x\}/\partial \bb\trans\x}{S(z,\bb\trans\x)E\{S_c( 
	z,\X_i)\mid\bb\trans\x\}}\\
&=&\frac{\partial
	f(z,\bb\trans\x)/\partial\bb\trans\x}{S(z,\bb\trans\x)}
-\frac{f(z,\bb\trans\x)\partial S(z,\bb\trans\x)/\partial \bb\trans\x}{S^2(z,\bb\trans\x)}\\
&=&\frac{\partial}{\partial\bb\trans\x}\lambda(z,\bb\trans\x)\\
&=&\blam_1(z,\bb\trans\x).
\ese
Combining $E\wh\blam_{11}$ and $E\wh\blam_{12}$, we readily obtain
\bse
\left|E\wh\blam_1(z,\bb\trans\x)-\frac{\partial}{\partial\bb\trans\x}\blam_1(z,\bb\trans\x)\right|
&=&\left|E\wh\blam_{11}+E\wh\blam_{12}-\frac{\partial}{\partial\bb\trans\x}\blam_1(z,\bb\trans\x)\right|\\
&=&O(n^{-1/2}h^{-3/2}+b^2+h^2)
\ese
uniformly under conditions \ref{assum:kernel}-\ref{assum:survivalfunction}.

The variance of $\wh\blam_1(z,\bb\trans\x)$ is given by
\bse
\var\left\{\wh\blam_1(z,\bb\trans\x)\right\}=\var\left\{\wh\blam_{11}+\wh\blam_{12}\right\}\le2\var(\wh\blam_{11})+2\var(\wh\blam_{12}).
\ese
We exam each term separately.
\bse
&&\var(\wh\blam_{11})\\
&=&\var\left\{\sum_{i=1}^nK_b(Z_i-z) \frac{-\Delta_iK_h'(\bb\trans\X_i-\bb\trans\x)}{\sum_j I(Z_j\ge
	z)K_h(\bb\trans\X_j-\bb\trans\x)}\right\}\\
&=&\var\left\{\frac{1}{n}\sumi  \frac{-K_b(Z_i-z)\Delta_iK_h'(\bb\trans\X_i-\bb\trans\x)}{f_{\bb\trans\X}(\bb\trans\x)S(Z_i,\bb\trans\x)E\{S_c(Z_i,\X_i)\mid\bb\trans\x\}+A}\right\}\\
&\le&\frac{2}{n}\var\left[\frac{-K_b(Z_i-z)\Delta_iK_h'(\bb\trans\X_i-\bb\trans\x)}{f_{\bb\trans\X}(\bb\trans\x)S(Z_i,\bb\trans\x)E\{S_c(Z_i,\X_i)\mid\bb\trans\x\}}\right]\\
&&+2\var\left[\frac{1}{n}\sumi\frac{-K_b(Z_i-z)\Delta_iK_h'(\bb\trans\X_i-\bb\trans\x)}{f_{\bb\trans\X}(\bb\trans\x)S(Z_i,\bb\trans\x)E\{S_c(Z_i,\X_i)\mid\bb\trans\x\}}O_p(A)\right]\\
\ese
The first part is given by
\bse
&&\frac{2}{n}\var\left[\frac{-K_b(Z_i-z)\Delta_iK_h'(\bb\trans\X_i-\bb\trans\x)}{f_{\bb\trans\X}(\bb\trans\x)S(Z_i,\bb\trans\x)E\{S_c(Z_i,\X_i)\mid\bb\trans\x\}}\right]\\
&=&\frac{2}{n}E\left[\frac{-K_b(Z_i-z)\Delta_iK_h'(\bb\trans\X_i-\bb\trans\x)}{f_{\bb\trans\X}(\bb\trans\x)S(Z_i,\bb\trans\x)E\{S_c(Z_i,\X_i)\mid\bb\trans\x\}}\right]^2\\
&&-\frac{2}{n}\left(E\left[\frac{-K_b(Z_i-z)\Delta_iK_h'(\bb\trans\X_i-\bb\trans\x)}{f_{\bb\trans\X}(\bb\trans\x)S(Z_i,\bb\trans\x)E\{S_c(Z_i,\X_i)\mid\bb\trans\x\}}\right]\right)^2\\
&=&\frac{2}{n}\iint\frac{K_b^2(z_i-z)K_h'^2(\bb\trans\x_i-\bb\trans\x)}{f_{\bb\trans\X}^2(\bb\trans\x)S^2(z_i,\bb\trans\x)E^2\{S_c(z_i,\X_i)\mid\bb\trans\x\}}\\
&&\times f(z_i,\bb\trans\x_i)S_c(z_i,\x_i)f_{\X\mid\bb\trans\X}(\x_i,\bb\trans\x_i)f_{\bb\trans\X}(\bb\trans\x_i)dz_id\x_id\bb\trans\x_i+O(1/n)\\
&=&\frac{2}{nbh^3}\iint\frac{K^2(v)K'^2(u)E\{S_c(z+bv,\X_i)\mid\bb\trans\x+hu\}}
{f_{\bb\trans\X}^2(\bb\trans\x)S^2(z+bv,\bb\trans\x)E^2\{S_c(z+bv,\X_i)\mid\bb\trans\x\}}\\
&&\times f(z+bv,\bb\trans\x+hu)f_{\bb\trans\X}(\bb\trans\x+hu)dvdu+O(1/n)\\
&=&\frac{2}{nbh^3}\iint\frac{K^2(v)K'^2(u)f(z,\bb\trans\x)}
{f_{\bb\trans\X}(\bb\trans\x)S^2(z,\bb\trans\x)E\{S_c(z,\X)\mid\bb\trans\x\}}dvdu\\
&&+\frac{b}{nh^3}\frac{\partial^2}{\partial z^2}\iint\frac{f(z^*,\bb\trans\x)}
{f_{\bb\trans\X}(\bb\trans\x)S^2(z^*,\bb\trans\x)E\{S_c(z^*,\X_i)\mid\bb\trans\x\}}v^2K^2(v)K'^2(u)dvdu\\
&&+\frac{1}{nbh}\frac{\partial^2}{\partial (\bb\trans\x)^2}\iint\frac{f(z,\bb\trans\x^*)f_{\bb\trans\X}(\bb\trans\x^*)E\{S_c(z,\X_i)\mid\bb\trans\x^*\}}
{f_{\bb\trans\X}^2(\bb\trans\x)S^2(z,\bb\trans\x)E^2\{S_c(z,\X_i)\mid\bb\trans\x\}}K^2(v)u^2K'^2(u)dvdu+O(1/n)\\
&\le&\frac{2}{nbh^3}\frac{f(z,\bb\trans\x)}
{f_{\bb\trans\X}(\bb\trans\x)S^2(z,\bb\trans\x)E\{S_c(z,\X)\mid\bb\trans\x\}}\left\{\int K^2(v)dv\right\}\left\{\int K'^2(u)du\right\}\\
&&+\frac{b}{nh^3}\sup_{z^*,\bb\trans\x}\left[\frac{\partial^2}{\partial z^2}\frac{f(z^*,\bb\trans\x)}
{f_{\bb\trans\X}(\bb\trans\x)S^2(z^*,\bb\trans\x)E\{S_c(z^*,\X_i)\mid\bb\trans\x\}}\right]\left\{\int v^2K^2(v)dv\right\}\left\{\int K'^2(u)du\right\}\\
&&+\frac{1}{nbh}\sup_{z,\bb\trans\x,\bb\trans\x^*}\left[\frac{\partial^2}{\partial (\bb\trans\x)^2}\frac{f(z,\bb\trans\x^*)f_{\bb\trans\X}(\bb\trans\x^*)E\{S_c(z,\X_i)\mid\bb\trans\x^*\}}
{f_{\bb\trans\X}^2(\bb\trans\x)S^2(z,\bb\trans\x)E^2\{S_c(z,\X_i)\mid\bb\trans\x\}}\right]\\
&&\qquad\times\left\{\int K^2(v)dv\right\}\left\{\int u^2K'^2(u)du\right\}+O(1/n)\\
&=&O\{1/(nbh^3)+b/(nh^3)+1/(nbh)+1/n\}\\
&=&O\{1/(nbh^3)\}
\ese
uniformly under condition \ref{assum:kernel}-\ref{assum:survivalfunction}.
Noting that $A=O_p(n^{-1/2}h^{-1/2}+h^2)$, the second part  is
\bse
&&2\var\left[\frac{1}{n}\sumi\frac{-K_b(Z_i-z)\Delta_iK_h'(\bb\trans\X_i-\bb\trans\x)}{f_{\bb\trans\X}(\bb\trans\x)S(Z_i,\bb\trans\x)E\{S_c(Z_i,\X_i)\mid\bb\trans\x\}}O_p(A)\right]\\
&\le&2E\left[\frac{1}{n}\sumi\frac{-K_b(Z_i-z)\Delta_iK_h'(\bb\trans\X_i-\bb\trans\x)}{f_{\bb\trans\X}(\bb\trans\x)S(Z_i,\bb\trans\x)E\{S_c(Z_i,\X_i)\mid\bb\trans\x\}}O_p(A)\right]^2\\
&\le&2E\left(\frac{1}{n}\sumi\left[\frac{-K_b(Z_i-z)\Delta_iK_h'(\bb\trans\X_i-\bb\trans\x)}{f_{\bb\trans\X}(\bb\trans\x)S(Z_i,\bb\trans\x)E\{S_c(Z_i,\X_i)\mid\bb\trans\x\}}\right]^2\sup \left\{O_p^2(A)\right\}\right)\\
&=&2E\left(\frac{1}{n}\sumi\left[\frac{-K_b(Z_i-z)\Delta_iK_h'(\bb\trans\X_i-\bb\trans\x)}{f_{\bb\trans\X}(\bb\trans\x)S(Z_i,\bb\trans\x)E\{S_c(Z_i,\X_i)\mid\bb\trans\x\}}\right]^2\right) O\{1/(nh)+h^4\}\\
&=&\left[\frac{2}{nbh^3}\frac{f(z,\bb\trans\x)}
{f_{\bb\trans\X}(\bb\trans\x)S^2(z,\bb\trans\x)E\{S_c(z,\X)\mid\bb\trans\x\}}\left\{\int K^2(v)dv\right\}\left\{\int K'^2(u)du\right\}+O\left(\frac{b}{nh^3}+\frac{1}{nbh}\right)\right]\\
&&\times O\{1/(nh)+h^4\}\\
&=& O\{n^{-2}b^{-1}h^{-4}+(nb)^{-1}h\}
\ese
uniformly under conditions \ref{assum:kernel}-\ref{assum:survivalfunction}.
Therefore 
\bse
\var(\wh\blam_{11})=O\{1/(nbh^3)\}
\ese
uniformly under conditions \ref{assum:kernel}-\ref{assum:survivalfunction}.

For $\wh\blam_{12}$,
\bse
&&\var(\wh\blam_{12})\\
&=&\var\left[\sum_{i=1}^nK_b(Z_i-z)\Delta_iK_h(\bb\trans\X_i-\bb\trans\x) \frac{\sumj I(Z_j\ge z)K_h'(\bb\trans\X_j-\bb\trans\x)}{\{\sumj I(Z_j\ge
	z)K_h(\bb\trans\X_j-\bb\trans\x)\}^2}\right]\\
&=&\var\left[\frac{1}{n}\sum_{i=1}^nK_b(Z_i-z)\Delta_iK_h(\bb\trans\X_i-\bb\trans\x) \frac{\partial\left[
	f_{\bb\trans\X}(\bb\trans\z)E\{I(Z\ge Z_i)\mid\bb\trans\x\}\right]/ \partial \bb\trans\x+B}{\{f_{\bb\trans\X}^2(\bb\trans\x)E^2\{I(Z\ge Z_i)\mid\bb\trans\x\}+A\}^2}\right]\\
&\le&\frac{2}{n}\var\left[K_b(Z_i-z)\Delta_iK_h(\bb\trans\X_i-\bb\trans\x) \frac{\partial\left[
	f_{\bb\trans\X}(\bb\trans\z)E\{I(Z\ge Z_i)\mid\bb\trans\x\}\right]/ \partial \bb\trans\x}{f_{\bb\trans\X}^2(\bb\trans\x)E^2\{I(Z\ge Z_i)\mid\bb\trans\x\}}\right]\\
&&+4\var\left[\frac{1}{n}\sumi K_b(Z_i-z)\Delta_iK_h(\bb\trans\X_i-\bb\trans\x) \frac{\partial\left[
	f_{\bb\trans\X}(\bb\trans\z)E\{I(Z\ge Z_i)\mid\bb\trans\x\}\right]/ \partial \bb\trans\x}{f_{\bb\trans\X}^2(\bb\trans\x)E^2\{I(Z\ge Z_i)\mid\bb\trans\x\}}O_p(B)\right]\\
&&+4\var\left[\frac{1}{n}\sumi K_b(Z_i-z)\Delta_iK_h(\bb\trans\X_i-\bb\trans\x) \frac{\partial\left[
	f_{\bb\trans\X}(\bb\trans\z)E\{I(Z\ge Z_i)\mid\bb\trans\x\}\right]/ \partial \bb\trans\x}{f_{\bb\trans\X}^2(\bb\trans\x)E^2\{I(Z\ge Z_i)\mid\bb\trans\x\}}O_p(A)\right].
\ese
The first part is given by
\bse
&&\frac{2}{n}\var\left[K_b(Z_i-z)\Delta_iK_h(\bb\trans\X_i-\bb\trans\x) \frac{\partial\left[
	f_{\bb\trans\X}(\bb\trans\z)E\{I(Z\ge Z_i)\mid\bb\trans\x\}\right]/ \partial \bb\trans\x}{f_{\bb\trans\X}^2(\bb\trans\x)E^2\{I(Z\ge Z_i)\mid\bb\trans\x\}}\right]\\
&=&\frac{2}{n}E\left[K_b(Z_i-z)\Delta_iK_h(\bb\trans\X_i-\bb\trans\x) \frac{\partial\left[
	f_{\bb\trans\X}(\bb\trans\z)E\{I(Z\ge Z_i)\mid\bb\trans\x\}\right]/ \partial \bb\trans\x}{f_{\bb\trans\X}^2(\bb\trans\x)E^2\{I(Z\ge Z_i)\mid\bb\trans\x\}}\right]^2\\
&&+\frac{2}{n}\left(E\left[K_b(Z_i-z)\Delta_iK_h(\bb\trans\X_i-\bb\trans\x) \frac{\partial\left[
	f_{\bb\trans\X}(\bb\trans\z)E\{I(Z\ge Z_i)\mid\bb\trans\x\}\right]/ \partial \bb\trans\x}{f_{\bb\trans\X}^2(\bb\trans\x)E^2\{I(Z\ge Z_i)\mid\bb\trans\x\}}\right]\right)^2\\
&=&\frac{2}{n}\iint K_b^2(z_i-z)K_h^2(\bb\trans\x_i-\bb\trans\x) \frac{\left(\partial\left[
	f_{\bb\trans\X}(\bb\trans\x)E\{I(Z\ge z_i)\mid\bb\trans\x\}\right]/ \partial \bb\trans\x\right)^2}{f_{\bb\trans\X}^4(\bb\trans\x)E^4\{I(Z\ge z_i)\mid\bb\trans\x\}}\\
&&\times f(z_i,\bb\trans\x_i)S_c(z_i,\x_i)f_{\X\mid\bb\trans\X}(\x_i,\bb\trans\x_i)f_{\bb\trans\X}(\bb\trans\x_i)dz_id\x_id\bb\trans\x_i+O(1/n)\\
&=&\frac{2}{n}\iint K_b^2(z_i-z)K_h^2(\bb\trans\x_i-\bb\trans\x) \frac{\left(\partial\left[
	f_{\bb\trans\X}(\bb\trans\x)S(z_i,\bb\trans\x)E\{S_c(z_i,X_i)\mid\bb\trans\x\}\right]/ \partial \bb\trans\x\right)^2}{f_{\bb\trans\X}^4(\bb\trans\x)S^4(z_i,\bb\trans\x)E^4\{S_c(z_i,X_i)\mid\bb\trans\x\}}\\
&&\times f(z_i,\bb\trans\x_i)E\{S_c(z_i,X_i)\mid\bb\trans\x_i\}f_{\bb\trans\X}(\bb\trans\x_i)dz_id\x_id\bb\trans\x_i+O(1/n)\\
&=&\frac{2}{nbh}\iint K^2(v)K^2(u) \frac{\left(\partial\left[
	f_{\bb\trans\X}(\bb\trans\x)S(z+bv,\bb\trans\x)E\{S_c(z+bv,X_i)\mid\bb\trans\x\}\right]/ \partial \bb\trans\x\right)^2}{f_{\bb\trans\X}^4(\bb\trans\x)S^4(z+bv,\bb\trans\x)E^4\{S_c(z+bv,X_i)\mid\bb\trans\x\}}\\
&&\times f(z+bv,\bb\trans\x+hu)E\{S_c(z+bv,X_i)\mid\bb\trans\x+hu\}f_{\bb\trans\X}(\bb\trans\x+hu)dvdu+O(1/n)\\
&=&\frac{2}{nbh}\iint K^2(v)K^2(u) \frac{\left(\partial\left[
	f_{\bb\trans\X}(\bb\trans\x)S(z,\bb\trans\x)E\{S_c(z,X_i)\mid\bb\trans\x\}\right]/ \partial \bb\trans\x\right)^2}{f_{\bb\trans\X}^3(\bb\trans\x)S^4(z,\bb\trans\x)E^3\{S_c(z,X_i)\mid\bb\trans\x\}}f(z,\bb\trans\x)dvdu\\
&&+\frac{b}{nh}\frac{\partial^2}{\partial z^2}\iint K^2(v)K^2(u) \frac{\left(\partial\left[
	f_{\bb\trans\X}(\bb\trans\x)S(z^*,\bb\trans\x)E\{S_c(z^*,X_i)\mid\bb\trans\x\}\right]/ \partial \bb\trans\x\right)^2}{f_{\bb\trans\X}^3(\bb\trans\x)S^4(z^*,\bb\trans\x)E^3\{S_c(z^*,X_i)\mid\bb\trans\x\}}f(z^*,\bb\trans\x)v^2dvdu\\
&&+\frac{h}{nb}\frac{\partial^2}{\partial (\bb\trans\x)^2}\iint K^2(v)K^2(u) \frac{\left(\partial\left[
	f_{\bb\trans\X}(\bb\trans\x)S(z,\bb\trans\x)E\{S_c(z,X_i)\mid\bb\trans\x\}\right]/ \partial \bb\trans\x\right)^2}{f_{\bb\trans\X}^4(\bb\trans\x)S^4(z,\bb\trans\x)E^4\{S_c(z,X_i)\mid\bb\trans\x\}}\\
&&\times f(z,\bb\trans\x)E\{S_c(z,X_i)\mid\bb\trans\x^*\}f_{\bb\trans\X}(\bb\trans\x^*)u^2dvdu+O(1/n)\\
&\le&\frac{2}{nbh}f(z,\bb\trans\x)\frac{\left(\partial\left[
	f_{\bb\trans\X}(\bb\trans\x)S(z,\bb\trans\x)E\{S_c(z,X_i)\mid\bb\trans\x\}\right]/ \partial \bb\trans\x\right)^2}{f_{\bb\trans\X}^3(\bb\trans\x)S^4(z,\bb\trans\x)E^3\{S_c(z,X_i)\mid\bb\trans\x\}}\left\{\int K^2(v)dv\right\}\left\{\int K^2(u) du\right\}\\
&&+\frac{b}{nh}\sup_{z^*,\bb\trans\x}\left|\frac{\partial^2}{\partial z^2} f(z^*,\bb\trans\x)\frac{\left(\partial\left[
	f_{\bb\trans\X}(\bb\trans\x)S(z^*,\bb\trans\x)E\{S_c(z^*,X_i)\mid\bb\trans\x\}\right]/ \partial \bb\trans\x\right)^2}{f_{\bb\trans\X}^3(\bb\trans\x)S^4(z^*,\bb\trans\x)E^3\{S_c(z^*,X_i)\mid\bb\trans\x\}}\right|\\
&&\times\left\{\int v^2K^2(v)dv\right\}\left\{\int K^2(u)du\right\}\\
&&+\frac{h}{nb}\sup_{z,\bb\trans\x,\bb\trans\x^*} \left|\frac{\left(\partial\left[
	f_{\bb\trans\X}(\bb\trans\x)S(z,\bb\trans\x)E\{S_c(z,X_i)\mid\bb\trans\x\}\right]/ \partial \bb\trans\x\right)^2}{f_{\bb\trans\X}^4(\bb\trans\x)S^4(z,\bb\trans\x)E^4\{S_c(z,X_i)\mid\bb\trans\x\}}\right.\\
&&\times f(z,\bb\trans\x)E\{S_c(z,X_i)\mid\bb\trans\x^*\}f_{\bb\trans\X}(\bb\trans\x^*)\Bigg|\left\{\int K^2(v)dv\right\}\left\{\int u^2K^2(u)du\right\}+O(1/n)\\
&=&O\{1/(nbh)+b/(nh)+h/(nb)+1/n\}
\ese
uniformly under conditions \ref{assum:kernel}-\ref{assum:survivalfunction}.
Noting that $A=O_p(n^{-1/2}h^{-1/2}+h^2)$ and $B=O_p(n^{-1/2}h^{-3/2}+h^2)$, the second part is
\bse
&&4\var\left[\frac{1}{n}\sumi K_b(Z_i-z)\Delta_iK_h(\bb\trans\X_i-\bb\trans\x) \frac{\partial\left[
	f_{\bb\trans\X}(\bb\trans\z)E\{I(Z\ge Z_i)\mid\bb\trans\x\}\right]/ \partial \bb\trans\x}{f_{\bb\trans\X}^2(\bb\trans\x)E^2\{I(Z\ge Z_i)\mid\bb\trans\x\}}O_p(B)\right]\\
&\le&4E\left(\frac{1}{n}\sumi \left[K_b(Z_i-z)\Delta_iK_h(\bb\trans\X_i-\bb\trans\x) \frac{\partial\left[
	f_{\bb\trans\X}(\bb\trans\z)E\{I(Z\ge Z_i)\mid\bb\trans\x\}\right]/ \partial \bb\trans\x}{f_{\bb\trans\X}^2(\bb\trans\x)E^2\{I(Z\ge Z_i)\mid\bb\trans\x\}}\right]^2\right)\\
&&\qquad\times O\{1/(nh)+h^4\}\\
&=&\left(\frac{4}{nbh}f(z,\bb\trans\x)\frac{\left(\partial\left[
	f_{\bb\trans\X}(\bb\trans\x)S(z,\bb\trans\x)E\{S_c(z,X_i)\mid\bb\trans\x\}\right]/ \partial \bb\trans\x\right)^2}{f_{\bb\trans\X}^3(\bb\trans\x)S^4(z,\bb\trans\x)E^3\{S_c(z,X_i)\mid\bb\trans\x\}}\right.\\
&&\left.\times \left\{\int K^2(v)dv\right\}\left\{\int K^2(u) du\right\}+O\{b/(nh)+h/(nb)\}\right)O(n^{-1}h^{-3}+h^4)\\
&=&O\{n^{-2}b^{-1}h^{-4}+(nb)^{-1}h^3\}
\ese
under conditions \ref{assum:kernel}-\ref{assum:survivalfunction} uniformly. The last part is
\bse
&&4\var\left[\frac{1}{n}\sumi K_b(Z_i-z)\Delta_iK_h(\bb\trans\X_i-\bb\trans\x) \frac{\partial\left[
	f_{\bb\trans\X}(\bb\trans\z)E\{I(Z\ge Z_i)\mid\bb\trans\x\}\right]/ \partial \bb\trans\x}{f_{\bb\trans\X}^2(\bb\trans\x)E^2\{I(Z\ge Z_i)\mid\bb\trans\x\}}O_p(A)\right]\\
&\le&4E\left(\frac{1}{n}\sumi \left[K_b(Z_i-z)\Delta_iK_h(\bb\trans\X_i-\bb\trans\x) \frac{\partial\left[
	f_{\bb\trans\X}(\bb\trans\z)E\{I(Z\ge Z_i)\mid\bb\trans\x\}\right]/ \partial \bb\trans\x}{f_{\bb\trans\X}^2(\bb\trans\x)E^2\{I(Z\ge Z_i)\mid\bb\trans\x\}}\right]^2\right)\\
&&\qquad\times O\{(nh)^{-1}+h^4\}\\
&=&\left(\frac{4}{nbh}f(z,\bb\trans\x)\frac{\left(\partial\left[
	f_{\bb\trans\X}(\bb\trans\x)S(z,\bb\trans\x)E\{S_c(z,X_i)\mid\bb\trans\x\}\right]/ \partial \bb\trans\x\right)^2}{f_{\bb\trans\X}^3(\bb\trans\x)S^4(z,\bb\trans\x)E^3\{S_c(z,X_i)\mid\bb\trans\x\}}\right.\\
&&\left.\times \left\{\int K^2(v)dv\right\}\left\{\int K^2(u) du\right\}+O\{b/(nh)+h/(nb)\}\right)O(n^{-1}h^{-1}+h^4)\\
&=&O\{(nh)^{-2}b^{-1}+h^3(nb)^{-1}\}
\ese
under conditions \ref{assum:kernel}-\ref{assum:survivalfunction} uniformly.
Therefore, 
\bse
\var(\wh\blam_{12})=O\{1/(nbh)\}
\ese
uniformly under conditions \ref{assum:kernel}-\ref{assum:survivalfunction}.

Summarizing the results above,
$\var\{\wh\blam_1(z,\bb\trans\x)\}=O\{1/(nbh^3)\}$ uniformly. Hence
the estimator $\wh\blam_1(z,\bb\trans\x)$ satisfies
\be
\wh\blam_1(x,\bb\trans\x)=\blam_1(x,\bb\trans\x)+O_p\{(nbh^3)^{-1/2}+h^2+b^2\}\n
\ee
uniformly under conditions \ref{assum:kernel}-\ref{assum:survivalfunction}.
\qed

\subsection{Proof of Theorem  \ref{th:consistency}}\label{sec:proofthconsistency}

For each $n$, let $\wh\bb_n$ satisfy
\bse
\frac{1}{n}\sumi \Delta_i\frac{\wh\blam_1(Z_i,\wh\bb_n\trans\X_i)}{\wh\lambda(Z_i,\wh\bb_n\trans\X_i)}
\otimes\left[\X_{li}-
\frac{\wh E\left\{\X_{li} 
	Y_i(Z_i)\mid\wh\bb_n\trans\X_i\right\}}
{\wh E\left\{Y_i(Z_i)\mid\wh\bb_n\trans\X_i\right\}}\right]=\0.
\ese
Under condition \ref{assum:bounded}, there exists a subsequence of
$\wh\bb_n, n=1, 2,\dots$, 
that converges. For notational simplicity,  we still write $\wh\bb_n,
n=1, 2, \dots, $
as the subsequence that  
converges and let the limit be $\bb^*$. 

From the uniform convergence in (\ref{eq:lemeq1}), (\ref{eq:lemeq2}),
(\ref{eq:lemeq5}), (\ref{eq:lemeq6}) given in
Lemma \ref{lem:pre}, we obtain
\bse
&&\frac{1}{n}\sumi \Delta_i\frac{\wh\blam_1(Z_i,\wh\bb_n\trans\X_i)}{\wh\lambda(Z_i,\wh\bb_n\trans\X_i)}
\otimes\left[\X_{li}-
\frac{\wh E\left\{\X_{li} 
	Y_i(Z_i)\mid\wh\bb_n\trans\X_i\right\}}
{\wh E\left\{Y_i(Z_i)\mid\wh\bb_n\trans\X_i\right\}}\right]\\
&=&\frac{1}{n}\sumi \Delta_i\frac{\blam_1(Z_i,\wh\bb_n\trans\X_i)+O_p\{(nbh^3)^{-1/2}+h^2+b^2\}}{\lambda(Z_i,\wh\bb_n\trans\X_i)+O_p\{(nbh)^{-1/2}+h^2+b^2\}}\\
&&
\otimes\left[\X_{li}-
\frac{E\left\{\X_{li} 
	Y_i(Z_i)\mid\wh\bb_n\trans\X_i\right\}+O_p\{(nh)^{-1/2}+h^2\}}
{E\left\{Y_i(Z_i)\mid\wh\bb_n\trans\X_i\right\}+O_p\{(nh)^{-1/2}+h^2\}}\right]\\
&=&\frac{1}{n}\sumi \Delta_i\left[\frac{\blam_1(Z_i,\wh\bb_n\trans\X_i)}{\lambda(Z_i,\wh\bb_n\trans\X_i)}+O_p\{(nbh^3)^{-1/2}+h^2+b^2\}\right]\\
&&
\otimes\left[\X_{li}-
\frac{E\left\{\X_{li} 
	Y_i(Z_i)\mid\wh\bb_n\trans\X_i\right\}}
{E\left\{Y_i(Z_i)\mid\wh\bb_n\trans\X_i\right\}}+O_p\{(nh)^{-1/2}+h^2\}\right]\\
&=&\frac{1}{n}\sumi \Delta_i\frac{\blam_1(Z_i,\wh\bb_n\trans\X_i)}{\lambda(Z_i,\wh\bb_n\trans\X_i)}
\otimes\left[\X_{li}-
\frac{E\left\{\X_{li} 
	Y_i(Z_i)\mid\wh\bb_n\trans\X_i\right\}}
{E\left\{Y_i(Z_i)\mid\wh\bb_n\trans\X_i\right\}}\right]+o_p(1).
\ese
Thus,  for sufficiently large $n$, we have 
\bse
&&\frac{1}{n}\sumi \Delta_i\frac{\blam_1(Z_i,\wh\bb_n\trans\X_i)}{\lambda(Z_i,\wh\bb_n\trans\X_i)}
\otimes\left[\X_{li}-
\frac{E\left\{\X_{li} 
	Y_i(Z_i)\mid\wh\bb_n\trans\X_i\right\}}
{E\left\{Y_i(Z_i)\mid\wh\bb_n\trans\X_i\right\}}\right]\\
&=&\frac{1}{n}\sumi \Delta_i\frac{\blam_1(Z_i,{\bb^*}\trans\X_i)}{\lambda(Z_i,{\bb^*}\trans\X_i)}
\otimes\left[\X_{li}-
\frac{E\left\{\X_{li} 
	Y_i(Z_i)\mid{\bb^*}\trans\X_i\right\}}
{E\left\{Y_i(Z_i)\mid{\bb^*}\trans\X_i\right\}}\right]+O_p(\wh\bb_n-\bb^*)\\
&=&\frac{1}{n}\sumi \Delta_i\frac{\blam_1(Z_i,{\bb^*}\trans\X_i)}{\lambda(Z_i,{\bb^*}\trans\X_i)}
\otimes\left[\X_{li}-
\frac{E\left\{\X_{li} 
	Y_i(Z_i)\mid{\bb^*}\trans\X_i\right\}}
{E\left\{Y_i(Z_i)\mid{\bb^*}\trans\X_i\right\}}\right]+o_p(1),
\ese
where the first equality is because the first derivative of the summation with respect to $\bb$ is
bounded uniformly under conditions \ref{assum:kernel}-\ref{assum:bandwidth} by Lemma \ref{lem:pre}, and
the last equality is because $\wh\bb_n$ converges to $\bb^*$.
In addition,

\bse
&&\frac{1}{n}\sumi \Delta_i\frac{\blam_1(Z_i,{\bb^*}\trans\X_i)}{\lambda(Z_i,{\bb^*}\trans\X_i)}
\otimes\left[\X_{li}-
\frac{E\left\{\X_{li} 
	Y_i(Z_i)\mid{\bb^*}\trans\X_i\right\}}
{E\left\{Y_i(Z_i)\mid{\bb^*}\trans\X_i\right\}}\right]\\
&=&E\left(\Delta\frac{\blam_1(Z,{\bb^*}\trans\X)}{\lambda(Z,{\bb^*}\trans\X)}
\otimes\left[\X_{l}-
\frac{E\left\{\X_{l} 
	Y(Z)\mid{\bb^*}\trans\X\right\}}
{E\left\{Y(Z)\mid{\bb^*}\trans\X_i\right\}}\right]\right)+o_p(1)
\ese
under conditions \ref{assum:kernel}-\ref{assum:bandwidth}. Thus, for sufficient large $n$ we have
\bse
\0&=&\frac{1}{n}\sumi \Delta_i\frac{\wh\blam_1(Z_i,\wh\bb_n\trans\X_i)}{\wh\lambda(Z_i,\wh\bb_n\trans\X_i)}
\otimes\left[\X_{li}-
\frac{\wh E\left\{\X_{li} 
	Y_i(Z_i)\mid\wh\bb_n\trans\X_i\right\}}
{\wh E\left\{Y_i(Z_i)\mid\wh\bb_n\trans\X_i\right\}}\right]\\
&=&\frac{1}{n}\sumi \Delta_i\frac{\blam_1(Z_i,\wh\bb_n\trans\X_i)}{\lambda(Z_i,\wh\bb_n\trans\X_i)}
\otimes\left[\X_{li}-
\frac{E\left\{\X_{li} 
	Y_i(Z_i)\mid\wh\bb_n\trans\X_i\right\}}
{E\left\{Y_i(Z_i)\mid\wh\bb_n\trans\X_i\right\}}\right]+o_p(1)\\
&=&\frac{1}{n}\sumi \Delta_i\frac{\blam_1(Z_i,{\bb^*}\trans\X_i)}{\lambda(Z_i,{\bb^*}\trans\X_i)}
\otimes\left[\X_{li}-
\frac{E\left\{\X_{li} 
	Y_i(Z_i)\mid{\bb^*}\trans\X_i\right\}}
{E\left\{Y_i(Z_i)\mid{\bb^*}\trans\X_i\right\}}\right]+o_p(1)\\
&=&E\left(\Delta\frac{\blam_1(Z,{\bb^*}\trans\X)}{\lambda(Z,{\bb^*}\trans\X)}
\otimes\left[\X_{l}-
\frac{E\left\{\X_{l} 
	Y(Z)\mid{\bb^*}\trans\X\right\}}
{E\left\{Y(Z)\mid{\bb^*}\trans\X_i\right\}}\right]\right)+o_p(1)
\ese
under conditions \ref{assum:kernel}-\ref{assum:bandwidth} and \ref{assum:bounded}. Note that 
\bse
E\left(\Delta\frac{\blam_1(Z,{\bb^*}\trans\X)}{\lambda(Z,{\bb^*}\trans\X)}
\otimes\left[\X_{l}-
\frac{E\left\{\X_{l} 
	Y(Z)\mid{\bb^*}\trans\X\right\}}
{E\left\{Y(Z)\mid{\bb^*}\trans\X_i\right\}}\right]\right)
\ese
is a nonrandom quantity that does not depend on $n$, hence it is zero.
Thus the uniqueness requirement in Condition \ref{assum:unique} ensures that
$\bb^*=\bb_0$.

We now show that the subsequence that converges includes all but a
finite number of $n$'s. Assume this is not the case, then we can
obtain an infinite sequence of $\wh\bb_n$'s that do not converge to
$\bb^*$. As an infinite sequence in a compact set $\cal B$, we can
thus obtain another subsequence that converges, say to
$\bb^{**}\ne\bb^*$. Identical derivation as before then leads to
$\bb^{**}=\bb_0$, which is a contradiction to 
$\bb^{**}\ne\bb^*$.
Thus we conclude
$
\wh\bb-\bb_0\to\0
$
in probability when $n\to\infty$ under condition
\ref{assum:kernel}-\ref{assum:unique}.
\qed

\subsection{Proof of Theorem \ref{th:eff}}\label{sec:prooftheff}
We first expand (\ref{eq:eff}) as
\be
\0
&=&n^{-1/2}\sumi \Delta_i\frac{\wh\blam_1(Z_i,\wh\bb\trans\X_i)}{\wh\lambda(Z_i,\wh\bb\trans\X_i)}
\otimes\left[\X_{li}-
\frac{\wh E\left\{\X_{li} 
	Y_i(Z_i)\mid\wh\bb\trans\X_i\right\}}
{\wh E\left\{Y_i(Z_i)\mid\wh\bb\trans\X_i\right\}}\right]\nonumber\\
&=&n^{-1/2}\sumi \Delta_i\frac{\wh\blam_1(Z_i,\bb_0\trans\X_i)}{\wh\lambda(Z_i,\bb_0\trans\X_i)}
\otimes\left[\X_{li}-
\frac{\wh E\left\{\X_{li} 
	Y_i(Z_i)\mid\bb_0\trans\X_i\right\}}
{\wh E\left\{Y_i(Z_i)\mid\bb_0\trans\X_i\right\}}\right]\label{eq:main}\\
&&+\frac{1}{n}\sumi \left\{\frac{\partial}{\partial(\X_i\trans\bb)}
\left(\Delta_i\frac{\wh\blam_1(Z_i,\bb\trans\X_i)}{\wh\lambda(Z_i,\bb\trans\X_i)}
\otimes\left[\X_{li}-
\frac{\wh E\left\{\X_{li} 
	Y_i(Z_i)\mid\bb\trans\X_i\right\}}
{\wh
	E\left\{Y_i(Z_i)\mid\bb\trans\X_i\right\}}\right]\right)\otimes\X_{li}\trans\right\}\Bigg|_{\bb=\wt\bb}\label{eq:easy}\notag\\\\
&&\times\sqrt{n}(\wh\bb-\bb_0),\nonumber
\ee
where $\wt\bb$ is on the line connecting $\bb_0$ and $\wh\bb$.

We first consider (\ref{eq:easy}). Because of Theorem
\ref{th:consistency} and Lemma \ref{lem:pre}, 
we have
\be
&&\frac{1}{n}\sumi \left\{\frac{\partial}{\partial(\X_i\trans\bb)}
\left(\Delta_i\frac{\wh\blam_1(Z_i,\bb\trans\X_i)}{\wh\lambda(Z_i,\bb\trans\X_i)}
\otimes\left[\X_{li}-
\frac{\wh E\left\{\X_{li} 
	Y_i(Z_i)\mid\bb\trans\X_i\right\}}
{\wh
	E\left\{Y_i(Z_i)\mid\bb\trans\X_i\right\}}\right]\right)\otimes\X_{li}\trans\right\}\Bigg|_{\bb=\wt\bb}\n\\
&=&\frac{1}{n}\sumi \left\{\frac{\partial}{\partial(\X_i\trans\bb_0)}
\left(\Delta_i\frac{\wh\blam_1(Z_i,\bb_0\trans\X_i)}{\wh\lambda(Z_i,\bb_0\trans\X_i)}
\otimes\left[\X_{li}-
\frac{\wh E\left\{\X_{li} 
	Y_i(Z_i)\mid\bb_0\trans\X_i\right\}}
{\wh
	E\left\{Y_i(Z_i)\mid\bb_0\trans\X_i\right\}}\right]\right)\otimes\X_{li}\trans\right\}\n\\
&&+o_p(1)\n\\
&=&-\frac{1}{n}\sumi \left(\Delta_i
\frac{\wh\blam_1^{\otimes2}(Z_i,\bb_0\trans\X_i)}{\wh\lambda^2(Z_i,\bb_0\trans\X_i)}
\otimes\left[\X_{li}-
\frac{\wh E\left\{\X_{li} 
	Y_i(Z_i)\mid\bb_0\trans\X_i\right\}}
{\wh
	E\left\{Y_i(Z_i)\mid\bb_0\trans\X_i\right\}}\right]\otimes\X_{li}\trans\right)\label{eq:easy1}\\
&&+\frac{1}{n}\sumi \frac{\Delta_i}{\wh\lambda(Z_i,\bb_0\trans\X_i)}
\frac{\partial}{\partial(\X_i\trans\bb_0)}
\left(\wh\blam_1(Z_i,\bb_0\trans\X_i)
\otimes\left[\X_{li}-
\frac{\wh E\left\{\X_{li} 
	Y_i(Z_i)\mid\bb_0\trans\X_i\right\}}
{\wh
	E\left\{Y_i(Z_i)\mid\bb_0\trans\X_i\right\}}\right]\right)\otimes\X_{li}\trans\n\notag\\\\
&&+o_p(1).\nonumber\label{eq:easy2}
\ee
Because of Lemma \ref{lem:pre}, (\ref{eq:easy1}) converges uniformly
in probability to
\bse
&&-E \left(\int_0^\infty 
\frac{\blam_1^{\otimes2}(s,\bb_0\trans\X)}{\lambda^2(s,\bb_0\trans\X)}
\otimes\left[\X_l-
\frac{ E\left\{\X_l 
	Y(s)\mid\bb_0\trans\X\right\}}
{E\left\{Y(s)\mid\bb_0\trans\X\right\}}\right]\otimes\X_l\trans
dN(s)\right)\\
&=&-E \left(\int_0^\infty 
\frac{\blam_1^{\otimes2}(s,\bb_0\trans\X)}{\lambda^2(s,\bb_0\trans\X)}
\otimes\left[\X_l-
\frac{ E\left\{\X_l 
	Y(s)\mid\bb_0\trans\X\right\}}
{E\left\{Y(s)\mid\bb_0\trans\X\right\}}\right]\otimes\X_l\trans 
Y(s)\lambda(s,\bb_0\trans\X)ds\right)\\
&=&-E \left(\int_0^\infty 
\frac{\blam_1^{\otimes2}(s,\bb_0\trans\X)}{\lambda(s,\bb_0\trans\X)}
\otimes\left[\X_l-
\frac{ E\left\{\X_l 
	Y(s)\mid\bb_0\trans\X\right\}}
{E\left\{Y(s)\mid\bb_0\trans\X\right\}}\right]\otimes
\left[\X_l
-\frac{ E\left\{\X_l 
	Y(s)\mid\bb_0\trans\X\right\}}
{E\left\{Y(s)\mid\bb_0\trans\X\right\}}\right]
\trans 
Y(s)ds\right)\\
&&-E \left(\int_0^\infty 
\frac{\blam_1^{\otimes2}(s,\bb_0\trans\X)}{\lambda(s,\bb_0\trans\X)}
\otimes\left[\X_l-
\frac{ E\left\{\X_l 
	Y(s)\mid\bb_0\trans\X\right\}}
{E\left\{Y(s)\mid\bb_0\trans\X\right\}}\right]\otimes
\frac{ E\left\{\X_l 
	Y(s)\mid\bb_0\trans\X\right\}}
{E\left\{Y(s)\mid\bb_0\trans\X\right\}}\trans 
Y(s)ds\right)\\
&=&-E\{\bS\eff(\Delta,Z,\X)^{\otimes2}\},
\ese
where the last equality is because the second term above is zero by 
first taking expectation conditional on $\bb_0\trans\X$.

Similarly, from Lemma \ref{lem:pre},  the term in (\ref{eq:easy2})
converges uniformly in probability to the
limit of
\bse
E \left\{\frac{\Delta_i}{\lambda(Z_i,\bb_0\trans\X_i)}
\frac{\partial}{\partial(\X_i\trans\bb_0)}
\left(\wh\blam_1(Z_i,\bb_0\trans\X_i)
\otimes \left[\X_{li}-
\frac{ E\left\{\X_{li} Y_i(Z_i)\mid\bb_0\trans\X_i\right\}}
{E\left\{Y_i(Z_i)\mid\bb_0\trans\X_i\right\}}
\right]
\right)\otimes\X_{li}\trans\right\}.
\ese
Now let $\wh\blam_{1,-i}(Z, \bb_0\trans\X)$ be the
leave-one-out version of $\wh\blam_{1}(Z, \bb_0\trans\X)$, i.e. it
is constructed the same as $\wh\blam_{1}(Z, \bb_0\trans\X)$ except
that the $i$th observation is not used.
Obviously, 
\bse
&&\frac{\Delta_i}{\lambda(Z_i,\bb_0\trans\X_i)}
\frac{\partial}{\partial(\X_i\trans\bb_0)}
\left(\wh\blam_1(Z_i,\bb_0\trans\X_i)
\otimes \left[\X_{li}-
\frac{ E\left\{\X_{li} 
	Y_i(Z_i)\mid\bb_0\trans\X_i\right\}}
{E\left\{Y_i(Z_i)\mid\bb_0\trans\X_i\right\}}\right]\right)\otimes\X_{li}\trans\\
&-&\frac{\Delta_i}{\lambda(Z_i,\bb_0\trans\X_i)}
\frac{\partial}{\partial(\X_i\trans\bb_0)}
\left(\wh\blam_{1,-i}(Z_i,\bb_0\trans\X_i)
\otimes \left[\X_{li}-
\frac{ E\left\{\X_{li} 
	Y_i(Z_i)\mid\bb_0\trans\X_i\right\}}
{E\left\{Y_i(Z_i)\mid\bb_0\trans\X_i\right\}}\right]\right)\otimes\X_{li}\trans\\
&=&o_p(1).
\ese
Now let $E_i$ mean taking expectation with respect to the $i$th
observation conditional on all other observations, then
\bse
&&
E_i\left\{\frac{\Delta_i}{\lambda(Z_i,\bb_0\trans\X_i)}
\frac{\partial}{\partial(\X_i\trans\bb_0)}
\left(\wh\blam_{1,-i}(Z_i,\bb_0\trans\X_i)
\otimes \left[\X_{li}-
\frac{ E\left\{\X_{li} 
	Y_i(Z_i)\mid\bb_0\trans\X_i\right\}}
{E\left\{Y_i(Z_i)\mid\bb_0\trans\X_i\right\}}\right]\right)\otimes\X_{li}\trans\right\}\\
&=&E_i \left\{\int \frac{1}{\lambda(s,\bb_0\trans\X_i)}
\frac{\partial}{\partial(\X_i\trans\bb_0)}
\left(\wh\blam_{1,-i}(s,\bb_0\trans\X_i)
\otimes \left[\X_{li}-
\frac{ E\left\{\X_{li} 
	Y_i(s)\mid\bb_0\trans\X_i\right\}}
{E\left\{Y_i(s)\mid\bb_0\trans\X_i\right\}}\right]\right)\otimes\X_{li}\trans
dN_i(s)\right\}\\
&=&E_i \left\{\int 
\frac{\partial}{\partial(\X_i\trans\bb_0)}
\left(\wh\blam_{1,-i}(s,\bb_0\trans\X_i)
\otimes \left[\X_{li}-
\frac{ E\left\{\X_{li} 
	Y_i(s)\mid\bb_0\trans\X_i\right\}}
{E\left\{Y_i(s)\mid\bb_0\trans\X_i\right\}}\right]\right)\otimes\X_{li}\trans
Y_i(s)ds\right\}\\
&=&E_i \left\{
\frac{\partial}{\partial\bb_0}\int 
\wh\blam_{1,-i}(s,\bb_0\trans\X_i) 
\otimes \left[\X_{li}-
\frac{ E\left\{\X_{li} 
	Y_i(s)\mid\bb_0\trans\X_i\right\}}
{E\left\{Y_i(s)\mid\bb_0\trans\X_i\right\}}\right]
E\{Y_i(s)\mid\X_i\}ds\right\}\\
&=&\frac{\partial}{\partial\bb_0} E_i \left\{
\int 
\wh\blam_{1,-i}(s,\bb_0\trans\X_i) 
\otimes \left[\X_{lI}-
\frac{ E\left\{\X_{li} 
	Y_i(s)\mid\bb_0\trans\X_i\right\}}
{E\left\{Y_i(s)\mid\bb_0\trans\X_i\right\}}\right]
E\{Y_i(s)\mid\X\}ds\right\}\\
&=&\frac{\partial}{\partial\bb_0} E_i \left\{
\int 
\wh\blam_{1,-i}(s,\bb_0\trans\X_i) 
\otimes \left[\X_{li}-
\frac{ E\left\{\X_{li} 
	Y_i(s)\mid\bb_0\trans\X_i\right\}}
{E\left\{Y_i(s)\mid\bb_0\trans\X_i\right\}}\right]
Y_i(s)ds\right\}\\
&=&\0.
\ese
Here, the last equality is because the integrand has expectation zero
conditional on $\bb_0\trans\X_i$ and all other observations, and the third last equality is because
the expectation is with respect to $\X_i$ and does not involve $\bb_0$. 
Therefore, the term in (\ref{eq:easy2}) converges in probability
uniformly to
\bse
E\left\{\frac{\Delta_i}{\lambda(Z_i,\bb_0\trans\X_i)}
\frac{\partial}{\partial(\X_i\trans\bb_0)}
\left(\wh\blam_{1,-i}(Z_i,\bb_0\trans\X_i)
\otimes \left[\X_{li}-
\frac{ E\left\{\X_{li} 
	Y_i(Z_i)\mid\bb_0\trans\X_i\right\}}
{E\left\{Y_i(Z_i)\mid\bb_0\trans\X_i\right\}}\right]\right)\otimes\X_{li}\trans\right\}
=0
\ese
Combining the results concerning (\ref{eq:easy1}) and
(\ref{eq:easy2}), we thus have obtained that the expression in (\ref{eq:easy}) is 
$-E\{\bS\eff(\Delta,Z,\X)^{\otimes2}\}+o_p(1)$.

Next we decompose (\ref{eq:main}) into 
\be\label{eq:Ts}
n^{-1/2}\sumi \Delta_i\frac{\wh\blam_1(Z_i,\bb_0\trans\X_i)}{\wh\lambda(Z_i,\bb_0\trans\X_i)}
\otimes\left[\X_{li}-
\frac{\wh E\left\{\X_{li} 
	Y_i(Z_i)\mid\bb_0\trans\X_i\right\}}
{\wh E\left\{Y_i(Z_i)\mid\bb_0\trans\X_i\right\}}\right]
=\T_1+\T_2+\T_3+\T_4,
\ee
where
\bse
\T_1
&=&n^{-1/2}\sumi \Delta_i\frac{\blam_1(Z_i,\bb_0\trans\X_i)}{\lambda(Z_i,\bb_0\trans\X_i)}
\otimes\left[\X_{li}-
\frac{E\left\{\X_{li} 
	Y_i(Z_i)\mid\bb_0\trans\X_i\right\}}
{E\left\{Y_i(Z_i)\mid\bb_0\trans\X_i\right\}}\right],\\
\T_2&=&n^{-1/2}\sumi \Delta_i\left\{\frac{\wh\blam_1(Z_i,\bb_0\trans\X_i)}{\wh\lambda(Z_i,\bb_0\trans\X_i)}
-\frac{\blam_1(Z_i,\bb_0\trans\X_i)}{\lambda(Z_i,\bb_0\trans\X_i)}\right\}\otimes\left[\X_{li}-
\frac{E\left\{\X_{li} 
	Y_i(Z_i)\mid\bb_0\trans\X_i\right\}}
{E\left\{Y_i(Z_i)\mid\bb_0\trans\X_i\right\}}\right],\\
\T_3&=&n^{-1/2}\sumi \Delta_i\frac{\blam_1(Z_i,\bb_0\trans\X_i)}{\lambda(Z_i,\bb_0\trans\X_i)}
\otimes\left[\frac{E\left\{\X_{li} 
	Y_i(Z_i)\mid\bb_0\trans\X_i\right\}}
{E\left\{Y_i(Z_i)\mid\bb_0\trans\X_i\right\}}-
\frac{\wh E\left\{\X_{li} 
	Y_i(Z_i)\mid\bb_0\trans\X_i\right\}}
{\wh
	E\left\{Y_i(Z_i)\mid\bb_0\trans\X_i\right\}}\right],\\
\T_4&=&n^{-1/2}\sumi \Delta_i
\left\{\frac{\wh\blam_1(Z_i,\bb_0\trans\X_i)}{\wh\lambda(Z_i,\bb_0\trans\X_i)}-
\frac{\blam_1(Z_i,\bb_0\trans\X_i)}{\lambda(Z_i,\bb_0\trans\X_i)}\right\}\\
&&\otimes\left[\frac{E\left\{\X_{li} 
	Y_i(Z_i)\mid\bb_0\trans\X_i\right\}}
{E\left\{Y_i(Z_i)\mid\bb_0\trans\X_i\right\}}-
\frac{\wh E\left\{\X_{li} 
	Y_i(Z_i)\mid\bb_0\trans\X_i\right\}}
{\wh
	E\left\{Y_i(Z_i)\mid\bb_0\trans\X_i\right\}}\right].
\ese

First, note  that 
\bse
\T_2&=&n^{-1/2}\sumi \int\left\{\frac{\wh\blam_1(s,\bb_0\trans\X_i)}{\wh\lambda(s,\bb_0\trans\X_i)}
-\frac{\blam_1(s,\bb_0\trans\X_i)}{\lambda(s,\bb_0\trans\X_i)}\right\}\otimes\left[\X_{li}-
\frac{E\left\{\X_{li} 
	Y_i(s)\mid\bb_0\trans\X_i\right\}}
{E\left\{Y_i(s)\mid\bb_0\trans\X_i\right\}}\right]dN_i(s)\\
&=&o_p\left(n^{-1/2}\sumi \int\left[\X_{li}-
\frac{E\left\{\X_{li} 
	Y_i(s)\mid\bb_0\trans\X_i\right\}}
{E\left\{Y_i(s)\mid\bb_0\trans\X_i\right\}}\right]Y_i(s)\lambda(s,\bb_0\trans\X_{li})ds\right)\\
&=&o_p(1),
\ese
where the last equality above is because the quantity inside the parenthesis is a mean zero normal random
quantity of order $O_p(1)$. Further,
\bse
\T_3&=&n^{-1/2}\sumi \Delta_i\frac{\blam_1(Z_i,\bb_0\trans\X_i)}{\lambda(Z_i,\bb_0\trans\X_i)}
\otimes\left(-
\frac{\wh E\left\{\X_{li} Y_i(Z_i)\mid\bb_0\trans\X_i\right\}}{E\left\{Y_i(Z_i)\mid\bb_0\trans\X_i\right\}}\right.\n\\
&&\left. +\frac{\wh E\left\{Y_i(Z_i)\mid\bb_0\trans\X_i\right\}E\left\{\X_{li}Y_i(Z_i)\mid\bb_0\trans\X_i\right\}}
{[E\left\{Y_i(Z_i)\mid\bb_0\trans\X_i\right\}]^2}\right) +o_p(1)\n\\
&=&
n^{-1/2}\sumi \Delta_i\frac{\blam_1(Z_i,\bb_0\trans\X_i)}{\lambda(Z_i,\bb_0\trans\X_i)}
\otimes\left(-
\frac{
	n^{-1}\sumj K_h(\bb_0\trans\X_j-\bb_0\trans\X_i)\X_{lj}I(Z_j\ge
	Z_i)}{f_{\bb_0\trans\X}(\bb_0\trans\X_i) E\left\{Y_i(Z_i)\mid\bb_0\trans\X_i\right\}}\right.\n\\
&&\left.+\frac{E\left\{\X_{li}Y_i(Z_i)\mid\bb_0\trans\X_i\right\}
	\{n^{-1}\sumj K_h(\bb_0\trans\X_j-\bb_0\trans\X_i)-f_{\bb_0\trans\X}(\bb_0\trans\X_i)
	\}}{f_{\bb_0\trans\X}(\bb_0\trans\X_i) E\left\{Y_i(Z_i)\mid\bb_0\trans\X_i\right\}}\right.\n\\
&&\left. +\frac{E\left\{\X_{li}Y_i(Z_i)\mid\bb_0\trans\X_i\right\}}
{[E\left\{Y_i(Z_i)\mid\bb_0\trans\X_i\right\}]^2}
\left[\frac{
	n^{-1}\sumj K_h(\bb_0\trans\X_j-\bb_0\trans\X_i)I(Z_j\ge
	Z_i)}{f_{\bb_0\trans\X}(\bb_0\trans\X_i) }\right.\right.\n\\
&&\left.\left.-\frac{E\left\{Y_i(Z_i)\mid\bb_0\trans\X_i\right\}
	\{n^{-1}\sumj K_h(\bb_0\trans\X_j-\bb_0\trans\X_i)-f_{\bb_0\trans\X}(\bb_0\trans\X_i)
	\}}{f_{\bb_0\trans\X}(\bb_0\trans\X_i) }\right]
\right) +o_p(1)\n\\
&=&n^{-3/2}\sumi\sumj \Delta_i\frac{\blam_1(Z_i,\bb_0\trans\X_i)}{\lambda(Z_i,\bb_0\trans\X_i)}
\otimes\left[
-\frac{K_h(\bb_0\trans\X_j-\bb_0\trans\X_i)\X_{lj}I(Z_j\ge 
	Z_i)}{f_{\bb_0\trans\X}(\bb_0\trans\X_i) E\left\{Y_i(Z_i)\mid\bb_0\trans\X_i\right\}}\right.\n\\
&&\left. +
\frac{E\left\{\X_{li}Y_i(Z_i)\mid\bb_0\trans\X_i\right\} K_h(\bb_0\trans\X_j-\bb_0\trans\X_i)I(Z_j\ge 
	Z_i)}{f_{\bb_0\trans\X}(\bb_0\trans\X_i) [E\left\{Y_i(Z_i)\mid\bb_0\trans\X_i\right\}]^2
}\right] +o_p(1)\n\\
&=&\T_{31}+\T_{32}+\T_{33}+o_p(1),
\ese
where
\bse
\T_{31}&=&n^{-1/2}\sumi\Delta_i
\frac{\blam_1(Z_i,\bb_0\trans\X_i)}{\lambda(Z_i,\bb_0\trans\X_i)}
\otimes E\left[
-\frac{K_h(\bb_0\trans\X_j-\bb_0\trans\X_i)\X_{lj}I(Z_j\ge 
	Z_i)}{f_{\bb_0\trans\X}(\bb_0\trans\X_i) E\left\{Y_i(Z_i)\mid\bb_0\trans\X_i\right\}}\right.\n\\
&&\left. +
\frac{E\left\{\X_{li}Y_i(Z_i)\mid\bb_0\trans\X_i\right\} K_h(\bb_0\trans\X_j-\bb_0\trans\X_i)I(Z_j\ge 
	Z_i)}{f_{\bb_0\trans\X}(\bb_0\trans\X_i) [E\left\{Y_i(Z_i)\mid\bb_0\trans\X_i\right\}]^2
}\mid \Delta_i, Z_i, \X_i\right]\\
\T_{32}&=&n^{-1/2}\sumj E \left(\Delta_i\frac{\blam_1(Z_i,\bb_0\trans\X_i)}{\lambda(Z_i,\bb_0\trans\X_i)}
\otimes\left[
-\frac{K_h(\bb_0\trans\X_j-\bb_0\trans\X_i)\X_{lj}I(Z_j\ge 
	Z_i)}{f_{\bb_0\trans\X}(\bb_0\trans\X_i) E\left\{Y_i(Z_i)\mid\bb_0\trans\X_i\right\}}\right.\right.\n\\
&&\left. \left.+
\frac{E\left\{\X_{li}Y_i(Z_i)\mid\bb_0\trans\X_i\right\} K_h(\bb_0\trans\X_j-\bb_0\trans\X_i)I(Z_j\ge 
	Z_i)}{f_{\bb_0\trans\X}(\bb_0\trans\X_i) [E\left\{Y_i(Z_i)\mid\bb_0\trans\X_i\right\}]^2
}\right]\mid \Delta_j, Z_j,\X_j\right)\\
\T_{33}&=&-n^{1/2}E\left(\Delta_i
\frac{\blam_1(Z_i,\bb_0\trans\X_i)}{\lambda(Z_i,\bb_0\trans\X_i)}
\otimes E\left[
-\frac{K_h(\bb_0\trans\X_j-\bb_0\trans\X_i)\X_{lj}I(Z_j\ge 
	Z_i)}{f_{\bb_0\trans\X}(\bb_0\trans\X_i) E\left\{Y_i(Z_i)\mid\bb_0\trans\X_i\right\}}\right.\right.\n\\
&&\left.\left. +
\frac{E\left\{\X_{li}Y_i(Z_i)\mid\bb_0\trans\X_i\right\} K_h(\bb_0\trans\X_j-\bb_0\trans\X_i)I(Z_j\ge 
	Z_i)}{f_{\bb_0\trans\X}(\bb_0\trans\X_i) [E\left\{Y_i(Z_i)\mid\bb_0\trans\X_i\right\}]^2
}\right]\right).
\ese
Here we used U-statistic property in the last equality above. Now when $nh^4\to0$,
\bse
\T_{31}&=&n^{-1/2}\sumi\Delta_i
\frac{\blam_1(Z_i,\bb_0\trans\X_i)}{\lambda(Z_i,\bb_0\trans\X_i)}
\otimes \left[
-\frac{E\{\X_{li}Y_i(Z_i)\mid \bb_0\trans\X_i\}}{E\left\{Y_i(Z_i)\mid\bb_0\trans\X_i\right\}}\right.\\
&&\left.+
\frac{E\left\{\X_{li}Y_i(Z_i)\mid\bb_0\trans\X_i\right\}E\left\{Y_i(Z_i)\mid\bb_0\trans\X_i\right\}}{
	[E\left\{Y_i(Z_i)\mid\bb_0\trans\X_i\right\}]^2
}\right]+O(n^{1/2}h^2)\\
&=&o_p(1).
\ese
Thus, $\T_{33}=o_p(1)$ as well.
To analyze $\T_{32}$, 
\bse
\T_{32}&=&n^{-1/2}\sumj E \left(\Delta_i\frac{\blam_1(Z_i,\bb_0\trans\X_i)}{\lambda(Z_i,\bb_0\trans\X_i)}
\otimes\left[
-\frac{K_h(\bb_0\trans\X_j-\bb_0\trans\X_i)\X_{lj}I(Z_j\ge 
	Z_i)}{f_{\bb_0\trans\X}(\bb_0\trans\X_i) E\left\{I(Z\ge
	Z_i)\mid\bb_0\trans\X=\bb_0\trans\X_i, Z_i\right\}}\right.\right.\n\\
&&\left. \left.+
\frac{E\left\{\X_{l}I(Z\ge Z_i)\mid\bb_0\trans\X=\bb_0\trans\X_i, Z_i\right\} K_h(\bb_0\trans\X_j-\bb_0\trans\X_i)I(Z_j\ge 
	Z_i)}{f_{\bb_0\trans\X}(\bb_0\trans\X_i) [E\left\{I(Z\ge
	Z_i)\mid\bb_0\trans\X=\bb_0\trans\X_i, Z_i\right\}]^2
}\right]\mid \Delta_j, Z_j,\X_j\right)\\
&=&n^{-1/2}\sumj E \left(\Delta_i\frac{\blam_1(Z_i,\bb_0\trans\X_i)}{\lambda(Z_i,\bb_0\trans\X_i)}
\otimes\left[
-\frac{K_h(\bb_0\trans\x_j-\bb_0\trans\X_i)\x_{lj}I(z_j\ge 
	Z_i)}{f_{\bb_0\trans\X}(\bb_0\trans\X_i) E\left\{I(Z\ge
	Z_i)\mid\bb_0\trans\X=\bb_0\trans\X_i, Z_i\right\}}\right.\right.\n\\
&&\left. \left.+
\frac{E\left\{\X_{l}I(Z\ge Z_i)\mid\bb_0\trans\X=\bb_0\trans\X_i, Z_i\right\} K_h(\bb_0\trans\x_j-\bb_0\trans\X_i)I(z_j\ge 
	Z_i)}{f_{\bb_0\trans\X}(\bb_0\trans\X_i) [E\left\{I(Z\ge
	Z_i)\mid\bb_0\trans\X=\bb_0\trans\X_i, Z_i\right\}]^2
}\right]\right)\\
&=&n^{-1/2}\sumj E\left\{E \left(\Delta_i\frac{\blam_1(Z_i,\bb_0\trans\X_i)}{\lambda(Z_i,\bb_0\trans\X_i)}
\otimes\left[
-\frac{\x_{lj}I(z_j\ge 
	Z_i)}{f_{\bb_0\trans\X}(\bb_0\trans\X_i) E\left\{I(Z\ge
	Z_i)\mid\bb_0\trans\X=\bb_0\trans\X_i, Z_i\right\}}\right.\right.\right.\n\\
&&\left. \left.\left.+
\frac{E\left\{\X_{l}I(Z\ge Z_i)\mid\bb_0\trans\X=\bb_0\trans\X_i, Z_i\right\} I(z_j\ge 
	Z_i)}{f_{\bb_0\trans\X}(\bb_0\trans\X_i) [E\left\{I(Z\ge
	Z_i)\mid\bb_0\trans\X=\bb_0\trans\X_i, Z_i\right\}]^2
}\right]\mid \bb_0\trans\X_i\right)
K_h(\bb_0\trans\x_j-\bb_0\trans\X_i)\right\}\\
&=&n^{-1/2}\sumj E \left(\Delta_i\frac{\blam_1(Z_i,\bb_0\trans\x_j)}{\lambda(Z_i,\bb_0\trans\x_j)}
\otimes\left[
-\frac{\x_{lj}I(z_j\ge 
	Z_i)}{E\left\{I(Z\ge
	Z_i)\mid\bb_0\trans\X=\bb_0\trans\x_j, Z_i\right\}}\right.\right.\n\\
&&\left.\left.+
\frac{E\left\{\X_{l}I(Z\ge Z_i)\mid\bb_0\trans\X=\bb_0\trans\x_j, Z_i\right\} I(z_j\ge 
	Z_i)}{[E\left\{I(Z\ge
	Z_i)\mid\bb_0\trans\X=\bb_0\trans\x_j, Z_i\right\}]^2
}\right]\mid\bb_0\trans\X_i=\bb_0\trans\x_j\right)+O_p(n^{1/2}h^2)\\
&=&n^{-1/2}\sumj E \left(\frac{\Delta_i I(z_j\ge Z_i)\blam_1(Z_i,\bb_0\trans\x_j)}{S(Z_i,\bb_0\trans\x_j)\lambda(Z_i,\bb_0\trans\x_j) E\left\{S_c(Z_i,\X)\mid\bb_0\trans\X=\bb_0\trans\x_j, Z_i\right\}}\right.\\
&&\left.\otimes\left[
\frac{E\left\{\X_lS_c(Z_i,\X)\mid\bb_0\trans\X=\bb_0\trans\x_j, Z_i\right\}
}{E\left\{S_c(Z_i,\X)\mid\bb_0\trans\X=\bb_0\trans\x_j, Z_i\right\}
}-\x_{lj}\right]\mid\bb_0\trans\X_i=\bb_0\trans\x_j\right) +O_p(n^{1/2}h^2)\\
&=&n^{-1/2}\sumj E\left(\int_0^{z_j}\frac{\blam_1(s,\bb_0\trans\x_j)}{E\left\{S_c(s,\X)\mid\bb_0\trans\X=\bb_0\trans\x_j\right\}}\right.\\
&&\left.\otimes\left[
\frac{E \left\{\X_lS_c(s,\X)\mid\bb_0\trans\X=\bb_0\trans\x_j\right\}
}{E\left\{S_c(s,\X)\mid\bb_0\trans\X=\bb_0\trans\x_j\right\}
}-\x_{lj}\right]
S_c(s,\X_i)ds
\mid\bb_0\trans\X_i=\bb_0\trans\x_j\right) +O_p(n^{1/2}h^2)\\
&=&n^{-1/2}\sumj \int_0^{z_j}\blam_1(s,\bb_0\trans\x_j)\otimes\left[
\frac{E\left\{\X_{lj}S_c(s,\X_j)\mid\bb_0\trans\x_j\right\}
}{E\left\{S_c(s,\X_j)\mid\bb_0\trans\x_j\right\}
}-\x_{lj}\right]ds +O_p(n^{1/2}h^2)\\
&=&n^{-1/2}\sumj \int Y_j(s)\lambda(s,\bb_0\trans\x_j)\frac{\blam_1(s,\bb_0\trans\x_j)}{\lambda(s,\bb_0\trans\x_j)}\otimes\left[
\frac{E\left\{\X_{lj}Y_j(s)\mid\bb_0\trans\x_j\right\}
}{E\left\{Y_j(s)\mid\bb_0\trans\x_j\right\}
}-\x_{lj}\right]ds +O_p(n^{1/2}h^2).
\ese
When $nh^4\to0$, plugging the results of $\T_1$ and $\T_{32}$ to
(\ref{eq:Ts}), we obtain that the expression in (\ref{eq:main}) is 
\bse
&&n^{-1/2}\sumi \Delta_i\frac{\wh\blam_1(Z_i,\bb_0\trans\X_i)}{\wh\lambda(Z_i,\bb_0\trans\X_i)}
\otimes\left[\X_{li}-
\frac{\wh E\left\{\X_{li} 
	Y_i(Z_i)\mid\bb_0\trans\X_i\right\}}
{\wh E\left\{Y_i(Z_i)\mid\bb_0\trans\X_i\right\}}\right]\\
&=&n^{-1/2}\sumi \int\frac{\blam_1(t,\bb_0\trans\X_i)}{\lambda(t,\bb_0\trans\X_i)}
\otimes\left[\X_{li}-
\frac{E\left\{\X_{li} 
	Y_i(t)\mid\bb_0\trans\X_i\right\}}
{E\left\{Y_i(t)\mid\bb_0\trans\X_i\right\}}\right]dM_i(t)+o_p(1)\\
&=&n^{-1/2}\sumi\bS\eff(\Delta_i,Z_i,\X_i)+o_p(1).
\ese

Finally,  note that 
\bse
\T_4&=&n^{-1/2}\sumi \Delta_i
\left\{\frac{\wh\blam_1(Z_i,\bb_0\trans\X_i)}{\wh\lambda(Z_i,\bb_0\trans\X_i)}-
\frac{\blam_1(Z_i,\bb_0\trans\X_i)}{\lambda(Z_i,\bb_0\trans\X_i)}\right\}\\
&&\times\left[\frac{E\left\{\X_{li} 
	Y_i(Z_i)\mid\bb_0\trans\X_i\right\}}
{E\left\{Y_i(Z_i)\mid\bb_0\trans\X_i\right\}}-
\frac{\wh E\left\{\X_{li} 
	Y_i(Z_i)\mid\bb_0\trans\X_i\right\}}
{\wh
	E\left\{Y_i(Z_i)\mid\bb_0\trans\X_i\right\}}\right]\\
&=&o_p\left(
n^{-1/2}\sumi \Delta_i
\left[\frac{E\left\{\X_{li} 
	Y_i(Z_i)\mid\bb_0\trans\X_i\right\}}
{E\left\{Y_i(Z_i)\mid\bb_0\trans\X_i\right\}}-
\frac{\wh E\left\{\X_{li} 
	Y_i(Z_i)\mid\bb_0\trans\X_i\right\}}
{\wh
	E\left\{Y_i(Z_i)\mid\bb_0\trans\X_i\right\}}\right]\right)\\
&=&o_p\left(n^{-1/2}\sumi \int Y_i(s)\lambda(s,\bb_0\trans\x_i)
\left[
\frac{E\left\{\X_{li}Y_i(s)\mid\bb_0\trans\x_i\right\}
}{E\left\{Y_i(s)\mid\bb_0\trans\x_i\right\}
}-\x_{li}\right]ds\right)+o_p(n^{1/2}h^2)\\
&=&o_p(1),
\ese
where the last equality is because the integrands have mean zero
conditional on $\bb_0\trans\X$, and the second last equality is obtained
following the same derivation of $\T_3$.
Using these results in (\ref{eq:main}), combined with the results on
(\ref{eq:easy}), it is now clear that the theorem holds. 
\qed

\bibliographystyle{imsart-nameyear}
\bibliography{censorsdr}

\end{document}